%% file: main.tex
\documentclass[a4paper, 10pt]{article}
\usepackage{graphicx}
\usepackage{amsmath, amssymb, amsthm}
\usepackage{color}
\usepackage{booktabs}
\usepackage{xcolor}
\usepackage{bm}
\usepackage{multirow}
\usepackage{cases}
\usepackage[linesnumbered,ruled,vlined]{algorithm2e}
\usepackage{subfig} % it is suggested to use subfloat rather than subfigure
\usepackage{hyperref}
\usepackage[margin=1in]{geometry}
\usepackage{booktabs}
\usepackage{threeparttable}
\usepackage{enumerate}
\usepackage[normalem]{ulem}
\graphicspath{{images/}}

\newtheorem{theorem}{Theorem}
\newtheorem{lemma}[theorem]{Lemma}

\newtheorem{remark}{Remark}

\newcommand\pd[2]{\dfrac{\partial {#1}}{\partial {#2}}}

\newcommand\bx{\boldsymbol{x}}

\newcommand{\imag}{\mathrm{i}}

\newcommand\bbR{\mathbb{R}}

\newcommand\bbZ{\mathbb{Z}}
\newcommand\bbS{\mathbb{S}}

\newcommand\dd{\,\mathrm{d}}

\newcommand\mO{\mathcal{O}}

\newcommand\mF{\mathcal{F}}
\newcommand\mG{\mathcal{G}}

\def\bx{\boldsymbol{x}}
\def\bbR{\mathbb{R}}

\def\bsOmega{\boldsymbol{\Omega}}

\numberwithin{equation}{section}

\definecolor{electricpurple}{rgb}{0.75,0.0,1.0}
\definecolor{darkred}{rgb}{0.65,0,0}
\definecolor{green}{rgb}{0.0, 0.5, 0.0}

\newcommand\II[3]{{{I}_{{#1}, {#2}}^{{#3}}}}

\title {An asymptotic-preserving IMEX $P_N$ method for the gray model of the radiative transfer equation}
 
\author{Jinxue
 Fu\thanks{Institute of Applied Physics and Computational Mathematics,
   Beijing, China, 100088, email: {\tt jinxue.fu@csrc.ac.cn}},
   ~~Juan Cheng\thanks{Institute of Applied Physics and Computational Mathematics, Beijing, China, 100088 and Center for Applied Physics and Technology, and College of Engineering, Peking University, Beijing,  China, 100871, email: \tt{cheng\_juan@iapcm.ac.cn}},
   ~~Weiming Li\thanks{ Institute of Applied Physics and Computational Mathematics, Beijing, China, 100088,  email: \tt{li\_weiming@iapcm.ac.cn}},
   ~~Tao Xiong\thanks{School of Mathematical Sciences, University of Science and Technology of China, Hefei, Anhui, China, 230026, email: \tt{taoxiong@ustc.edu.cn}},
   ~~Yanli Wang\thanks{Beijing Computational Science Research Center, Beijing,
   China, 100193, email: {\tt ylwang@csrc.ac.cn}}.}

\begin{document}
\maketitle
%\tableofcontents
%\clearpage

\input{abs_intro}
\input{model}

\input{method}

\input{analysis}

\input{numerical_results}

\input{article_conclusion}

\bibliographystyle{plain}
\bibliography{reference}
\end{document}

%% file: abs_intro.tex
\begin{abstract}
An asymptotic-preserving (AP) implicit-explicit $P_N$ numerical scheme is proposed for the gray model of the radiative transfer equation, where the first- and second-order numerical schemes are discussed for both the linear and nonlinear models. The AP property of this numerical scheme is proved theoretically and numerically, while the numerical stability of the linear model is verified by Fourier analysis. Several classical benchmark examples are studied to validate the efficiency of this numerical scheme.

\end{abstract}


\section{Introduction}

The radiative transfer equation (RTE) describes the transport of photons and energy exchange with different background materials, including several important processes such as photon transmission, absorption, scattering, and emission \cite{1960Radiative, 1973Pomraning}. RTE has been widely utilized in fields such as astrophysics \cite{Drake2006HighEnergyDensityPF, Lan2022}, biomedical optics \cite{Kim2010}, inertial confinement fusion (ICF) \cite{LANGER2003275, Chen2022}, and solving it efficiently is full of importance for both theoretical analysis and practical applications. However, due to its integro-differential form, high dimensionality which involves time, frequency, spatial and angular variables, and multi-scale features, it is impossible to solve it analytically, and numerically solving it poses significant challenges.

Generally, numerical methods for solving RTE are divided mainly into two categories, the stochastic and the deterministic methods. The most commonly used stochastic method is the implicit Monte Carlo (IMC) method, first proposed by Fleck and Cummings \cite{fleck1971}, which approximates the complex energy exchange process between radiation and matter as a linear system by introducing the Fleck factor. The IMC method is highly parallelizable, adaptable to arbitrary geometries, and performs efficiently in optically thin regions. However, it suffers from statistical noise, and a large number of scattering events are needed to calculate in optically thick regions significantly reducing the computational efficiency. Several significant attempts have been made to increase the computational efficiency of the IMC method, including the transport-diffusion hybrid method \cite{Gentile2001, DENSMORE2007485, Densmore2012}, the moment-based scale-bridging method \cite{dens2015, Park2012} and essentially implicit Monte Carlo (EIMC) method \cite{Shiyi2019, SHI2023112552}. Moreover, the unified gas-kinetic wave particle (UGKWP) methods have been developed, where the non-equilibrium part is simulated by the Monte Carlo method \cite{liweiming2020, LI2024112663, HuYuan2024, Liuchang2023}, and the equilibrium part is approximated using the macroscopic equations, which greatly improve the computational efficiency in the optically thick regions.

For deterministic methods, one of the popular methods is the discrete ordinate ($S_N$) method \cite{koch1995, Lath1965}, which directly discretizes the angular variable using a quadrature rule. The $S_N$ method is known to suffer from ray effects \cite{Math1999, Morel2003}. To mitigate ray effects, several methods have been proposed, such as the rotated $S_N$ method \cite{CAMMINADY2019105}, artificial scattering $S_N$ method \cite{Frank2020Rayeffectartifical2020} and random ordinate method \cite{TangMin2024}. Moreover, the unified gas-kinetic method (UGKS) is proposed in the framework of $S_N$ for RTE, which handles the smooth transition from microscopic transport to macroscopic diffusion by constructing interface fluxes using the integral solution of the transport equation along with collision terms \cite{sun2015asymptotic1, sun2018asymptotic}. Besides, the spatial second-order positive method is proposed in \cite{Xuxiaojing2023} and a uniformly unconditionally stable scheme is proposed in \cite{Xiong2023}, both in the framework of UGKS.

Another commonly used deterministic method is the spherical harmonic ($P_N$) method \cite{Kershaw1976, Lewis1993}, which employs spherical harmonic expansion to approximate the distribution function effectively. The $P_N$ method exhibits rotational invariance, thereby mitigating the occurrence of ray effects. However, since the $P_N$ method is essentially a truncated spectral method, it can lead to non-physical oscillations and even negative energy density \cite{McCl2008, McCl2010}. To solve this problem, several approaches have been attempted, such as adding artificial scattering terms \cite{Olson2009}, filtering methods \cite{implicit2016Laboure, mcclarren2010robust, Xuxiaojing2024}, using different closure methods \cite{ABDELMALIK2023116454, MINERBO1978541, ZHENG2016682} and hybrid methods \cite{McClarren2021, LiQi2024}.

In practical simulations, the optical parameters of RTE can vary by several orders of magnitude across different regions, resulting in multiscale challenges. To address this issue, the asymptotic-preserving (AP) schemes have been proposed \cite{Jin1993, Jin1991, Larsen1989, Larsen1987}. An AP scheme converges to a stable and consistent discretization of the diffusion limit equation as the mean free path tends to zero while keeping the time step and mesh size fixed \cite{Larsen1989, Larsen1987}. In recent decades, the design of AP schemes has become a popular research topic and several AP schemes have been developed. The micro-macro decomposition, by explicitly separating the multiscale characteristics of a system into macroscopic and microscopic components, has become a primary approach for constructing AP schemes \cite{Lemou2010}. An AP semi-Lagrangian discontinuous Galerkin (DG)  method is proposed in \cite{CAI2024113190}, and the implicit-explicit (IMEX) temporal discretization with spatially DG discretization with AP property is proposed in \cite{xiong2020, xiong2015, peng2020}. The three-state update method is utilized to capture the correct front propagation in the diffusion limit in \cite{tang2021}, while the multiscale high/low order (HOLO) method is adopted to build the AP scheme in \cite{Hammer2019, Maginot2016} for the frequency-dependent RTE, and the AP scheme for frequency-dependent RTE is also proposed in  \cite{ZHANG2023112368}. Besides, through the Chapman-Enskog expansion, the order of the radiation intensity moment coefficients $I_l$ in the $P_N$ system with respect to $\epsilon$ was obtained in \cite{IMEX2022}, specifically $I_l = \mathcal{O}(\epsilon^l)$. Based on this, an AP IMEX numerical scheme in the framework of $P_N$ is discussed in \cite{IMEX2022}, where the higher-order expansion coefficients with respect to $\epsilon$ are treated explicitly, and the lower-order expansion coefficients are treated implicitly. This IMEX numerical scheme is proposed to solve the $P_N$ system with the numerical cost of an explicit scheme. However, when the multiscale parameter $\epsilon$ approaches zero, this scheme converges to an explicit form of the diffusion limit equation, which results in a constraint of the time-step length of $\Delta x^2$ when $\epsilon$ is small, thus affecting computational efficiency greatly.

In this paper, we propose a new AP IMEX numerical method, which is an enhancement of the numerical scheme in \cite{IMEX2022}, especially in the diffusive regime. To improve the computational efficiency in the diffusion regime, the convection term in the first governing equation of the $P_N$ system is treated implicitly, while it is explicitly treated in \cite{IMEX2022}. With this improvement,  this new AP scheme converges to an implicit scheme of the diffusion limit equation when $\epsilon$ reaches zero, which greatly enlarges the time step size in the diffusive regime. Then, the mid-point scheme for the explicit terms and the Crank-Nicolson scheme for the implicit terms are utilized to derive the second-order IMEX scheme. Furthermore, the stability properties of the new scheme for linear RTE are studied using the Fourier analysis method, and the AP property of this IMEX scheme is proved theoretically. A series of benchmark examples are studied to test the effectiveness of this scheme, including the AP test, the plane source problem, the Marshak wave problems, the line source problem, the lattice problem, and the Riemann problem. Finally, the efficiency of this new AP IMEX method is validated by comparing with the IMEX method proposed in \cite{IMEX2022}. 

The rest of this paper is organized as follows. In Sec. \ref{sec:pre}, the RTE model and the $P_N$ method are introduced. The temporal discretization for the linear RTE and nonlinear RTE is discussed in Sec. \ref{sec:linear} and \ref{sec:RTE}, respectively. In Sec. \ref{sec:spatial_AP}, the spatial discretization together with the proof of the AP property and the analysis of the numerical stability by Fourier analysis is presented. The numerical experiments together and the study of efficiency are presented in Sec. \ref{sec:num} with some conclusions in Sec. \ref{sec:conclusion}. 

%% file: model.tex
% vim: wrap:tw=70:spell:colorcolumn=101

\section{Gray model of RTE and the \texorpdfstring{$P_N$}{PN} method}
\label{sec:pre}
In this section, the gray model of the radiative transfer equations will be introduced with its several basic properties. Additionally, the $P_N$-based method to discretize the angular space will be discussed.

\subsection{The gray approximation to RTE}
\label{sec:rte}
The radiative transfer equations describe the transport of radiation and the energy exchange with the material. Under the assumption of isotropic scattering and local thermal equilibrium, the gray approximation to the radiative transfer equations in the absence of scattering and external sources has the form as 
\begin{subequations}
  \label{eq:RTE}
  \begin{align}
    \label{eq:I_rt}
    & \dfrac{\epsilon^2}{c} \pd{I}{t} + \epsilon \bsOmega \cdot \nabla I
    = \sigma \left(\frac{1}{4\pi}acT^4 - I\right), \\
    \label{eq:T}
   &  \epsilon^2 C_{v} \pd{T}{t} \equiv \epsilon^2 \pd{U}{t} =
    \sigma\left(\int_{\bbS^2}  I  \dd \bsOmega - ac T^4\right).
      \end{align}
\end{subequations}
Here, $I(t, \bx, \bsOmega)$ is the radiation density with $\bsOmega \in \bbS^2$ the surface of the unit sphere. $\bx = (x_1, x_2, x_3)\in \bbR^3$ is the spatial variable. $T(t, \bx)$ is the material temperature and $\sigma$ is the opacity of the material, which may depend on the material temperature $T$ and vary greatly with $T$. $\epsilon$ is a dimensionless parameter, which is similar to the Knudsen number in the Boltzmann equations \cite{Mie2013, Larsen1987}. $c$ is the speed of light, and $a$ is the radiation constant given by 
\begin{equation}
    \label{eq:coe_a}
    a = \frac{8 \pi k^4}{15 h^3 c^3},
\end{equation}
with $h$ the Planck constant and $k$ the Boltzmann constant. In the governing equation \eqref{eq:T} of the material temperature $T$, the material energy density  $U$ is related to the material temperature $T$ through the equation of state as 
\begin{equation}
    \label{eq:eos}
    \pd{U}{T} = C_v,
\end{equation}
where $C_v(\bx, T)$ is the specific heat capacity of the material, which may depend on the spatial variable $\bx$ and the material temperature $T$. The radiation energy density $E_r$ \cite{dens2015} is defined as 
\begin{equation}
    \label{eq:energy_Er}
    E_r = \frac{1}{c}\int_{\bsOmega \in \bbS^2} I(t, \bx, \bsOmega) \dd \bsOmega, 
\end{equation}
which is the zeroth angular moment of the radiation intensity, and the radiation temperature $T_r$  is given by the following relation
\begin{equation}
    \label{eq:radiation_tem}
    E_r = aT_r^4. 
\end{equation}

In the RTE system, \eqref{eq:I_rt} describes the transport of radiation through background media, \eqref{eq:T} governs the change in material energy from radiation being absorbed and emitted by the material. By integrating the radiation equation \eqref{eq:I_rt} over all directions, and adding the equation of $T$ \eqref{eq:T} together, we can obtain the following total energy conservation equation 
\begin{equation}
    \label{eq:energy}
    \pd{}{t} (E_r + U) + \frac{1}{\epsilon} \nabla \cdot F_r = 0, 
\end{equation}
where the radiation flux $F_r$ is defined as 
\begin{equation}
    \label{eq:flux}
    F_r = \int_{\bbS^2} \bsOmega I \dd \bsOmega.
\end{equation}

When the evolution of the transport and interaction of the radiation is in the same medium, \eqref{eq:RTE} is reduced to a linear transport equation \cite{Mie2013,seibold2014starmap,Hu2021lowrank} as 
\begin{equation}
  \label{eq:linear_RTE}
  \dfrac{\epsilon^2}{c} \pd{I}{t} + \epsilon \bsOmega \cdot \nabla I
  = \sigma_{s}\left(\frac{1}{4\pi} \int_{\bbS^2} I \dd \bsOmega - I\right)- \epsilon^2 \sigma_{a}I +\epsilon^2 \frac{G}{4\pi},
\end{equation}
where $\sigma_s(\bx)$ is the scattering coefficient, $\sigma_a(\bx)$ is the absorption coefficient, and $G(\bx)$ is a given source term.

For the RTE system \eqref{eq:RTE} and \eqref{eq:linear_RTE}, the parameter $\epsilon$ can be very small in some zones and very large elsewhere. When $\epsilon$ approaches zero, the diffusion limit of the RTE system is derived. Precisely, the radiation intensity $I$ in \eqref{eq:RTE} will approach a Planckian distribution at the local temperature \cite{sun2018asymptotic, sun2015asymptotic1} as 
\begin{equation}
    \label{eq:I0}
     I^{(0)} = \frac{1}{4\pi}ac\left(T^{(0)}\right)^4,
\end{equation}
and the diffusion limit of the RTE system has the below form, 
\begin{equation}
  \label{eq:limit}
 C_v \pd{T^{(0)}}{t} + a \pd{}{t}\left(T^{(0)}\right)^4 = \nabla \cdot \left(\frac{a
    c}{3 \sigma} \nabla \left(T^{(0)}\right)^4\right).
\end{equation}
Additionally, when $\epsilon$ approaches zero, the radiation intensity $I$ in \eqref{eq:linear_RTE} reduces to 
\begin{equation}
    \label{eq:limit_rte_I}
    I^{(0)} = \frac{1}{4 \pi} \int_{\bbS^2} I \dd \bsOmega, 
\end{equation}
and 
the diffusion limit of \eqref{eq:linear_RTE} becomes 
\begin{equation}
  \label{eq:linear_limit}
   \frac{1}{c}\pd{I^{(0)}}{t} - \nabla \cdot\left( \frac{1}{3 \sigma_s} \nabla I^{(0)} \right) = -\sigma_a I^{(0)} + \frac{G}{4 \pi}. 
\end{equation}

The RTE system forms a nonlinear system of equations with the radiation density $I$ and the material temperature $T$, which is also tightly coupled. Solving \eqref{eq:RTE} is difficult due to the high dimensionality of the angular space and the stiff coupling between the radiation and material energy. Moreover, due to the multi-scale property led by the variation of the opacity $\sigma$, the standard numerical methods to solve it are quite difficult, for they should resolve the smallest microscopic scale, and an asymptotic preserving (AP) scheme is demanded. In this work, we will propose an AP IMEX  numerical scheme based on \cite{IMEX2022}, and the $P_N$ method is first introduced in the next section. 

\subsection{\texorpdfstring{$P_N$}{PN} method to discretize the angular space}
\label{sec:Pn}
For the $P_N$ method, the series of basis functions is utilized to approximate the radiation intensity $I$, and we refer to \cite{Mie2013, semi2008Ryan, McCl2008} and the reference therein for more details. For simplicity, we will introduce the $P_N$ method for the spatially 1D RTE system. The governing equations of RTE \eqref{eq:RTE} reduce to the following form in the 1D case as 
\begin{subequations}
    \label{eq:1D_RTE}
    \begin{align}
        \label{eq:1D_RTE_I}
         & \frac{\epsilon^2}{c} \pd{I}{t} +\epsilon \mu \pd{I}{x} =
    \sigma\left(\frac{1}{2}a c   T^4 - I\right),   \\
    \label{eq:1D_RTE_T}
    & \epsilon^2 C_{v} \pd{T}{t} = \sigma \left(\int_{-1}^1 I \dd \mu - ac
      T^4\right),
    \end{align}
\end{subequations}
where $\mu = \cos \theta\in [-1, 1]$ is the internal coordinate associated with the angular $\theta \in [0, \pi]$. The other parameters remain the same as in \eqref{eq:RTE}. The diffusion limit \eqref{eq:limit} is changed into
\begin{equation}
    \label{eq:1D_limit}
    C_v \pd{T^{(0)}}{t} + a \pd{}{t} \left(T^{(0)} \right)^4 = \pd{}{x} \left(\frac{ac}{3 \sigma} \pd{}{x}\left( T^{(0)}\right)^4\right). 
\end{equation}
Additionally, the 1D form of \eqref{eq:linear_RTE} without the source term reduces to 
\begin{equation}
    \label{eq:1D_linear_RTE}
    \frac{\epsilon^2}{c}\pd{I}{t} + \epsilon \mu \pd{I}{x} = \sigma_s \left(\frac{1}{2}\int_{-1}^{1} I \dd \mu - I \right) -\epsilon^2 \sigma_a I,
\end{equation}
with its 1D diffusion limit \eqref{eq:linear_limit} reducing to 
\begin{equation}
    \label{eq:1D_linear_limit}
    \frac{1}{c}\pd{I^{(0)}}{t} -\pd{}{x} \left( \frac{1}{3 \sigma_s} \pd{}{x} I^{(0)} \right) = -\sigma_a I^{(0)}. 
\end{equation}

For the 1D RTE system, the Legendre polynomials are utilized as the basis functions, and the radiation intensity $I$ is approximated as 
\begin{equation}
    \label{eq:approx}
    I(t, x, \mu) \approx \sum_{l = 0}^M \frac{2l +1}{2} I_l(t, x) P_l(\mu), 
\end{equation}
where $P_l(\mu)$ is the $l$-th Legendre polynomial, and $M$ is the expansion truncation order. With the orthogonality of the basis function 
\begin{equation}
\label{eq:orth}
     \int_{-1}^{1}P_{m}(\mu)P_{n}(\mu)\dd\mu = \frac{2\delta_{mn}}{2m + 1},
\end{equation}
the expansion coefficient $I_l$ is calculated as 
\begin{equation}
    \label{eq:coe}
    I_l(t, x) = \int_{-1}^1 P_l(\mu) I(t, x, \mu) \dd \mu. 
\end{equation}
Multiplying $P_l(\mu) $ on both sides of \eqref{eq:1D_RTE_I}, and integrating over $\mu$, with the recursion property of the Legendre polynomials 
 \begin{equation}
 \label{eq:rec}
   \mu P_{n}(\mu) =  \frac{n + 1}{2n + 1}P_{n+1}(\mu) + \frac{n}{2n + 1}P_{n - 1}(\mu),\qquad P_{0}(\mu) = 1, \qquad P_{1}(\mu) = \mu,
 \end{equation}
we obtain the $P_N$ equations as 
\begin{subnumcases}
    {\label{eq:PN_RTE}}
    \label{eq:PN_0}
     \frac{\epsilon^2}{c} \pd{I_0}{t} +
    \epsilon \pd{I_{1}}{x} = \sigma(acT^{4} - I_{0}),\\
    \label{eq:PN_l}
     \frac{\epsilon^2}{c} \pd{I_{l}}{t} + \epsilon a_{l-1} 
    \pd{I_{l-1}}{x} + \epsilon b_{l+1}\pd{I_{l+1}}{x} = -\sigma
    I_{l},\quad l = 1,\cdots,M-1, \\
    \label{eq:PN_M}
     \frac{\epsilon^2}{c} \pd{I_{M}}{t} + \epsilon a_{M-1} 
    \pd{I_{M-1}}{x} = -\sigma
    I_{M},
\end{subnumcases}
with 
\begin{equation}
  \label{eq:coe_moment1D}
 a_{l-1} = \frac{l}{2l + 1},\quad b_{l+1} = \frac{l + 1}{2l + 1}, \qquad l = 1, \cdots M.
\end{equation}
In the framework of $P_N$ method, \eqref{eq:1D_RTE_T} reduces to 
\begin{equation}
    \label{eq:PN_T}
     \epsilon^2 C_{v} \pd{T}{t} = \sigma \left(I_0- ac
      T^4\right). 
\end{equation}
 Together with \eqref{eq:PN_RTE}, \eqref{eq:PN_T} and the closure relation 
  \begin{equation}
     \label{eq:closure}
     I_{M+1} = 0,
 \end{equation}
 we obtain the $P_N$ system of the RTE \eqref{eq:1D_RTE}. The radiation energy density $E_r$ \eqref{eq:energy_Er} and the radiation flux $F_r$ \eqref{eq:flux} are related to the expansion coefficients $I_l$ as 
 \begin{equation}
     \label{eq:E_r_Il}
     E_r = \frac{I_0}{c}, \qquad F_r = I_1. 
 \end{equation}
Additionally, the $P_N$ system of the 1D linear RTE \eqref{eq:1D_linear_RTE} has the form below
\begin{subnumcases}  
 { \label{eq:final_Pn_linear}}
  \label{eq:Pn_linear_0}
\frac{\epsilon^2}{c}\pd{I_0}{t} + \epsilon \pd{I_1}{x} = - \epsilon^2 \sigma_a I_0,\\
\label{eq:Pn_linear_1}
\frac{\epsilon^2}{c}\pd{I_1}{t} + \epsilon a_{0} \pd{I_{0}}{x} +  \epsilon b_{2} \pd{I_{2}}{x}= - (\epsilon^2 \sigma_a  + \sigma_s)I_{1}, \\
\label{eq:Pn_linear_l}
\frac{\epsilon^2}{c}\pd{I_l}{t} + \epsilon a_{l-1} \pd{I_{l-1}}{x} +  \epsilon b_{l+1} \pd{I_{l+1}}{x}= - (\epsilon^2 \sigma_a  + \sigma_s)I_{l},\qquad l = 2,\cdots,M - 1, \\
\label{eq:Pn_linear_M}
\frac{\epsilon^2}{c}\pd{I_M}{t} + \epsilon a_{M-1} \pd{I_{M-1}}{x} = - (\epsilon^2 \sigma_a  + \sigma_s)I_{M},
\end{subnumcases}
where the coefficients $a_l$ and $b_l$ are the same as \eqref{eq:coe_moment1D}.

For high-dimensional RTE, the spherical harmonics are utilized to obtain the $P_N$ system, and we refer to \cite{McCl2010, mcclarren2010robust} and the references therein for more details. To solve the $P_N$ system of the RTE \eqref{eq:PN_RTE} and \eqref{eq:PN_T}, as well as the $P_N$ system of the linear RTE \eqref{eq:final_Pn_linear}, an AP numerical scheme is proposed here, which can resolve the diffusion limit when $\epsilon$ approaches zero with the same time step length, which we will introduce in detail in the following sections.

%% file: method.tex
\section{Temporal discretization for the linear RTE}
\label{sec:linear}
In this section, the temporal discretization for the linear radiation transfer equations will be introduced. 
We want to emphasize that in \cite{IMEX2022}, a similar AP IMEX numerical scheme is proposed based on a similar order analysis. However, in the numerical scheme therein, it will converge to an explicit five-point difference scheme of the diffusion limit \eqref{eq:1D_limit}, which restricts the time step length to $\mathcal{O}(\Delta x^2)$. In the following sections, we will propose an AP temporal scheme that can resolve the diffusion limit with a consistent time-step length independent of $\epsilon$ and the time-step length is constrained to $\mathcal{O}(\Delta x)$. Without loss of generality, we will begin with the 1D linear RTE \eqref{eq:1D_linear_RTE} with constant scattering and absorption coefficients. The first- and second-order semi-discrete numerical schemes are introduced in Sec. \ref{sec:linear_first} and \ref{sec:linear_second}, respectively. 

\subsection{First-order temporal discretization}
\label{sec:linear_first}
In this section, the first-order temporal discretization for the 1D linear RTE \eqref{eq:1D_linear_RTE} is proposed. To design an AP temporal scheme that can resolve the diffusion limit with a time step length constrained to $\mO(\Delta x)$, it should satisfy the fact that the numerical scheme converges to an implicit numerical scheme of the diffusion limit. 

We will begin with the order analysis of the expansion coefficients $I_l, l = 0, 1, \cdots, M$, based on which the temporal discretization will be proposed. Following the analysis in \cite{IMEX2022}, the order analysis is similar to the Chapman-Enskog method \cite{Chapman1990}, and is adopted on the parameter $\epsilon$. The main result is listed in Lem. \ref{thm:CE} and we refer to \cite{IMEX2022} for the detailed proof. 

\begin{lemma}
\label{thm:CE}
For the RTE \eqref{eq:1D_RTE} and the linear system \eqref{eq:1D_linear_RTE}, when the parameter $\epsilon$ approaches zero, the expansion coefficient $I_l, l = 0, \cdots, M$ in \eqref{eq:approx} has the property below 
\begin{equation}
    \label{eq:CE}
    I_l = \mathcal{O}(\epsilon^l), \qquad l = 0, \cdots, M.
\end{equation}
\end{lemma}
% \begin{proof}[Proof of Lemma \ref{thm:CE}]
%     When parameter $\epsilon$ is small, the radiation intensity $I$ is expanded into the power series of $\epsilon$ by the Chapman-Enskog method as 
%     \begin{equation}
%     \label{eq:CE_I}
%     I = I^{(0)} + \epsilon I^{(1)} + \epsilon^2 I^{(2)} +
%     \cdots, \qquad I^{(l)} = \mathcal{O}(1). 
%   \end{equation}
% Define  
%   \begin{equation}
%   \label{eq:I_lk}
%     I_{l}^{(k)} =
%     \int_{-1}^{1} P_{l} I^{(k)}(
%     t, x, \mu) \dd \mu.
%   \end{equation}
%   Following the mathematical induction in \cite{IMEX2022}, it is proved that 
%   \begin{equation}
%       \label{eq:I_lk1}
%        I^{(k)}_l = 0, \qquad  {\rm for}\quad l > k. 
%   \end{equation}
%   Therefore, substituting \eqref{eq:CE_I} into \eqref{eq:coe}, with \eqref{eq:I_lk} and \eqref{eq:I_lk1},   
%   it holds for $I_l$ that 
%   \begin{equation}
%       \label{eq:I_lk2}
%       I_l(t, x) =  \int_{-1}^{1} P_{l}(\mu) I(
%     t, x, \mu) \dd \mu =  \int_{-1}^{1} P_{l}(\mu)  \left(\epsilon^l I^{(l)}  + \epsilon^{l+1} I^{(l+1)} + \cdots \right) \dd \mu  = \mathcal{O}(\epsilon^l).
%   \end{equation}
%   Then, we complete the proof. 
% \end{proof}

To capture the diffusion limit, only the temporal discretization of $I_0$ and $I_1$ should be treated specially. We begin with the order analysis of the $P_N$ system \eqref{eq:Pn_linear_0} and \eqref{eq:Pn_linear_1}.  Based on the order analysis in Lem. \ref{thm:CE}, it holds that 
\begin{align}
\label{eq:linear_rte_0}
& \epsilon^2 \pd{I_0}{t} = \mO(\epsilon^2), \qquad   \epsilon \pd{I_{1}}{x} = \mO(\epsilon^2), \qquad \epsilon^2 \sigma_a  I_0 = \mO(\epsilon^2), \\
\label{eq:linear_rte_1}
& \epsilon^2 \pd{I_1}{t} = \mO(\epsilon^3), \qquad   \epsilon \pd{I_{0}}{x} = \mO(\epsilon), \qquad \epsilon  b_{2} \pd{I_{2}}{x} = \mO(\epsilon^{3}), \qquad (\epsilon^2 \sigma_a +\sigma_s) I_1 = \mO(\epsilon). 
\end{align}
Following the routine in \cite{IMEX2022}, the terms at the high-order of $\epsilon$ are treated explicitly, while other terms are treated implicitly. Therefore, the detailed discretization for $I_0$ and $I_1$ is 
\begin{itemize}
    \item $I_0$: the convection and right-hand side terms in \eqref{eq:Pn_linear_0} are treated implicitly since they are all at the order of $\mO(\epsilon^2)$.
    \item $I_1$: the convection term $\pd{I_{0}}{x}$ and the right-hand side terms in \eqref{eq:Pn_linear_1} are treated implicitly, since they are at the order of $\mO(\epsilon)$, while the convection term $\pd{I_2}{x}$ is treated explicitly since it is at the order of $\mO(\epsilon^3)$. 
\end{itemize}
Precisely, supposing the numerical solution at $n$-th time level is $I_l^n, l = 0, \cdots, M$, with the forward Euler scheme utilized for the time derivatives, the first-order semi-discrete numerical scheme for $I_0$ and $I_1$ has the form below 
{\small 
\begin{subnumcases}
 {   \label{eq:first_linear}}
    \label{eq:first_linear_I0}
        \frac{\epsilon^2}{c} \frac{I_{0}^{n+1} - I_{0}^n }{\Delta t}
      +   \epsilon\left(\pd{I_{1}}{x}\right)^{n+1} =-\epsilon^2 \sigma_{a} I_{0}^{n+1}, \\
      \label{eq:first_linear_I1}
     \frac{\epsilon^2}{c} \frac{I_{1}^{n+1} - I_{1}^n }{\Delta t} +  \epsilon a_0 \left(\pd{I_{0}}{x}\right)^{n+1} + \epsilon b_2\left(\pd{I_{2}}{x}\right)^{n} = -\left(\epsilon^2\sigma_a + \sigma_s\right)I_{1}^{n+1}.
\end{subnumcases}
}Moreover, when $\epsilon$ approaches zero, \eqref{eq:first_linear_I0} will degenerate to zero and \eqref{eq:first_linear_I1} will degenerate to 
\begin{equation}
    \label{eq:deg_I1}
    \sigma_s I_1 = 0,
\end{equation}
in which case, we can not resolve the diffusion limit. To avoid this, we introduce a new variable as 
\begin{equation}
    \label{eq:hat_I1}
    \hat{I}_1 = \frac{I_1}{\epsilon} = \mO(1). 
\end{equation}
Substituting \eqref{eq:hat_I1} into \eqref{eq:first_linear}, it reduces to 
{\small 
\begin{subnumcases}
  {  \label{eq:first_linear_I}}
       \label{eq:first_linear_I0_1}
        \frac{1}{c} \frac{I_{0}^{n+1} - I_{0}^n }{\Delta t}
      +   \left(\pd{\hat{I}_{1}}{x}\right)^{n+1} =- \sigma_{a} I_{0}^{n+1}, \\
      \label{eq:first_linear_I1_1}
       \frac{\epsilon^2}{c} \frac{\hat{I}_{1}^{n+1} - \hat{I}_{1}^n }{\Delta t} +   a_0 \left(\pd{I_{0}}{x}\right)^{n+1} +  b_2\left(\pd{I_{2}}{x}\right)^{n} = -\left(\epsilon^2\sigma_a + \sigma_s\right)\hat{I}_{1}^{n+1}.
  \end{subnumcases}
}For the expansion coefficients $I_l, l >1$, based on the order analysis, it holds for \eqref{eq:Pn_linear_l} and \eqref{eq:Pn_linear_M} that 
\begin{align}
\label{eq:linear_rte_l}
\epsilon  a_{l-1} \pd{I_{l-1}}{x} = \mO(\epsilon^l), \qquad \epsilon  b_{l+1} \pd{I_{l+1}}{x} = \mO(\epsilon^{l+2}), \qquad (\epsilon^2 \sigma_a +\sigma_s) I_l = \mO(\epsilon^{l}), \qquad l = 2, \cdots, M,
\end{align}
with $I_{M+1} = 0$. 
Following the same routine that the terms at the high-order of $\epsilon$ are treated explicitly, while other terms are treated implicitly, the first-order semi-discrete numerical scheme for $I_l, l > 1$ has the following form
{\small 
\begin{subnumcases}
  {\label{eq:first_linear_II}}
   \label{eq:first_linear_I2}
      \frac{\epsilon^2}{c} \frac{I_{2}^{n+1} - I_{2}^n
      }{\Delta t} + \epsilon a_{1}\left(\pd{I_{1}}{x}\right)^{n+1} +
      \epsilon b_{3}\left(\pd{I_{3}}{x}\right)^{n}=
      -\left( \epsilon^{2} \sigma_a +  \sigma_s\right) I_2^{n+1}, \\ 
         \label{eq:first_linear_Il}
       \frac{\epsilon^2}{c} \frac{I_{l}^{n+1} - I_{l}^n
      }{\Delta t} + \epsilon a_{l-1}\left(\pd{I_{l-1}}{x}\right)^{n+1} +
      \epsilon b_{l+1}\left(\pd{I_{l+1}}{x}\right)^{n}=
      -\left( \epsilon^{2} \sigma_a +  \sigma_s\right) I_l^{n+1}, \qquad l = 3, \cdots, M-1, \\
         \label{eq:first_linear_IM}
       \frac{\epsilon^2}{c} \frac{I_{M}^{n+1} - I_{M}^n
      }{\Delta t} + \epsilon a_{M-1}\left(\pd{I_{M-1}}{x}\right)^{n+1} 
      =      -\left( \epsilon^{2} \sigma_a +  \sigma_s\right) I_M^{n+1}.
\end{subnumcases}
}Together with \eqref{eq:first_linear_I} and \eqref{eq:first_linear_II}, we obtain the first-order temporal discretization for the linear RTE \eqref{eq:1D_linear_RTE}. For the first part \eqref{eq:first_linear_I}, the linear equation system of $I_0$ and $I_1$ is solved, while for the second part \eqref{eq:first_linear_II}, the equations for $I_l$ are solved successfully. In this case, when updating the terms $I_l$, the implicit terms $\left(\pd{I_{l-1}}{x}\right)^{n+1}$ are already obtained. Thus, the IMEX scheme \eqref{eq:first_linear_II} can be solved at the numerical cost of an explicit scheme. Moreover, we will verify this in Sec. \ref{sec:AP} that this scheme will converge to an implicit scheme of the diffusion limit \eqref{eq:1D_linear_limit}. 
Additionally, this first-order semi-discrete scheme can be extended to a second-order scheme, which we will introduce in detail in Sec. \ref{sec:linear_second}.

\subsection{Second-order temporal discretization}
\label{sec:linear_second}
To obtain the second-order temporal discretization, the second-order IMEX numerical scheme \cite{xiong2015} is utilized. Precisely, we employ the mid-point scheme for the explicit terms, while the Crank-Nicolson scheme is adopted for the implicit terms. Once, we have obtained the numerical solution $I_l^n$ at the time level $n$, the numerical solution $I_l^{n+1/2}$ at $(n+1/2)$-th time level is obtained with the first-order scheme \eqref{eq:first_linear_I} and \eqref{eq:first_linear_II} for time-step length $\Delta t/2$. Thus, the second-order temporal discretization for \eqref{eq:1D_linear_RTE} has the form as  
{\small
  \begin{align}
   \label{eq:first_order_scheme_linear_time_2_order}
   \left\{
   \begin{aligned}     
      & \frac{1}{c} \frac{I_{0}^{n+1} - I_{0}^n }{\Delta t}
      +   \frac{1}{2}\left[\left(\pd{\hat{I}_{1}}{x}\right)^{n+1} + \left(\pd{\hat{I}_{1}}{x}\right)^{n}\right] =-\frac{\sigma_{a}}{2} \left(I_{0}^{n+1} + I_{0}^{n}\right), \\
      & \frac{\epsilon^2}{c} \frac{\hat{I}_{1}^{n+1} - \hat{I}_{1}^n }{\Delta t} + \frac{a_0}{2} \left[\left(\pd{I_{0}}{x}\right)^{n+1} + \left(\pd{I_{0}}{x}\right)^{n}\right] +  b_2 \left(\pd{I_{2}}{x}\right)^{n + \frac{1}{2}} = -\frac{\epsilon^2\sigma_a + \sigma_s}{2}\left(\hat{I}_{1}^{n+1} + \hat{I}_{1}^{n}\right),\\
    & \frac{\epsilon^2}{c} \frac{I_{l}^{n+1} - I_{l}^n
      }{\Delta t} +\frac{ \epsilon a_{l-1}}{2}\left[\left(\pd{I_{l-1}}{x}\right)^{n+1} + \left(\pd{I_{l-1}}{x}\right)^{n}\right]
      + \epsilon b_{l+1}\left(\pd{I_{l+1}}{x}\right)^{n + \frac{1}{2}}=
      -\frac{ \epsilon^{2} \sigma_a +  \sigma_s}{2} \left(I_l^{n+1} + I_l^{n}\right), \qquad l = 2, \cdots, M-1,\\
      & \frac{\epsilon^2}{c} \frac{I_{M}^{n+1} - I_{M}^n
      }{\Delta t} +\frac{ \epsilon a_{M-1}}{2}\left[\left(\pd{I_{M-1}}{x}\right)^{n+1} + \left(\pd{I_{M-1}}{x}\right)^{n}\right]
      =
      -\frac{ \epsilon^{2} \sigma_a +  \sigma_s}{2} \left(I_M^{n+1} + I_M^{n}\right).
      \end{aligned} \right.
    \end{align}
}For now, the second-order semi-discrete system for the linear RTE \eqref{eq:1D_linear_RTE} is derived, which is a second-order IMEX scheme, and can be extended to the high-order IMEX scheme naturally, and we refer \cite{xiong2015,IMEX2022} for more details.

\section{Temporal discretization for the gray model of RTE}
\label{sec:RTE}
In this section, an AP scheme is proposed for the temporal discretization of the gray model of RTE \eqref{eq:1D_RTE}. Compared to linear RTE \eqref{eq:1D_linear_RTE}, a special design is made due to the non-linear interaction of the radiation and the background material. We first assume that the opacity is a constant in Sec. \ref{sec:RTE_con_sig}, and then a discussion on the design of the AP scheme for the nonlinear opacity is presented in Sec. \ref{sec:RTE_Marshak}. 

\subsection{RTE with constant opacity}
\label{sec:RTE_con_sig}
For the gray model of RTE, the non-linear interaction of the radiation and the background material makes it more difficult to design the AP numerical scheme. We first assume that the opacity $\sigma = \mO(1)$ is a constant. We begin with the order analysis of the expansion coefficients. With the results in Lem. \ref{thm:CE}, it holds for the $P_N$ system \eqref{eq:PN_RTE} of RTE \eqref{eq:1D_RTE}
\begin{subequations}
    \label{eq:order_PN_RTE}
    \begin{align}
   &  \epsilon^2 \pd{I_0}{t} = \mO(\epsilon^2), \qquad  \epsilon \pd{I_1}{x} = \mO(\epsilon^2), \qquad   \sigma(ac T^4 - I_0)  = \mO(\epsilon^2), \\
   &  \epsilon^2 \pd{I_1}{t} = \mO(\epsilon^3), \qquad \epsilon a_0\pd{I_0}{x} = \mO(\epsilon), \qquad \epsilon b_2 \pd{I_2}{x} = \mO(\epsilon^3), \qquad \sigma I_1 = \mO(\epsilon), \\
   &  \epsilon^2 \pd{I_l}{t} = \mO(\epsilon^{l+2}), \qquad \epsilon a_{l-1}\pd{I_{l-1}}{x} = \mO(\epsilon^{l}), \qquad \epsilon b_{l+1} \pd{I_{l+1}}{x} = \mO(\epsilon^{l+2}), \qquad \sigma I_l = \mO(\epsilon^l). 
    \end{align}
\end{subequations}
With the same principle that the terms at the high order of $\epsilon$ are treated explicitly, while those at the low order of $\epsilon$ are treated implicitly, the temporal discretization of \eqref{eq:1D_RTE} has the form below 
{\small 
\begin{subequations}
    \label{eq:first_RTE}
    \begin{numcases}
        {\label{eq:first_RTE_I0}}
       \frac{\epsilon^2}{c} \frac{I_{0}^{n+1} - I_{0}^n }{\Delta t}
      +  \epsilon^2\left(\pd{\hat{I}_{1}}{x}\right)^{n+1} =\sigma \left[ac (T^{n+1})^4 - I_{0}^{n+1} \right],\\ 
        {\label{eq:first_RTE_I1}}
       \frac{\epsilon^2}{c} \frac{\hat{I}_{1}^{n+1} - \hat{I}_{1}^n }{\Delta t} + a_0\left(\pd{I_{0}}{x}\right)^{n+1} + b_2\left(\pd{I_{2}}{x}\right)^{n} = - \sigma\hat{I}_{1}^{n+1},\\
         {\label{eq:first_RTE_T}}
       \epsilon^2 C_{v} \frac{T^{n+1} - T^{n}}{\Delta t}  =
    \sigma\left[ I_{0}^{n+1} - ac (T^{n+1})^4\right],
  \end{numcases}
  \begin{equation}
   \label{eq:first_RTE_Il}
    \quad  \frac{\epsilon^2}{c} \frac{I_{l}^{n+1} - I_{l}^n
      }{\Delta t} + \epsilon a_{l-1}\left(\pd{I_{l-1}}{x}\right)^{n+1} +
      \epsilon b_{l+1}\left(\pd{I_{l+1}}{x}\right)^{n}=
      -\sigma I_l^{n+1}, \qquad l = 2, \cdots, M.
  \end{equation}
\end{subequations}
}Here, $I_l^n, l = 0, \cdots, M$ are the numerical solutions at the $n$-th time level, while $\hat{I}_1^n$ is defined the same as in \eqref{eq:hat_I1}, with $I_{M+1} = 0$.  In \eqref{eq:first_RTE}, the same strategy as in \eqref{eq:hat_I1} is applied to $I_1$ to avoid the reduction of the scheme when $\epsilon$ approaches zero. 

Unlike the $P_N$ system of the linear RTE \eqref{eq:first_linear}, the system \eqref{eq:first_RTE} results in a non-linearly implicit numerical scheme to update $I_0, I_1$ and $T$. Moreover, the non-linearity is due to the term $T^4$. To handle this non-linear term, we introduce a new variable as 
\begin{equation}
    \label{eq:new_T}
    \psi = T^4.
\end{equation}
Then, the iteration method is utilized to solve this nonlinear system, and \eqref{eq:first_RTE_I0}, \eqref{eq:first_RTE_I1}, \eqref{eq:first_RTE_T} are reduced to the form as 
\begin{subnumcases}
   { \label{eq:first_RTE_1}}
    \label{eq:first_RTE_1_I0}
       \frac{\epsilon^2}{c} \frac{I_{0}^{n+1, k+1} - I_{0}^n }{\Delta t}
      +  \epsilon^2\left(\pd{\hat{I}_{1}}{x}\right)^{n+1, k+1} =\sigma \left[ac \psi^{n+1, k+1} - I_{0}^{n+1, k+1} \right],\\
      \label{eq:first_RTE_1_I1}
       \frac{\epsilon^2}{c} \frac{\hat{I}_{1}^{n+1, k+1} - \hat{I}_{1}^n }{\Delta t} + a_0\left(\pd{I_{0}}{x}\right)^{n+1, k+1} + b_2\left(\pd{I_{2}}{x}\right)^{n} = - \sigma\hat{I}_{1}^{n+1, k+1},\\
       \label{eq:first_RTE_1_T}
       \epsilon^2 C_{v} \frac{G^{n+1, k+1} - T^{n}}{\Delta t}  =
    \sigma\left[ I_{0}^{n+1, k+1} - ac \psi^{n+1, k+1}\right], \qquad G^{n+1, k+1} = \frac{\psi^{n+1, k+1}}{(T^{n+1, k})^3}, 
\end{subnumcases}
where $k = 0, 1, \cdots,$ is the superscript for the iteration, while  $()^{n+1, k+1}$ is the numerical solution after the $(k+1)$-th iteration of \eqref{eq:first_RTE_1} with $()^{n+1, 0} = ()^n$. With the new variable $\psi$, \eqref{eq:first_RTE_1} is changed into a linear system of $I_0, I_1$ and $\psi$, and the computational cost is greatly reduced. Moreover, we want to emphasize that when solving \eqref{eq:first_RTE_1}, the coefficient matrix degenerates into a singular matrix when $\epsilon$ approaches zero. To avoid this, we first substitute \eqref{eq:first_RTE_1_T} into \eqref{eq:first_RTE_1_I0} to eliminate the variable $I_0^{n+1, k+1}$. Then,  the equations related to $\hat{I}_1^{n+1, k+1}$ \eqref{eq:first_RTE_1_I1} and $T^{n+1, k+1}$ \eqref{eq:first_RTE_1_T} are solved to update the numerical solution. Finally, $I_0^{n+1, k+1}$ is obtained with \eqref{eq:first_RTE_1_T}.  

For the expansion coefficients $\{I_l, 2\leqslant l \leqslant M\}$, similar to \eqref{eq:first_linear_II}, $I_l^{n+1}$ are obtained successively according to \eqref{eq:first_RTE_Il}, which is solved at the numerical cost of an explicit scheme. 
\begin{remark}
    In the numerical scheme \eqref{eq:first_RTE_1}, the iteration method is utilized to renew the numerical solution. In the simulation, the criteria to stop the iteration are as below 
    \begin{itemize}
        \item $|T^{n+1, k+1} - T^{n+1, k}| < \bar{\epsilon}$, 
        \item $k < N_0$,
    \end{itemize}
    where $\bar{\epsilon}$ and $N_0$ are problem dependent. To summarize, we have Alg. \ref{alg:first_RTE}, the one-time step update of RTE. 

\begin{algorithm}[htpb]
\SetKwInput{KwInput}{Input}
\SetKwInput{KwOutput}{Output}
\caption{One step of first-order semi-discrete update for \eqref{eq:1D_RTE}}
\label{alg:first_RTE}
\KwInput{$\{I_l^n, l = 0, \cdots, M\},T^n$ }
\KwOutput{$\{I_l^{n+1}, l = 0, \cdots, M\},T^{n+1}$}
\begin{enumerate}
    \item obtain $\psi^n$ from \eqref{eq:new_T};
    \item obtain $I_1^{n+1}, \psi^{n+1}$ and $I_0^{n+1}$ from \eqref{eq:first_RTE_I0}, \eqref{eq:first_RTE_I1}, \eqref{eq:first_RTE_T};
        \begin{enumerate}
            \item Let $k \gets 0$, and $I_0^{n+1, k} = I_0^{n}$, $I_1^{n+1, k} = I_1^{n}$, $\psi^{n+1, k} = \psi^n$, and $T^{n+1, k} = T^n$;
            \item Obtain $I_1^{n+1, k+1}$ and $\psi^{n+1, k+1}$ from \eqref{eq:first_RTE_1};
            \item Obtain $I_0^{n+1, k+1}$ from \eqref{eq:first_RTE_1_T}; 
            \item Obtain $T^{n+1, k+1}$ from \eqref{eq:new_T};
            \item Let $k \gets k+1$ and go to Step (b) until $k < N_0$ or $|T^{n+1, k+1} - T^{n+1, k}| < \bar{\epsilon}$;
        \end{enumerate}
    \item obtain $I_l^{n+1}, 2\leqslant l \leqslant M$ from \eqref{eq:first_RTE_Il}. 
\end{enumerate}
\end{algorithm}

\end{remark}

\paragraph{Second-order temporal discretization for RTE}
Additionally, this first-order numerical scheme can be extended to a second-order scheme. Supposing $I_l^{n}$ and $T^n$ are the numerical solutions at the $n$-th time level, and the numerical solutions at $(n+1/2)$-th time level have already been obtained through the first order scheme \eqref{eq:first_RTE}, then the second-order numerical scheme has the form below 
{\small 
\begin{subnumcases}
{  \label{eq:second_RTE}}
    \label{eq:second_RTE_I0}
      \frac{\epsilon^2}{c} \frac{I_{0}^{n+1} - I_{0}^n }{\Delta t}
      + \frac{\epsilon^2}{2}\left[\left(\pd{\hat{I}_{1}}{x}\right)^{n+1} + \left(\pd{\hat{I}_{1}}{x}\right)^{n}\right]=\frac{\sigma}{2} \big[\left(ac (T^{4})^{n+1} - I_{0}^{n+1}\right) + \left( ac (T^4)^n - I_{0}^{n}\right)\big],\\
      \label{eq:second_RTE_I1}
      \frac{\epsilon^2}{c} \frac{\hat{I}_{1}^{n+1} - \hat{I}_{1}^n }{\Delta t} + \frac{1}{2}a_0\left[\left(\pd{I_{0}}{x}\right)^{n+1} + \left(\pd{I_{0}}{x}\right)^{n}\right] +  b_2\left(\pd{I_{2}}{x}\right)^{n + \frac{1}{2}}
      = - \frac{\sigma}{2}\left[\hat{I}_{1}^{n+1} +  \hat{I}_{1}^{n}\right], \\
       \label{eq:second_RTE_T}
      \epsilon^2 C_{v} \frac{T^{n+1} - T^{n}}{\Delta t}  =
    \frac{\sigma}{2}\left[\left( I_{0}^{n+1} - ac (T^{n+1})^4\right) + \left( I_{0}^{n} - ac (T^{n})^4\right)\right],\\
     \label{eq:second_RTE_Il}
     \frac{\epsilon^2}{c} \frac{I_{l}^{n+1} - I_{l}^n
      }{\Delta t} +\frac{ \epsilon a_{l-1}}{2}\left[\left(\pd{I_{l-1}}{x}\right)^{n+1} + \left(\pd{I_{l-1}}{x}\right)^{n}\right]
     + \epsilon b_{l+1}\left(\pd{I_{l+1}}{x}\right)^{n + \frac{1}{2}} =
     -\frac{ \sigma}{2} \left[I_l^{n+1} + I_l^{n}\right].
\end{subnumcases}
}For this second-order scheme, the same strategies as introducing the new variable $\psi = T^4$ \eqref{eq:new_T} and solving the linear equation system of $I_1$ and $T$ \eqref{eq:first_RTE_1} are also employed. For now, the AP scheme for the constant opacity is proposed, and we will discuss the numerical scheme for the non-linear opacity. 

\subsection{The nonlinearity in the opacity}
\label{sec:RTE_Marshak}
When the opacity $\sigma$ depends on $T$ non-linearly, a numerical scheme similar to \eqref{eq:first_RTE} can be utilized, and we consider a specific example with 
\begin{equation}
    \label{eq:sigma}
    \sigma = \frac{\rho}{T^3},
\end{equation}
where $\rho$ is constant, which is also the model corresponding to the classical Marshak-wave type problem. 

Substituting \eqref{eq:sigma} into \eqref{eq:1D_RTE}, the interaction term $\sigma(acT^4/2 - I)$ is changed into 
\begin{equation}
    \label{eq:right} 
    \sigma(acT^4/2-I) = \frac{ac\rho T}{2} - \frac{\rho I}{T^3}. 
\end{equation}
Then $I_0, I_1$ and $T$ are updated in the same way as in \eqref{eq:first_RTE}, and the IMEX scheme is utilized to update $I_l, l > 1$. The first-order scheme has the form below
{\small 
\begin{subnumcases}
        {\label{eq:first_RTE_non}}
        \label{eq:first_RTE_non_I0}
       \frac{\epsilon^2}{c} \frac{I_{0}^{n+1} - I_{0}^n }{\Delta t}
      +  \epsilon^2\left(\pd{\hat{I}_{1}}{x}\right)^{n+1} = ac \rho T^{n+1} - \frac{\rho I_0^{n+1}}{(T^{n+1})^3}, \\
              \label{eq:first_RTE_non_I1}
             \frac{\epsilon^2}{c} \frac{\hat{I}_{1}^{n+1} - \hat{I}_{1}^n }{\Delta t} + a_0\left(\pd{I_{0}}{x}\right)^{n+1} + b_2\left(\pd{I_{2}}{x}\right)^{n} = - \frac{\rho \hat{I}_{1}^{n+1}}{(T^{n+1})^3},\\
               \label{eq:first_RTE_non_T}
       \epsilon^2 C_{v} \frac{T^{n+1} - T^{n}}{\Delta t}  = \frac{\rho I_0^{n+1}}{(T^{n+1})^3} - ac \rho T^{n+1}, \\
               \label{eq:first_RTE_non_Il}
    \frac{\epsilon^2}{c} \frac{I_{l}^{n+1} - I_{l}^n
      }{\Delta t} + \epsilon a_{l-1}\left(\pd{I_{l-1}}{x}\right)^{n+1} +
      \epsilon b_{l+1}\left(\pd{I_{l+1}}{x}\right)^{n}= -\frac{\rho I_l^{n+1}}{(T^{n+1})^3}, \qquad l = 2, \cdots, M.
  \end{subnumcases}
}This is also a non-linear system, which is expensive to solve. To reduce the computational cost, an iteration method similar to \eqref{eq:first_RTE_1} is utilized to update $I_0, I_1$ and $T$, and the resulting scheme of \eqref{eq:first_RTE_non_I0}, \eqref{eq:first_RTE_non_I1} and \eqref{eq:first_RTE_non_T} is as follows 
{\small 
\begin{subnumcases}
  { \label{eq:first_RTE_non_1}}
        \label{eq:first_RTE_non_1_I0}
       \frac{\epsilon^2}{c} \frac{I_{0}^{n+1, k+1} - I_{0}^n }{\Delta t}
      +  \epsilon^2\left(\pd{\hat{I}_{1}}{x}\right)^{n+1, k+1} = ac \rho T^{n+1, k+1} - \frac{\rho I_0^{n+1, k+1}}{(T^{n+1, k})^3}, \\
        \label{eq:first_RTE_non_1_I1}
             \frac{\epsilon^2}{c} \frac{\hat{I}_{1}^{n+1, k+1} - \hat{I}_{1}^n }{\Delta t} + a_0\left(\pd{I_{0}}{x}\right)^{n+1, k+1} + b_2\left(\pd{I_{2}}{x}\right)^{n} = - \frac{\rho \hat{I}_{1}^{n+1, k+1}}{(T^{n+1, k})^3},\\
        \label{eq:first_RTE_non_1_T}
       \epsilon^2 C_{v} \frac{T^{n+1, k+1} - T^{n}}{\Delta t}  = \frac{\rho I_0^{n+1, k+1}}{(T^{n+1, k})^3} - ac \rho T^{n+1, k+1}.
  \end{subnumcases}
}In this case, \eqref{eq:first_RTE_non_1} reduces to a linear system of $I_0^{n+1,k+1}, I_1^{n+1, k+1}$ and $T^{n+1, k+1}$, resulting in the decreasing of the computational cost.  Moreover, the same strategy to obtain $I_1^{n+1, k+1}$ and $T^{n+1, k+1}$ as in Sec. \ref{sec:RTE_con_sig} is employed. When updating $I_l^{n+1, k+1}$ in \eqref{eq:first_RTE_non_Il}, the updated $T^{n+1}$ is utilized. For now, the first-order temporal discretization is introduced, and the one-time update is summarized in Alg. \ref{alg:first_RTE_non}. 

\begin{algorithm}[htpb]
\SetKwInput{KwInput}{Input}
\SetKwInput{KwOutput}{Output}
\caption{One step of first-order semi-discrete update for \eqref{eq:1D_RTE} with nonlinear opacity}
\label{alg:first_RTE_non}
\KwInput{$\{I_l^n, l = 0, \cdots, M\},T^n$ }
\KwOutput{$\{I_l^{n+1}, l = 0, \cdots, M\},T^{n+1}$}
\begin{enumerate}
    \item obtain $I_1^{n+1}, T^{n+1}$ and $I_0^{n+1}$ from \eqref{eq:first_RTE_non_1};
        \begin{enumerate}
            \item Let $k \gets 0$, and $I_0^{n+1, k} = I_0^{n}$, $I_1^{n+1, k} = I_1^{n}$, and $T^{n+1, k} = T^n$;
            \item Obtain $I_1^{n+1, k+1}$ and $T^{n+1, k+1}$ from \eqref{eq:first_RTE_non_1};
            \item Obtain $I_0^{n+1, k+1}$ from \eqref{eq:first_RTE_non_1_T}; 
            \item Let $k \gets k+1$ and go to Step (b) until $k < N_0$ or $|T^{n+1, k+1} - T^{n+1, k}| < \bar{\epsilon}$;
        \end{enumerate}
    \item obtain $I_l^{n+1}, 2\leqslant l \leqslant M$ from \eqref{eq:first_RTE_non_Il}. 
\end{enumerate}
\end{algorithm}

\paragraph{Second-order temporal discretization for RTE with nonlinear opacity}
This first-order scheme can be extended into the second-order scheme naturally. Supposing $I_l^{n}$ and $T^n$ are the numerical solutions at the $n$-th time level, and the numerical solutions at the $(n+1/2)$-th time level have already been obtained through the first-order scheme \eqref{eq:first_RTE_non}, then the second-order numerical scheme has the form below
{\small 
\begin{subnumcases}
        {\label{eq:second_RTE_non}}
       \frac{\epsilon^2}{c} \frac{I_{0}^{n+1} - I_{0}^n }{\Delta t}
      +  \frac{\epsilon^2}{2}
            \left[ \left(\pd{\hat{I}_{1}}{x}\right)^{n+1}+ \left(\pd{\hat{I}_{1}}{x}\right)^{n} \right]= ac \rho \frac{T^{n+1} + T^n}{2} - \frac{\rho (I_0^{n+1} + I_0^n)}{(T^{n+1})^3 + (T^n)^3}, \\
             \frac{\epsilon^2}{c} \frac{\hat{I}_{1}^{n+1} - \hat{I}_{1}^n }{\Delta t} + \frac{a_0}{2}
             \left[
             \left(\pd{I_{0}}{x}\right)^{n+1}  +  \left(\pd{I_{0}}{x}\right)^{n} 
                 \right]+ b_2\left(\pd{I_{2}}{x}\right)^{n+1/2} = - \frac{\rho (\hat{I}_{1}^{n+1} + \hat{I}_{1}^{n})}{(T^{n+1})^3+ (T^{n})^3},\\
       \epsilon^2 C_{v} \frac{T^{n+1} - T^{n}}{\Delta t}  = \frac{\rho (I_0^{n+1} + I_0^n)}{(T^{n+1})^3 + (T^n)^3} - ac \rho\frac{T^{n+1} + T^n}{2}, \\
    \frac{\epsilon^2}{c} \frac{I_{l}^{n+1} - I_{l}^n
      }{\Delta t} + \frac{\epsilon a_{l-1}}{2}\left[ \left(\pd{I_{l-1}}{x}\right)^{n+1}+ \left(\pd{I_{l-1}}{x}\right)^{n}\right] +
      \epsilon b_{l+1}\left(\pd{I_{l+1}}{x}\right)^{n+1/2}= -\frac{\rho (I_l^{n+1}+I_l^{n}) }{  (T^{n+1})^{3} +(T^n)^{3} }, \qquad l = 2, \cdots, M.
  \end{subnumcases}
}

For now, we have introduced the temporal discretization for the gray model of the RTE, and the spatial discretization will be proposed in the following section.

%% file: analysis.tex
%%%%%%%%%%%%%%%%%%%%%%%%%%%%%%%%%%%%%%%%%%%%%%%%%%%%%%%%%%%%%%%%%%%%%%%%%%%%%%%%%%%%%%%%%%%%%%%%%%%%%%%%%
\section{Spatial discretization and formal asymptotic analysis}
\label{sec:spatial_AP}
In this section, the spatial discretization of the $P_N$ system is first proposed in Sec. \ref{sec:spatial}, and the asymptotic preserving (AP) property of the related numerical scheme for RTE \eqref{eq:1D_RTE} is discussed in Sec. \ref{sec:AP}, with the numerical stability by Fourier analysis presented in Sec. \ref{sec:Fourier}

\subsection{Spatial discretization}
\label{sec:spatial}
For spatial discretization, the finite volume method is utilized for the $P_N$ system. Take the spatially 1D problem as an example where the uniform mesh in the spatial space is utilized. Let $\{I_{l, i}^n, l = 0, \cdots, M; i = 0, \cdots, N-1\}$ and $T_i^n$ be the numerical solution on the cell $\{x: x_{i-1/2} < x < x_{i+1/2}\}$ at the $n$-th time step level, where $N$ is the total mesh number. Let \eqref{eq:first_linear_I} and \eqref{eq:first_linear_II} as an example and the local Lax-Friedrichs flux is adopted here. We first discuss the spatial discretization of \eqref{eq:first_linear_I}. The convection term is discretized as 
\begin{equation}
    \label{eq:flux1}
    \begin{aligned}
        \left(\pd{\hat{I}_1}{x} \right)_i^{n+1} & \approx  \frac{\mF^{n+1}_{-1, i+\frac{1}{2}} - \mF^{n+1}_{-1, i-\frac{1}{2}}}{\Delta x}  +  \frac{\mG^{n+1}_{1, i+\frac{1}{2}} - \mG^{n+1}_{1, i-\frac{1}{2}}}{\Delta x}, \\
        a_0 \left(\pd{I_0}{x} \right)_i^{n+1}  +  b_2 \left(\pd{I_2}{x} \right)_i^{n} & \approx \frac{\mF^{n+1}_{0, i+\frac{1}{2}} - \mF^{n+1}_{0, i-\frac{1}{2}}}{\Delta x} + \frac{\mG^{n}_{2, i+\frac12} - \mG^{n}_{2, i-\frac12}}{\Delta x}, 
    \end{aligned}
\end{equation}
with 
\begin{subequations}
    \label{eq:FG1}
    \begin{align}
    \label{eq:FG1_0}
    &  \mathcal{F}_{-1,i+\frac{1}{2}}^{n+1} = -\frac{v_{\max}}{2}(I_{0,i+1}^{n + 1} - I_{0,i}^{n + 1}), \qquad \mathcal{G}_{1,i+\frac{1}{2}}^{n+1} = \frac{1}{2}(\hat{I}_{1,i+1}^{n+1} + \hat{I}_{1,i}^{n+1}), \\
    \label{eq:FG1_1}
    & \mathcal{F}_{0,i+\frac{1}{2}}^{n+1} = \frac{a_{0}}{2}(I_{0,i+1}^{n+1} + I_{0,i}^{n+1}) - \frac{\epsilon v_{\max}}{2}(\hat{I}_{1,i+1}^{n + 1} - \hat{I}_{1,i}^{n + 1}), \qquad 
    \mathcal{G}_{2,i+\frac{1}{2}}^{n} = \frac{b_{2}}{2}(I_{2,i+1}^{n} + I_{2,i}^{n}),
    \end{align}
\end{subequations}
where $v_{\max}$ is the local wave speed. For the second part \eqref{eq:first_linear_II}, since the convection term is split into the implicit and explicit parts, the detailed spatial discretization has the form below. 
\begin{equation}
\label{eq:convection}
\epsilon a_{l-1}\left(\pd{I_{l-1}}{x}\right)_i^{n+1} + \epsilon b_{l+1}\left(\pd{I_{l+1}}{x}\right)_i^{n} \approx  \frac{\mathcal{F}_{l-1,i+\frac{1}{2}}^{n+1} - \mathcal{F}_{l-1,i-\frac{1}{2}}^{n+1}}{\Delta x} +  \frac{\mathcal{G}_{l+1,i+\frac{1}{2}}^{n} - \mathcal{G}_{l+1,i-\frac{1}{2}}^{n}}{\Delta x}, \qquad l = 2, \cdots, M. 
\end{equation}
with 
\begin{equation}
\label{eq:flux_LF}
\mathcal{F}_{l-1,i+\frac{1}{2}}^{n+1} = \frac{\epsilon a_{l-1}}{2}(I_{l - 1,i+1}^{n+1} + I_{l - 1,i}^{n+1}),\qquad 
\mathcal{G}_{l + 1,i+\frac{1}{2}}^{n} = \frac{\epsilon b_{l+1}}{2}(I_{l + 1,i+1}^{n} + I_{l + 1,i}^{n}) - \frac{\epsilon v_{\max}}{2}(I_{l,i+1}^{n} - I_{l,i}^{n}).
\end{equation}

In the $P_N$ system \eqref{eq:final_Pn_linear}, the maximum characteristic velocity is the maximum positive root of the $(M+1)$-th order Legendre polynomial with $M$ the expansion order of the $P_N$ system, which is always smaller than $1$. Therefore, we set the wave speed in the local Lax-Friedrichs flux as $v_{\max} = 1$. Here, we want to emphasize that since we need to solve a linear equation system to obtain $I_0$ and $I_1$, the diffusion term in the Lax-Friedrichs flux  \eqref{eq:FG1} is treated implicitly. Moreover, $I_{M+1}$ is set as zero in the numerical flux \eqref{eq:flux_LF}. Additionally, reconstruction in the spatial space such as the linear reconstruction is utilized to reduce the computational cost.

In the reconstruction, the numerical solution at $t^{n}$ is utilized to obtain the reconstruction slope for the numerical solution at $t^{n+1}$ for linear radiation transfer equations, while the numerical solution at the $k-$th iteration at $t^{n+1}$ is adopted to obtain the reconstruction slope for the numerical solution at the $(k+1)-$th iteration for the gray model of radiation transfer equations.

Combining the temporal discretization in Sec. \ref{sec:RTE} and the spatial discretization  \eqref{eq:convection}, we obtain an IMEX numerical scheme for the gray model of RTE \eqref{eq:1D_RTE} and the linear system \eqref{eq:1D_linear_RTE} in the framework of the $P_N$ method. In the following Sec. \ref{sec:AP}, the asymptotic preserving (AP) property of these numerical schemes will be discussed.

\subsection{Formal asymptotic analysis}
\label{sec:AP}
In this section, the AP property of this IMEX numerical scheme is discussed. The AP property is of great importance for the multi-scale problems. If the numerical scheme has the AP property, it will reduce to a numerical scheme of the diffusion limit equation when the multi-scale parameter $\epsilon$ goes to zero, with the mesh size and the time step length remaining unchanged \cite{sun2015asymptotic1, Klar1998, Larsen1987}. Without the AP property, grids and time step length consistent with the parameter $\epsilon$ are needed to resolve the diffusion limit equation \cite{Mie2013, Larsen1989}, which will be quite expensive when $\epsilon$ is small.

For the gray model of RTE, when the parameter $\epsilon$ approaches zero, it will achieve the diffusion limit \eqref{eq:1D_limit}. Therefore, it is expected this IMEX numerical scheme will reduce to a numerical scheme of the diffusion limit with $\epsilon$ approaching zero. The main results for the linear RTE \eqref{eq:1D_linear_RTE} and RTE \eqref{eq:1D_RTE} are listed in the two theorems below. The proofs of these two theorems are similar, and only the proof of Thm. \ref{thm:AP_RTE} is presented here.

\begin{theorem} [AP property of the IMEX scheme for linear RTE \eqref{eq:1D_linear_RTE}]
\label{thm:AP_RTE_linear}
  With parameter $\epsilon$ approaching zero, the numerical scheme \eqref{eq:first_linear_I} and \eqref{eq:first_linear_II} together with the spatial discretization \eqref{eq:flux1} and \eqref{eq:convection} reduces to an implicit five-point scheme for the diffusion limit \eqref{eq:1D_linear_limit}. 
\end{theorem}
  
\begin{theorem} [AP property of the IMEX scheme for RTE \eqref{eq:1D_RTE}]
\label{thm:AP_RTE}
  With parameter $\epsilon$ approaching zero, the numerical scheme \eqref{eq:first_RTE} together with the spatial discretization \eqref{eq:flux1} and \eqref{eq:convection} reduces to an implicit five-point scheme for the diffusion limit \eqref{eq:1D_limit}. 
\end{theorem}

\begin{proof}[Proof of Thm. \ref{thm:AP_RTE}]
Without loss of generality, we assume $\sigma$ is constant, and the periodic boundary condition is employed. In the semi-discretization \eqref{eq:first_RTE}, summarizing \eqref{eq:first_RTE_I0} and \eqref{eq:first_RTE_T} up, with the spatial discretizations \eqref{eq:flux1} and \eqref{eq:FG1}, we can derive that 
\begin{equation}
\label{eq:AP_1}
    \frac{1}{c} \frac{I_{0, i}^{n+1} - I_{0, i}^{n}}{\Delta t} +  C_v \frac{T_{i}^{n+1} - T_{i}^n}{\Delta t} +    
    \frac{1}{\Delta x}(\hat{I}_{1, i+1}^{n+1} - \hat{I}_{1, i-1}^{n+1}) -\frac{1}{2\Delta x}(I_{0,i+1}^{n+1} - 2I_{0, i}^{n+1} + I_{0, i-1}^{n+1}) = 0.
%    \frac{1}{c} \frac{I_{0, i}^{n+1} - I_{0,i}^{n}}{\Delta t} + \frac{\mF_{-1, i+1/2}^{n+1} - \mF_{-1, i-1/2}^{n+1}}{\Delta x}  + \frac{\mG_{1, i+1/2}^{n+1} - \mG_{1, i-1/2}^{n+1}}{\Delta x} = 0.
\end{equation}
From the order analysis \eqref{eq:order_PN_RTE}, it holds for $I_2$ 
\begin{equation}
    \label{eq:AP_2} 
    I_2 = \mO(\epsilon^2). 
\end{equation}
Then, we can derive from \eqref{eq:first_RTE_I1} and the spatial discretization \eqref{eq:flux1} and \eqref{eq:FG1} 
\begin{equation}
    \label{eq:AP_3}
        -\sigma \hat{I}_{1, i}^{n+1} = \frac{\mF_{0, i+1/2}^{n+1} - \mF_{0, i-1/2}^{n+1}}{\Delta x} + \mO(\epsilon^2) =  \frac{a_0}{2\Delta x} (I_{0, i+1}^{n+1} - I_{0, i-1}^{n+1}) +
        \mO(\Delta x) + \mO(\epsilon^2), \qquad a_0 = \frac{1}{3}.
\end{equation}
Substituting \eqref{eq:AP_3} into \eqref{eq:AP_1},  we can derive 
\begin{equation}
\label{eq:AP_4}
\begin{aligned}
      \frac{1}{c} \frac{I_{0, i}^{n+1} - I_{0, i}^{n}}{\Delta t} & + C_v \frac{T_{i}^{n+1} - T_{i}^n}{\Delta t}  \\
        & +  \frac{1}{2\Delta x}\left(\frac{-a_0}{2\sigma  \Delta x} (I_{0, i+2}^{n+1} -I_{0, i}^{n+1}) - \frac{-a_0}{2\sigma  \Delta x} (I_{0, i}^{n+1} -I_{0, i-2}^{n+1})\right) = \mO(\Delta x) + \mO(\epsilon^2). 
  \end{aligned}
\end{equation}
When $\epsilon$ approaches zero, form \eqref{eq:first_RTE_I0}, we can deduce 
\begin{equation}
    \label{eq:AP_5}
    I_0 = ac T^4 \triangleq 2 I^{(0)} = ac \left(T^{(0)}\right)^4. 
\end{equation}
Substituting \eqref{eq:AP_5} into \eqref{eq:AP_4}, \eqref{eq:AP_4} is changed into
\begin{equation}
    \label{eq:AP_6}
    \begin{aligned}
      a \frac{G_i^{n+1} - G_i^{n}}{\Delta t} & + C_v \frac{\left(T_{i}^{(0)}\right)^{n+1} - \left(T_{i}^{(0)}\right)^n}{\Delta t}  \\
        & =  \frac{1}{2\Delta x}\left(\frac{ac}{3\sigma} \frac{G_{i+2}^{n+1} -G_{i}^{n+1}}{2\Delta x} - \frac{a c}{3\sigma}\frac{G_{i}^{n+1} -G_{i-2}^{n+1}}{2\Delta x}\right), \qquad G = \left(T^{(0)}\right)^4,
  \end{aligned}
\end{equation}
with the small terms $\mO(\Delta x)$ and $\mO(\epsilon^2)$ omitted, and the proof is completed. 
\end{proof}

For the numerical schemes with second-order temporal discretization proposed in Sec. \ref{sec:linear} and \ref{sec:RTE}, they all preserve the AP property with the spatial discretization \eqref{eq:flux1} and \eqref{eq:convection}, and the proof is omitted here. 

\subsection{Numerical stability by Fourier analysis}
\label{sec:Fourier}
In this section, the numerical stability of this AP numerical scheme for the linear system \eqref{eq:1D_linear_RTE} is studied by Fourier analysis. We follow the method proposed in \cite{jang2014analysis, peng2020} and a similar analysis is proposed in \cite{IMEX2022} for linear RTE. However, compared to the numerical scheme proposed in \cite{IMEX2022}, the convection term in \eqref{eq:first_linear_I0_1} is treated implicitly, which is supposed to enlarge the stable region of the numerical scheme, and will be shown through Fourier analysis. 

\begin{figure}[!hptb]
  \centering
  \subfloat[$M = 7$]{
    \includegraphics[width=0.3\textwidth]{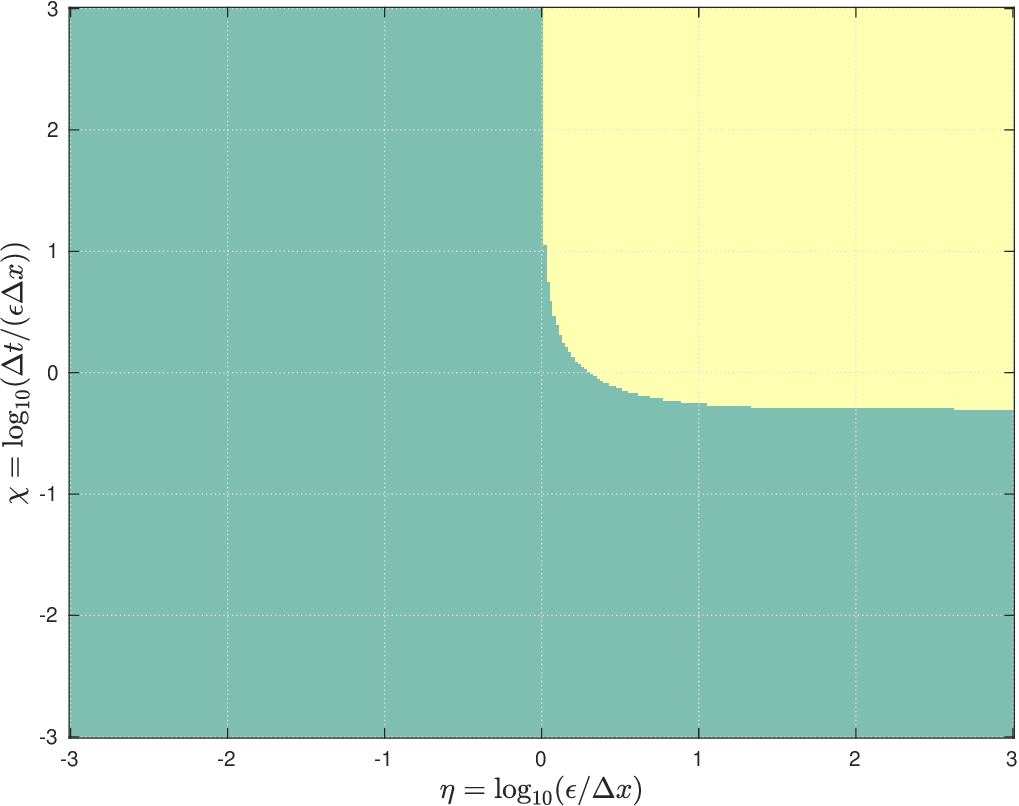}}
  \hfill
  \subfloat[$M = 15$]{
    \includegraphics[width=0.3\textwidth]{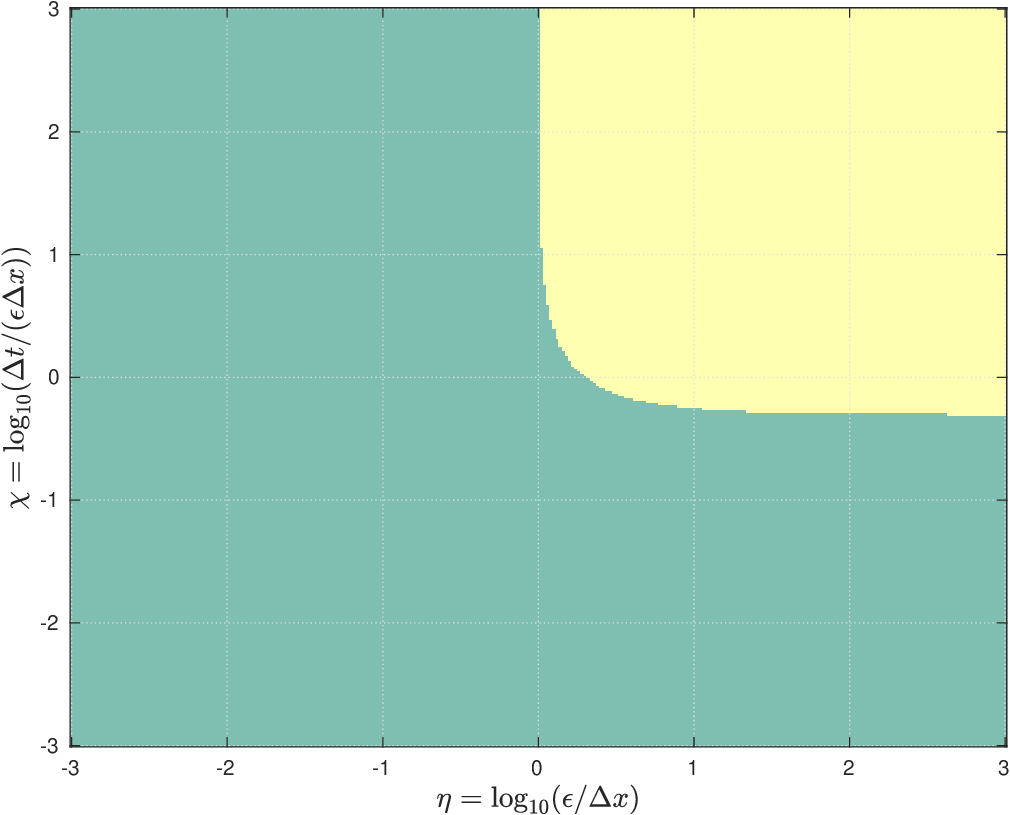}}
     \hfill
  \subfloat[$M = 31$]{
    \includegraphics[width=0.3\textwidth]{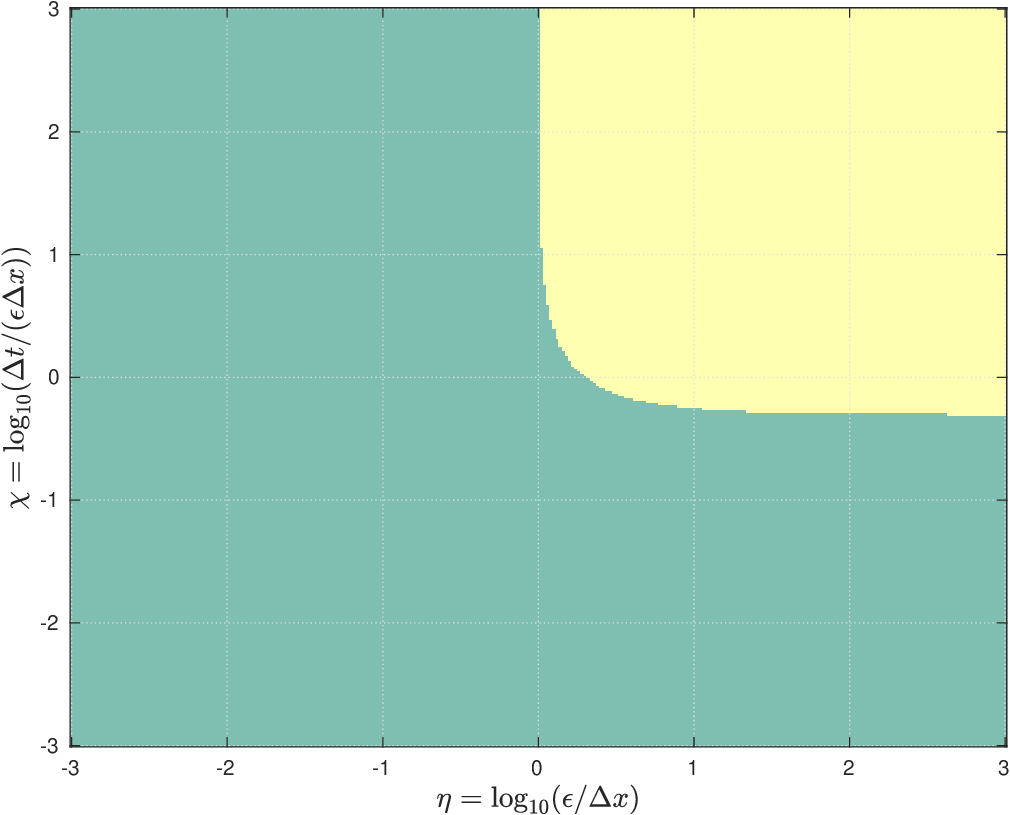}}
  \caption{(Numerical stability by Fourier analysis in Sec. \ref{sec:Fourier}) Numerical stability of the numerical scheme for the linear system \eqref{eq:1D_linear_RTE} by Fourier analysis for the expansion number $M = 7, 15$ and $31$. Here the blue area is the stable region and the yellow area is the unstable region. }
  \label{fig:stability_Fourier}
\end{figure}

To carry out the Fourier analysis, we assume the periodic boundary condition is applied to the spatial space and the mesh is uniform with the parameter $c = 1, \sigma_a = 0$, and $\sigma_s = 1$. With the spatial discretization in Sec. \ref{sec:spatial} and $\hat{I}_1$ substituted with $I_1/\epsilon$, the numerical scheme \eqref{eq:first_linear_I} and \eqref{eq:first_linear_II} is reduced to 
{\small 
\begin{subequations}
    \label{eq:linear_Fourier}
    \begin{align}
    \label{eq:linear_Fourier_I0}
        &\frac{\epsilon^2}{\Delta t}(I_{0, i}^{n+1} - I_{0, i}^{n})   + \frac{\epsilon}{2 \Delta x} (\II{1}{i+1}{n+1} - \II{1}{i-1}{n+1})  - \frac{\epsilon^2}{2\Delta x}\left(I_{0,i+1}^{n+1} - 2I_{0,i}^{n+1} + \II{0}{i-1}{n+1}\right)= 0, \\\nonumber
    \label{eq:linear_Fourier_I1}
        &\frac{\epsilon^2}{\Delta t}(\II{1}{i}{n+1} - \II{1}{i}{n})   + \frac{\epsilon a_0}{2 \Delta x} (\II{0}{i+1}{n+1} - \II{0}{i-1}{n+1}) - \frac{\epsilon}{2\Delta x}(\II{1}{i+1}{n+1}- 2\II{1}{i}{n+1} + \II{1}{i-1}{n+1})\\
         &  \hspace{8cm}  
         + \frac{\epsilon b_2}{2 \Delta x} (\II{2}{i+1}{n} - \II{2}{i-1}{n}) = -\II{1}{i}{n+1},  \\\nonumber
    \label{eq:linear_Fourier_Il}
        &\frac{\epsilon^2}{\Delta t}(\II{l}{i}{n+1} - \II{l}{i}{n})   + \frac{\epsilon a_{l-1}}{2 \Delta x} (\II{l-1}{i+1}{n+1} - \II{l-1}{i-1}{n+1}) - \frac{\epsilon}{2\Delta x}(\II{l}{i+1}{n} - 2\II{l}{i}{n} + \II{l}{i-1}{n}) \\
        &\hspace{8cm}
         + \frac{\epsilon b_{l+1}}{2 \Delta x} (\II{l+1}{i+1}{n} - \II{l+1}{i-1}{n}) = -\II{l}{i}{n+1}.
    \end{align}
\end{subequations}
}Taking the Fourier ansatz
\begin{equation}
  \label{eq:Fourier_scheme}
   I_{l,i}^n =  \tilde{I}_{l,k}^n \exp(\imag k i \Delta x),  \qquad l = 0, 1, \dots, M ,\quad i = 0, \cdots, N-1, \quad k \in \bbZ,
\end{equation}
where $k$ is the index for the Fourier mode, $i$ is the index in the $x-$axis, and $\Delta x$ is the mesh size. Then, the linear scheme \eqref{eq:first_linear_I} and \eqref{eq:first_linear_II} will render 
\begin{equation}
  \label{eq:fourier_mat}
 {\bf A}(\epsilon,  \Delta t, \Delta x, \xi) \tilde{\bm{I}}^{n+1}_k = {\bf B}(\epsilon,  \Delta t, \Delta x, \xi) \tilde{\bm I}^{n}_k,
  \end{equation}
where  $\tilde{\bm{I}}_k^{n+1} = (\tilde{I}_{0, k}^{n+1},\tilde{I}_{1, k}^{n+1}, \cdots, \tilde{I}_{M, k}^{n+1})^T$, and ${\bf A}, {\bf B}$ are $(M+1) \times (M+1)$ coefficient matrices dependent on the model parameter $\epsilon$, the mesh size $\Delta x$, the time step length $\Delta t$, and the discrete wave number $\xi = k \Delta x \in [0, 2\pi]$. The detailed forms of the matrices ${\bf A}$ and ${\bf B}$ are as below 
{\small
\begin{equation}
    \label{eq:fourier_A_B}
    \left( \begin{array}[c]{cccccc}
  c_0 - \epsilon c_1& c_2 & 
  0 &\cdots& 0 & 0\\ 
  a_0 c_2 & c_0 + 1 - c_1 & 0 &\cdots& 0 & 0\\
  0 & a_1 c_2  & c_0 + 1  &\cdots& 0 & 0\\
  \vdots&\vdots&\vdots&\vdots&\vdots&\vdots\\
  0&0&0&\cdots&a_{M-1} c_2 & c_0 + 1
  \end{array}\right), \qquad 
  \left( \begin{array}[c]{ccccccc}
  c_0& 0 & 0&
  0 &\cdots& 0 & 0\\ 
  0& c_0  & -b_2 c_2 &0 &\cdots& 0 & 0\\
  0 &0&c_0 + c_1 & -b_3 c_2 &\cdots& 0 & 0\\
  \vdots&\vdots&\vdots&\vdots&\vdots&\vdots&\vdots\\
  0&0&0&0&\cdots&0& c_0 + c_1
  \end{array}\right),
\end{equation}
}
with 
\begin{equation}
    \label{eq:coe_c}
    c_0 = \frac{\epsilon^2}{\Delta t}, \qquad c_1 = \frac{\epsilon (\cos(\xi) - 1)}{\Delta x}, \qquad c_2 = \frac{\imag \epsilon \sin(\xi)}{\Delta x}. 
\end{equation}
\begin{remark}
Due to the special form of $I_0$ and $I_1$ \eqref{eq:first_linear_I}, and the numerical flux \eqref{eq:flux1} and \eqref{eq:FG1}, the first two diagonal elements of matrix ${\bf A}$ are different from others. This is due to the implicit viscous terms of the numerical flux in \eqref{eq:flux1}. Moreover, the parameter $\epsilon$ in the first diagonal element $A_{11}$ of $\bf A$ is led by the new variable $\hat{I}_1$. In this case, the parameter $\epsilon$ does not appear in the semi-discrete equation \eqref{eq:first_linear_I0_1} of $I_0$, and will not appear in the numerical flux \eqref{eq:FG1_0}. Therefore, when \eqref{eq:linear_Fourier_I0} is deduced, $\epsilon^2$ exists in the implicit viscous term of the numerical flux, making the different form of $A_{11}$. 
% Notably, the first element of matrix A is $c_0 - \epsilon c_1$, where $c_1$ is preceded by the parameter $\epsilon$. The reason for this is that the coefficient in front of the numerical viscosity term in scheme \eqref{eq:linear_Fourier_I0}  is $\frac{\epsilon^2}{\Delta x}$.
\end{remark}
Then the amplification matrix ${\bf C}$ is obtained by ${\bf C} = {\bf A}^{-1} {\bf B}$. To analyze the numerical stability of this scheme by Fourier analysis,  the same stable principle as in \cite{peng2020} is utilized here as 
\paragraph{Principle for numerical stability:} If for all $\xi \in [0, 2\pi]$, the eigenvalues of ${\bf C}(\epsilon, \Delta t, \Delta x, \xi)$ satisfy either of the two conditions below, the numerical scheme is stable. 
\begin{enumerate}[Cond. 1:]
\item \label{linear_stability_1} $\max\limits_{i=1,\cdots,M+1}\{|\lambda_i(\xi)|\} < 1$, 
\item \label{linear_stability_2} $\max\limits_{i=1,\cdots,M+1}\{|\lambda_i(\xi)|\} =1$ and $\bf C$ is real and diagonalizable.
\end{enumerate}
As stated in \cite{peng2020}, this principle is a necessary condition for the standard $l^2$ energy stability. From the specific form of ${\bf C}$,  we can find that it is only related to $\frac{\epsilon}{\Delta x}, \frac{\epsilon^2}{\Delta t} = \frac{\epsilon}{\Delta x} \frac{\Delta x \epsilon}{\Delta t}$, $\epsilon$ and $\xi$. Two parameters $\eta = \log_{10}(\frac{\epsilon}{\Delta x})$ and $\chi = \log_{10}(\frac{\Delta t}{\epsilon \Delta x})$ are then introduced. The stability area for $M = 7, 15, 31$ is plotted in Fig. \ref{fig:stability_Fourier}, where the computational area for $\eta$ and $\chi$ is $[-4, 4]\times [-3, 3]$ with the discretization $\Delta \eta = \Delta \chi = 0.02$. The discretization for $\xi$ is $\Delta \xi = 0.01 \pi$. As to $\epsilon$, we have tested $\epsilon \in [0.001, 1]$ with the discretization $\Delta \epsilon = 0.001$. The numerical results indicate that $\epsilon$ does not affect the stability area since it only appears in $A(1,1)$  as shown in \eqref{eq:fourier_A_B}. Therefore, the stability area with $\epsilon = 1$ is shown in Fig. \ref{fig:stability_Fourier}.

It indicates that the stability area does not change according to the expansion order $M$. When $\eta < \eta_0$, the numerical scheme \eqref{eq:first_linear_I} and \eqref{eq:first_linear_II} is unconditionally stable. When $\eta > \eta_0$, the scheme is stable if $\chi < \chi_{0}$. In the simulation, we set $\eta_0 = 0$ and $\chi_0 = -0.5$, and based on this, the time step length is chosen as 
\begin{equation}
\label{eq:time_step}
\left\{
\begin{aligned}
&\Delta t = C \epsilon \Delta x / c,\qquad & \epsilon > \Delta x,\\
&\Delta t = C  \Delta x / c,\qquad &  \epsilon < \Delta x,
\end{aligned}
\right.
\end{equation}
where $C< 1$ is the CFL condition number. \\

From Fig. \ref{fig:stability_Fourier}, we can find that the stability region of the proposed scheme significantly extends that of the method in  \cite{IMEX2022}, enabling stability with larger time steps. Specifically, for
$\eta < \eta_0$, the proposed scheme demonstrates unconditional stability, while the stability of the scheme in \cite{IMEX2022} requires a time-step limit proportional to the square of the spatial step size, indicating its improved robustness under varying time-step constraints.

%% file: numerical_results.tex
\section{Numerical results}
\label{sec:num}
% three 1d example 
In this section, several numerical examples are tested to verify this $P_N$-based IMEX AP scheme (IMEX-IM) for the RTE. The AP property is first validated by a periodic 1D problem. Then, the classical 1D plane source problem and Marshak wave problem are studied. The 2D problems such as the 2D line source problem, lattice problem, and Riemann problem are tested to validate the high efficiency of IMEX-IM. The CFL number is set as $C = 0.4$ in both 1D and 2D numerical tests.

\subsection{The verification of the AP property}
\label{sec:ex1}
In this section, the AP property of IMEX-IM is tested. A similar initial condition as in \cite{IMEX2022} is adopted here, where the initial temperature and the radiation density are set as 
\begin{equation}
  \label{eq:initial_ex1}
  T =  (3 + \sin(\pi  x) ) / 4, \qquad I = \frac{1}{2}ac T^4, \qquad x \in L, \qquad L = [0, 2],
\end{equation}
with the parameters set as 
\begin{equation}
    \label{eq:ex1_param}
    a = c = 1.0, \qquad C_{v} = 0.1, \qquad \sigma = 10.    
\end{equation}

\begin{figure}[!hptb]
  \centering
  \subfloat[radiation temperature $T_r$]{
    \includegraphics[width=0.4\textwidth]{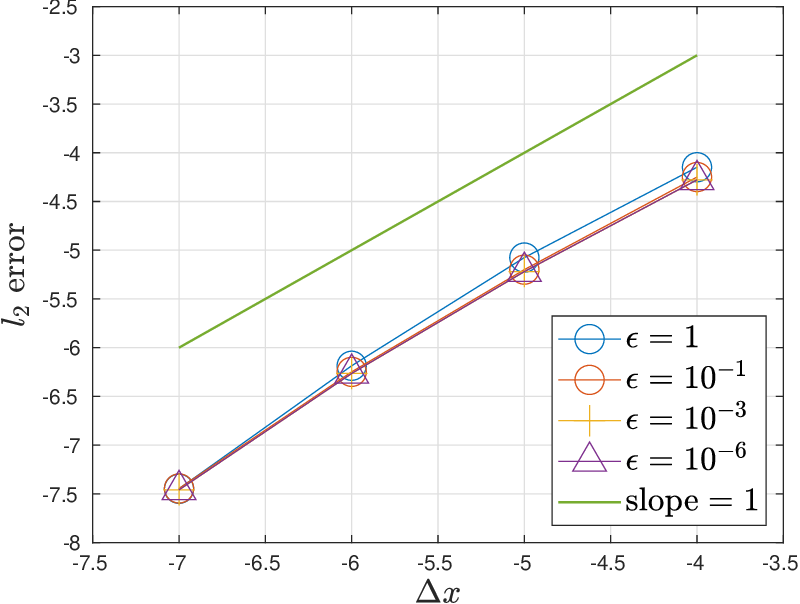}}
  \hfill
  \subfloat[material temperature $T$]{
    \includegraphics[width=0.4\textwidth]{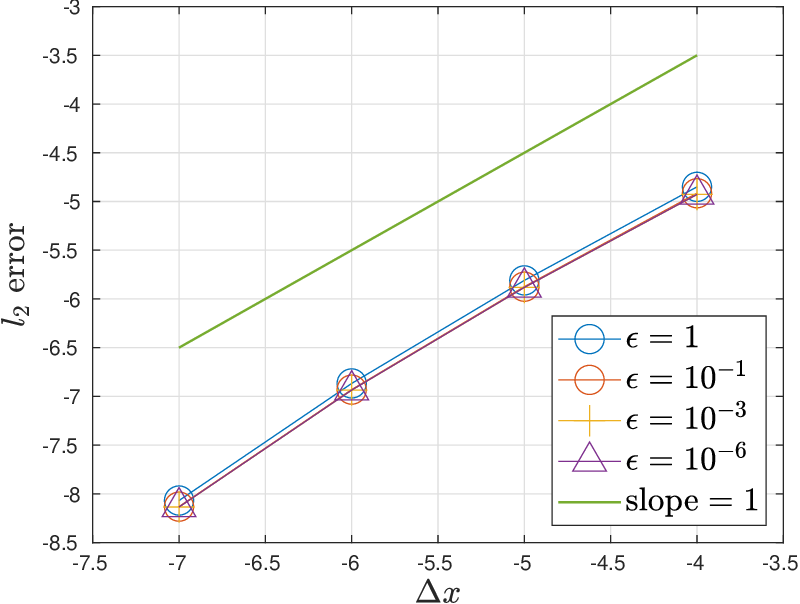}}
  \caption{(The verification of the AP property  in Sec. \ref{sec:ex1}) 
  The $l_2$ error of the numerical solution obtained by the first-order scheme \eqref{eq:first_RTE} with mesh sizes $N = 32, 64, 128$ and $256$ and the reference solution. The reference solution is obtained by the same numerical method with mesh size $N = 1024$. The parameters $\epsilon$ tested are $\epsilon = 1, 10^{-1}, 10^{-3}$ and $10^{-6}$. (a) The $l_2$ error of the radiative temperature $T_r$. (b) The $l_2$ error of the material temperature $T$. }
  \label{fig:order_error}
\end{figure}

\begin{figure}[!hptb]
  \centering
  \subfloat[radiation temperature $T_r$]{
    \includegraphics[width=0.4\textwidth]{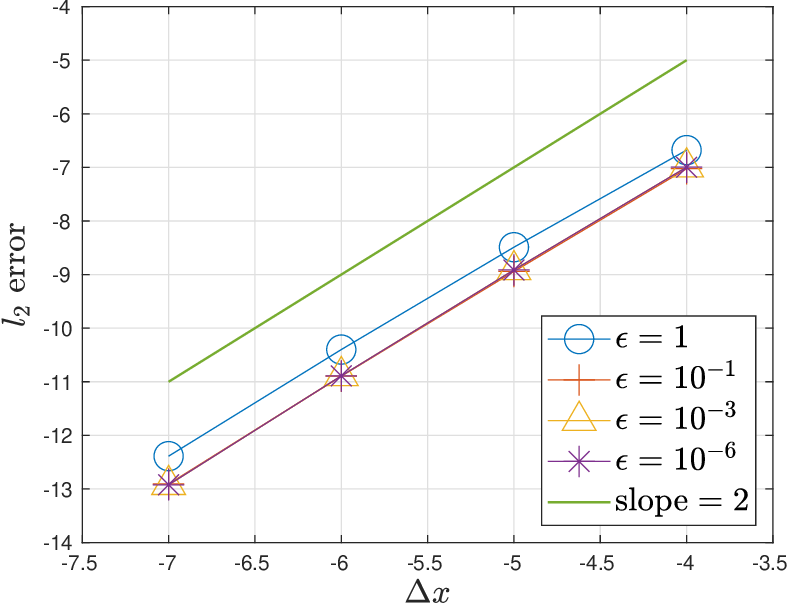}}
  \hfill
  \subfloat[material temperature $T$]{
    \includegraphics[width=0.4\textwidth]{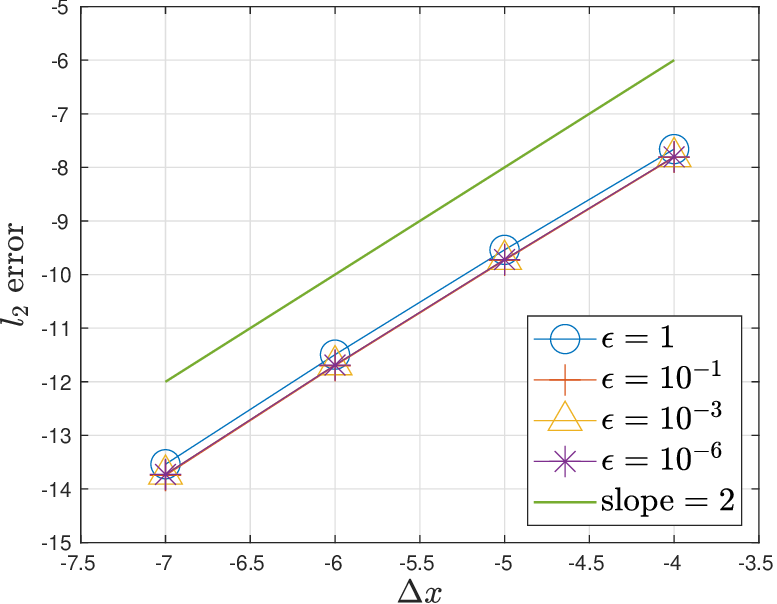}}
  \caption{(The verification of the AP property in Sec. \ref{sec:ex1}) The $l_2$ error of the numerical solution obtained by the second-order scheme \eqref{eq:second_RTE}  with mesh sizes $N = 32, 64, 128$ and $256$ and the reference solution. The reference solution is obtained by the same numerical method with mesh size $N = 1024$. The parameters $\epsilon$ tested are $\epsilon = 1, 10^{-1}, 10^{-3}$ and $10^{-6}$. (a) The $l_2$ error of the radiative temperature $T_r$. (b) The $l_2$ error of the material temperature $T$.}
  \label{fig:order_error2}
\end{figure}

In the simulation, the expansion order of the $P_N$ method is set as $M = 7$, and the mesh sizes are chosen as $N = 32, 64, 128$, and $256$ with the linear reconstruction utilized. We first consider the first-order scheme \eqref{eq:first_RTE}, the $l_2$ error of the numerical solution with different mesh sizes, and the reference solution at $t = 0.5$ is plotted in Fig. \ref{fig:order_error}, where the reference solution is obtained by the same method with the mesh size $N = 1024$. To verify the AP property of IMEX-IM, the behavior of the numerical solution with different parameter $\epsilon$ is studied, where the parameter $\epsilon$ is set as $\epsilon = 1, 10^{-1}, 10^{-3}$ and $10^{-6}$. Fig. \ref{fig:order_error} shows that for both the radiation temperature $T_r$ defined in \eqref{eq:radiation_tem} and the material temperature $T$, they are all converging to the reference solution, and the convergence rate is first order. Moreover, for different $\epsilon$, the convergence behavior is almost the same, which indicates the AP property of the first-order numerical scheme \eqref{eq:first_RTE}. For the second-order scheme \eqref{eq:second_RTE}, the numerical solution of $T_r$ and $T$ at $t = 0.5$ is illustrated in Fig. \ref{fig:order_error2}, where the mesh sizes and parameter $\epsilon$ are the same as those in Fig. \ref{fig:order_error}. The reference solution is obtained by this second-order scheme with mesh size $N = 1024$.  Fig. \ref{fig:order_error2} shows that for the second-order scheme, the numerical solution of $T_r$ and $T$ converge to the reference solution at the rate of the second order, which is consistent with the numerical scheme. Moreover, the converging behavior of the numerical solution is also the same for different $\epsilon$, indicating the AP property of the second-order scheme \eqref{eq:second_RTE}.

\subsection{1D plane source problem}
\label{sec:ex2}
In this section, the 1D plane source problem is studied, one of the classical problems for the linear radiative transport equation \eqref{eq:1D_linear_RTE}. It describes an initial particle pulse emitted in a pure scattering medium, studied in \cite{Ganapol2001, Kusch2023, PENG2021110672, PENG2020109735}. First, the Delta function is set as the equilibrium initial condition, and as stated in \cite{Kusch2023}, a Gaussian function is utilized to approximate the initial Delta function, which has the form below
\begin{equation}
    \label{eq:ex2_ini}
    I = \frac{1}{2\sqrt{2 \pi \theta}}e^{-\frac{x^2}{2 \theta}}, \qquad \theta = 1.28  \times 10^{-3}.
\end{equation}
The computational region is $L = [-6, 6]$, with the extrapolation boundary conditions enforced. The other parameters are set as 
\begin{equation}
    \label{eq:ex2_ini_1}
    \sigma_s = 1, \qquad \sigma_a = 0, \qquad G = 0, \qquad c = 1.
\end{equation}

\begin{figure}[!hptb]
\centering
  \subfloat[$t= 0$]{\includegraphics[width=0.3\textwidth]{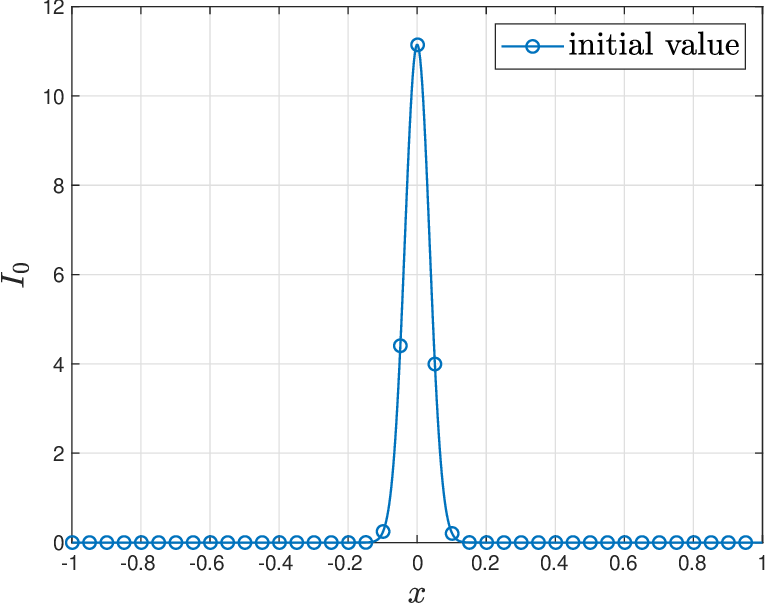}
  \label{fig:initial}}
  \hfill
  \subfloat[$t = 1$]{
    \includegraphics[width=0.3\textwidth]{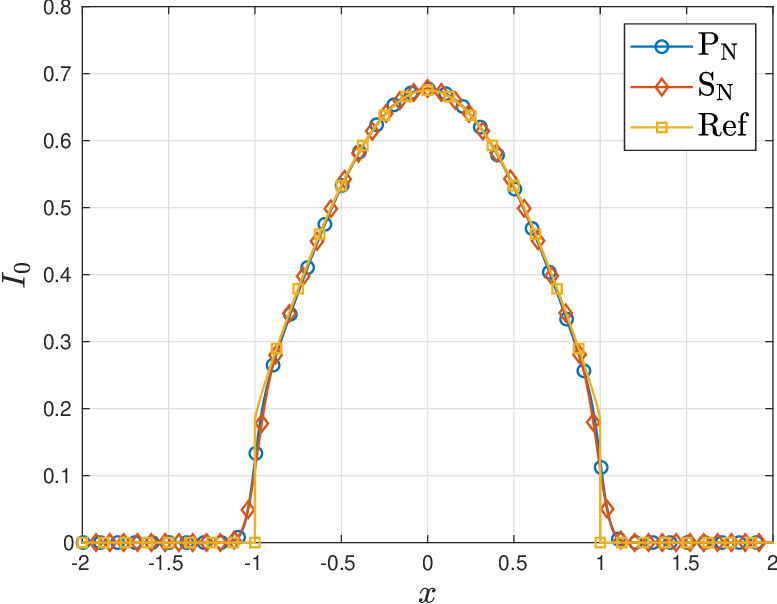}
    \label{fig:plane_source_t_1}} 
\hfill
  \subfloat[$t = 5$]{
  \label{fig:plane_source_t_5}
    \includegraphics[width=0.3\textwidth]{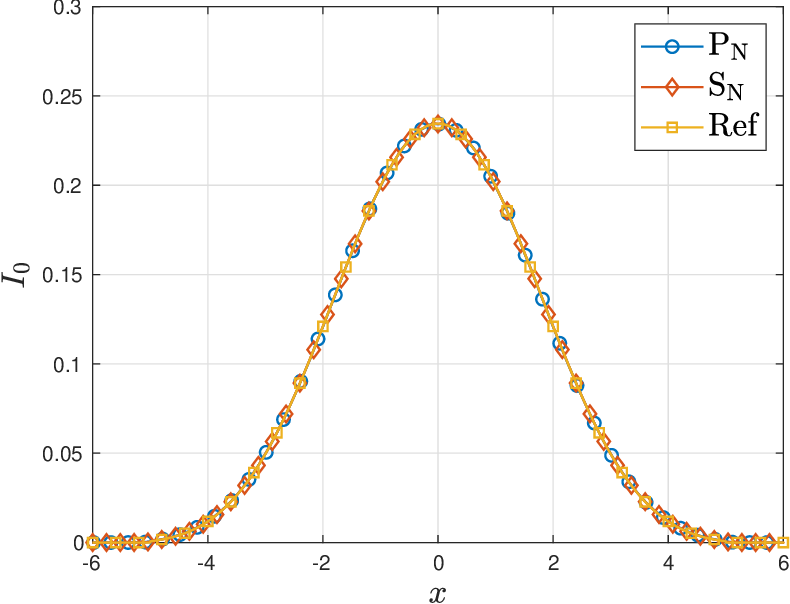}}    
  \caption{(1D plane source problem in Sec. \ref{sec:ex2}) $I_{0}$ of the plane source problem for different times. Here, the blue circle lines are obtained from IMEX-IM, the orange diamond lines are the reference solution obtained by the $S_N$ method, and the yellow square lines are the analytical solution in \cite{Ganapol2001}. (a) Initial condition $t = 0$. (b) $t = 1$. (c) $t = 5$.}
  \label{fig:plane_source}
\end{figure}

In the simulation, the expansion order of the $P_N$ method is set as $M =  59$. Here, we want to emphasize that when solving the plane source problem with the $P_N$, the non-physical oscillation affect the numerical solution \cite{PENG2021110672, PENG2020109735}. There are several methods proposed to mitigate this phenomenon \cite{Olson2009, mcclarren2010robust}. However, we only utilize a large expansion number to undermine the non-physical oscillation here. We first set $\epsilon = 1$, and the mesh size as $N = 1200$ with the linear reconstruction utilized. The numerical solution of $I_0$ at $t = 1$ is shown in Fig. \ref{fig:plane_source_t_1}, as well as the reference solution obtained by the $S_N$ method, and the analytical benchmark solution given in \cite{Ganapol2001}. It shows that the three solutions fit well with each other. Moreover, the support region of the numerical solution is $[-1.5, 1.5]$. Compared to the initial condition plotted in Fig. \ref{fig:initial}, which is an approximation to the Delta function, $I_0$ keeps diffusing. The numerical solution of $I_0$ at $t = 5$ is shown in Fig. \ref{fig:plane_source_t_5}, where it also matches well with the reference solution obtained by the $S_N$ method, and the analytical solution by \cite{Ganapol2001}. Compared to the initial condition shown in Fig. \ref{fig:initial} and the numerical solution at $t=1$, the diffusion property of the numerical solution is shown more clearly. 

\begin{figure}[!hptb]
  \centering
   \subfloat[$\epsilon = 1$]{\includegraphics[width=0.3\textwidth, height = 0.24\textwidth]{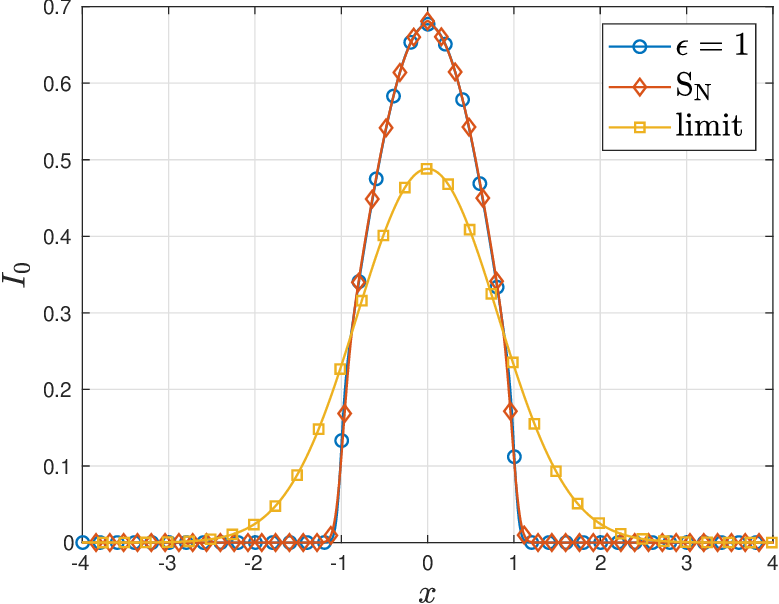}
  \label{fig:planesource_ep_1}}
  \hfill
  \subfloat[$\epsilon = 0.8$]{\includegraphics[width=0.3\textwidth, height = 0.24\textwidth]{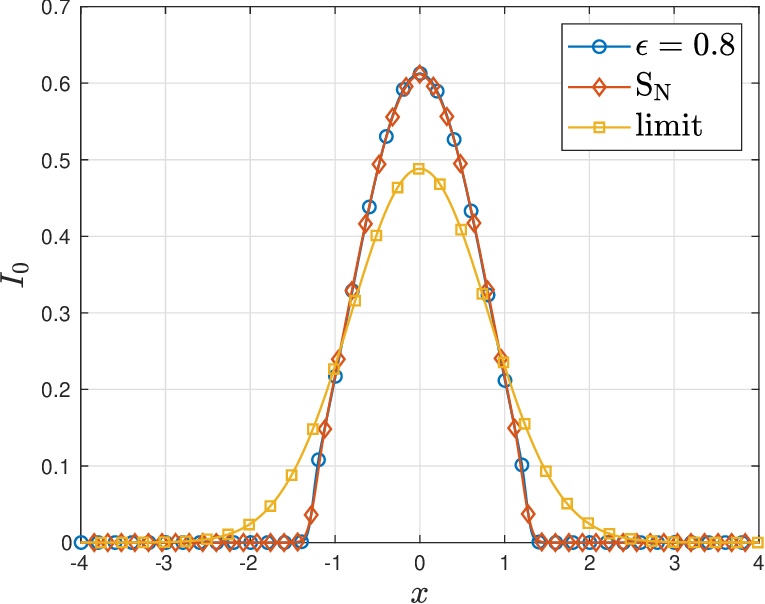}
  \label{fig:planesource_ep_08}}
  \hfill
  \subfloat[$\epsilon = 0.6$]{\includegraphics[width=0.3\textwidth, height = 0.24\textwidth]{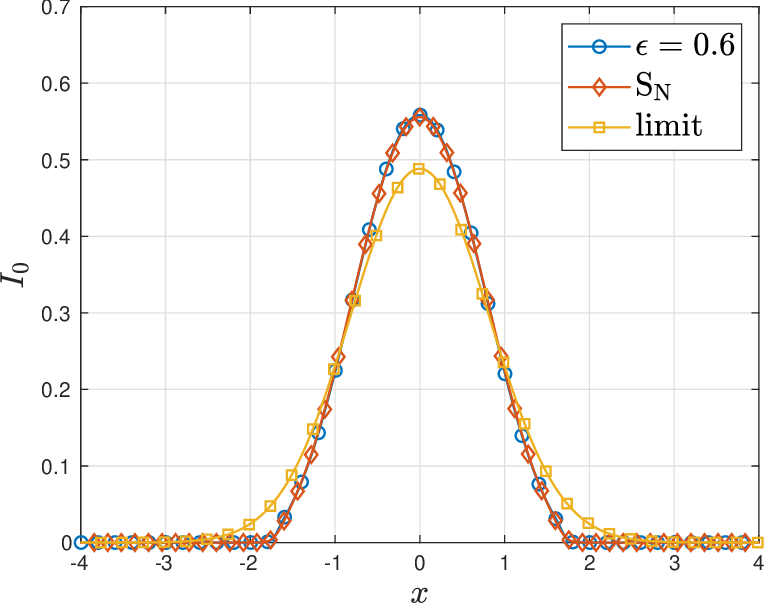}
  \label{fig:planesource_ep_06}
  }\hfill
  \subfloat[$\epsilon = 0.4$]{\includegraphics[width=0.3\textwidth, height = 0.24\textwidth]{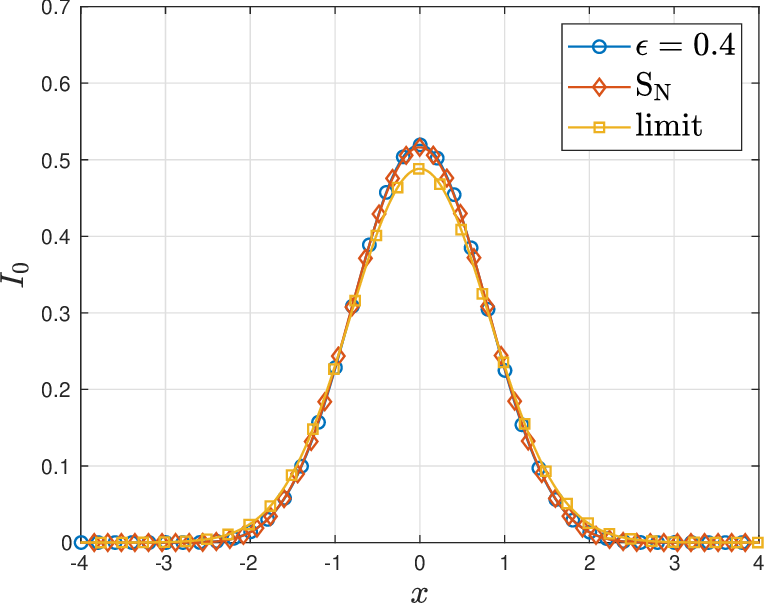}
  \label{fig:planesource_ep_04}}
  \hfill
  \subfloat[$\epsilon = 10^{-3}$]{\includegraphics[width=0.3\textwidth, height = 0.24\textwidth]{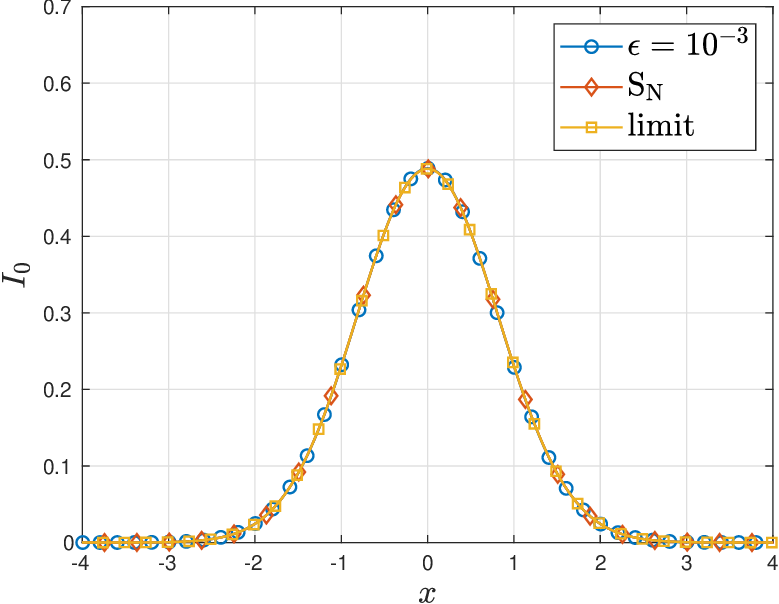}
  \label{fig:planesource_ep_02}}
  \hfill
  \subfloat[$\epsilon = 10^{-6}$]{\includegraphics[width=0.3\textwidth, height = 0.24\textwidth]{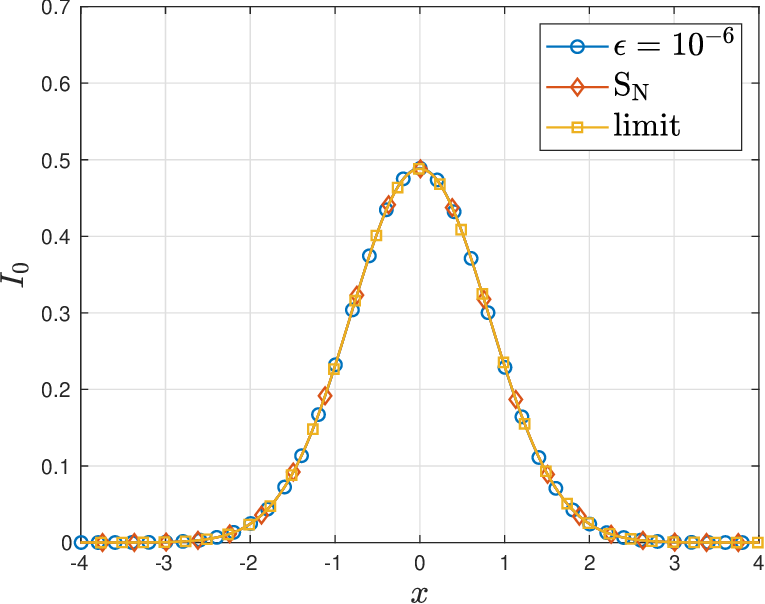}
  \label{fig:planesource_ep_1e6}}
  \caption{(1D plane source problem in Sec. \ref{sec:ex2}) $I_{0}$ of the plane source problem for different $\epsilon$. Here, the blue circle lines are obtained from IMEX-IM, the orange diamond lines are the reference solution obtained by the $S_N$ method, and the yellow square lines are the numerical solution of the diffusion limit equation \eqref{eq:1D_linear_limit}. (a) $\epsilon = 1$. (b) $\epsilon = 0.8$. (c) $\epsilon = 0.6$. (d) $\epsilon = 0.4$. (e) $\epsilon = 10^{-3}$. (f) $\epsilon = 10^{-6}$. }
  \label{fig:plane_source_limit}
\end{figure}

To verify the AP property of IMEX-IM, the behavior of $I_0$ with different parameter $\epsilon$ is studied, with the other parameters such as the expansion order $M$ and the time step length $\Delta t$ remaining the same. Fig. \ref{fig:plane_source_limit} shows the numerical solutions of $I_0$ for $\epsilon = 1, 0.8, 0.6, 0.4, 10^{-3}$ and $10^{-6}$ at $t = 1$, as well as the reference solutions obtained by the $S_N$ method, and the numerical solutions of the diffusion limit equation \eqref{eq:1D_linear_limit}. It illustrates that for different $\epsilon$, the numerical solutions of $I_0$ all match well with the reference solutions. Moreover, the smaller $\epsilon$ is, the quicker the diffusion of $I_0$ is. Besides, when $\epsilon$ is large as $\epsilon = 1$, there exists a large difference between the numerical solution of RTE \eqref{eq:1D_linear_RTE} and the diffusion limit equation \eqref{eq:1D_linear_limit}. With the decreasing of $\epsilon$, the discrepancy between these two solutions is getting smaller. When $\epsilon$ is decreasing to $\epsilon = 10^{-3}$, the numerical solution of RTE \eqref{eq:1D_linear_RTE} and that of the diffusion limit equation \eqref{eq:1D_linear_limit} are almost on top of each other. We want to emphasize that the time step length even for $\epsilon = 10^{-6}$ is the same as that for $\epsilon = 1$, which also indicates the AP property of IMEX-IM.

\subsection{1D Marshak wave problems}
\label{sec:ex3}

In this section, the classical Marshak wave problem is studied. It consists of two examples with the absorption/emission coefficients $\sigma$ 
depending on the temperature $T$. This problem is also studied in \cite{sun2015asymptotic1, Larsen2013, semi2008Ryan}. The related parameters are set as the speed of light $ c = 29.98 \rm{cm/ns}$, the radiation constant $a = 0.01372  \rm{GJ/cm^3}-\rm{keV^4} $, and the specific heat $0.1 \rm {GJ/g/keV}$. The initial temperature $T$ is chosen as $T = 10^{-6} {\rm keV}$, with the radiation density at equilibrium 
\begin{equation}
    \label{eq:ex3_ini} 
    I = \frac{1}{2}a c T^4.
\end{equation}
A constant isotropic incident radiation intensity with a Planckian distribution as  $T = 1 \rm{keV}$ is set at the left boundary, where the Marshak type boundary condition \cite{semi2008Ryan} is utilized.

\paragraph{Marshak wave-2B}

\begin{figure}[!hptb]
  \centering 
    \subfloat[Radiation temperature $T_r$]{
    \includegraphics[width=0.4\textwidth]{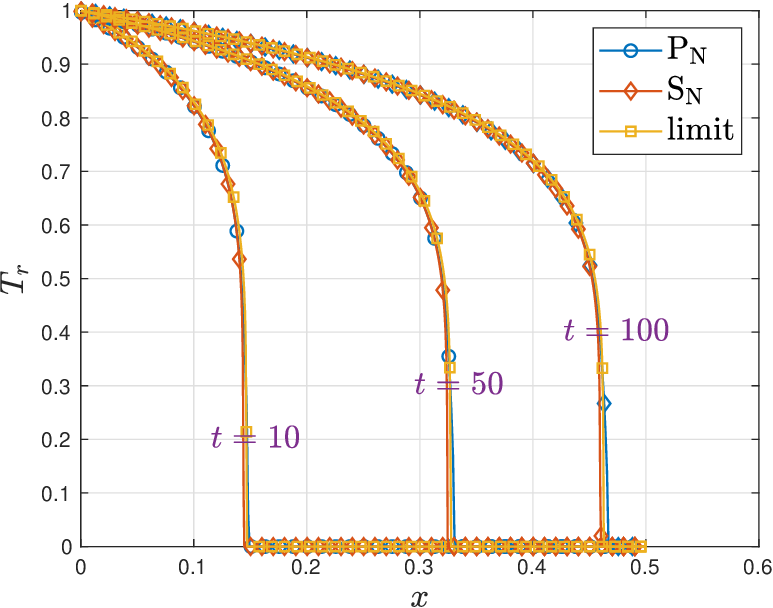}
    \label{fig:2B_Tr}} \hfill
   \subfloat[Material temperature $T$]{
    \includegraphics[width=0.4\textwidth]{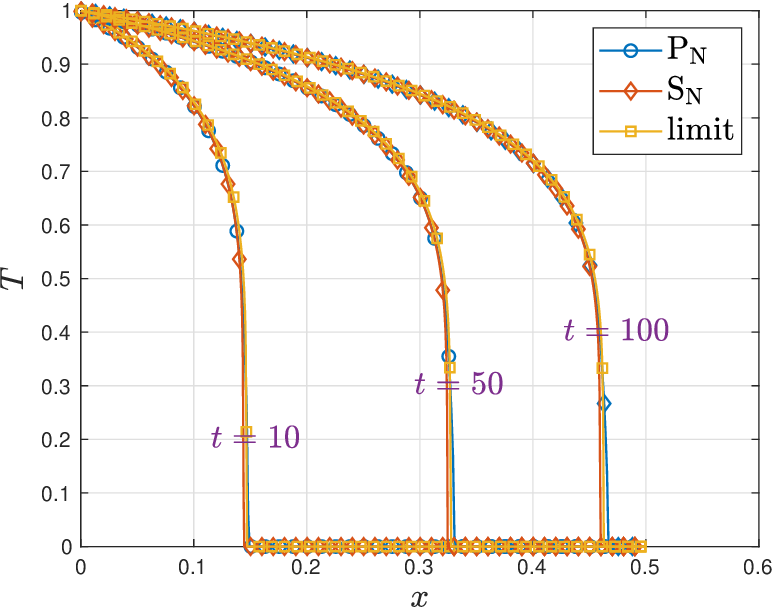}
    \label{fig:2B_1}}   
  \caption{(Marshak wave-2B problem in Sec. \ref{sec:ex3}) The radiation temperature $T_r$ and material temperature $T$ of the Marshak wave-2B problem at times $t = 10, 50$ and $100$. The blue circle lines are the numerical solution obtained by IMEX-IM, the orange diamond lines are the reference solution obtained by the $S_N$ method, and the yellow square lines are the solution of the diffusion limit equation. (a) Radiation temperature $T_r$. (b) Material temperature $T$. }
  \label{fig:2B}
\end{figure}

We first set the absorption/emission coefficient as is $\sigma = 100 \rho/T^3 {\rm cm^2/g}$ and the density $\rho$ as $\rho = 3.0 \rm {g/cm^3}$. In the simulation, the computational region is taken as $L = [0, 0.6]$, and the expansion order of the $P_N$ method is set as $M = 7$, with the spatial grid size as $N = 400$. In the Marshak wave-2B problem, since the absorption/emission coefficient is large enough, the numerical solution of RTE \eqref{eq:1D_RTE} is almost the same as that of the diffusion limit \eqref{eq:1D_limit}.

Fig. \ref{fig:2B} presents the numerical solution of the radiation temperature $T_r$ and the material temperature $T$ at times $t = 10, 50$ and $100$. The radiation wave propagates forward at a certain speed. The reference solution obtained by the $S_N$ method and the numerical solution to the diffusion limit \eqref{eq:1D_limit} are also plotted. It shows that for both temperatures, the numerical solution matches well with both the reference solution and the solution to diffusion limit. This is also consistent with the theoretical analysis.

\paragraph{Marshak wave-2A} 
In the Marshak wave-2A problem, a smaller absorption/emission coefficient is utilized as $\sigma =10/ T^3 {\rm cm^2/g}$ and the density $\rho = 3.0 \rm {g/cm^3}$. We first set $\epsilon = 1$. Since the absorption/emission $\sigma$ is not large enough, the solution to RTE is not consistent with that of the diffusion limit. The same expansion order and the grid size as the Marshak wave-2B problem are adopted as $M = 7$ and $N = 400$. 
\begin{figure}[!hptb]
  \centering
   \subfloat[Radiation temperature $T_r$]{
       \includegraphics[width=0.4\textwidth]{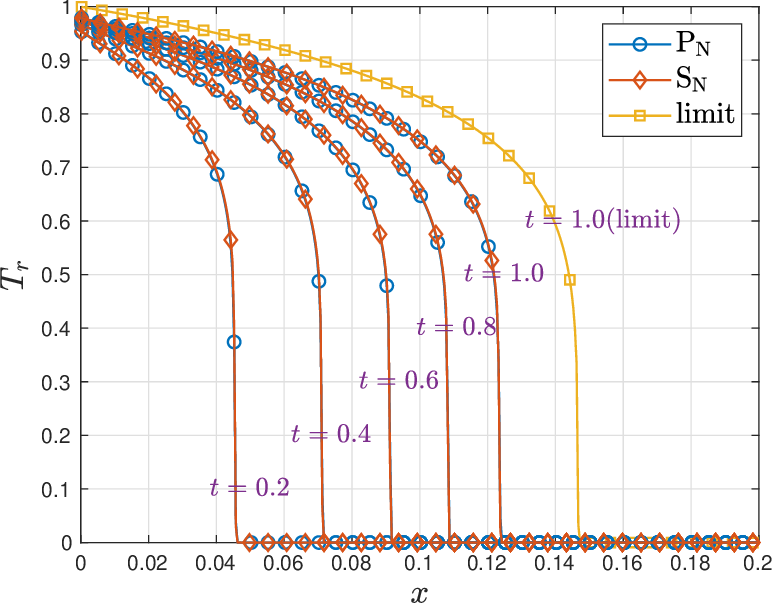}} 
       \hfill
  \subfloat[Material temperature $T$]{
       \includegraphics[width=0.4\textwidth]{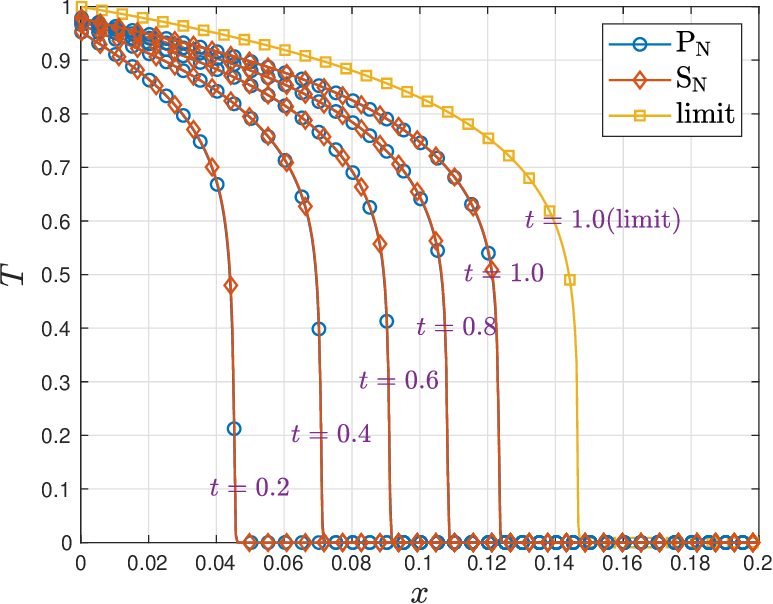}}   
   \caption{(Marshak wave-2A problem in Sec. \ref{sec:ex3}) 
   The radiation temperature $T_r$ and material temperature $T$ of the Marshak wave-2A problem with $\epsilon = 1$
   at times $t = 0.2, 0.4, 0.6, 0.8$ and $1.0$. The blue circle lines are the numerical solution obtained by IMEX-IM,
    the orange diamond lines are the reference solution obtained by the $S_N$ method, and the yellow square lines are the solution of the diffusion limit equation at $t = 1$. (a) Radiation temperature $T_r$. (b) Material temperature $T$.}
 \label{fig:2A}
\end{figure}

 Fig. \ref{fig:2A} illustrates the numerical solution of the radiation temperature $T_r$ and the material temperature $T$ at times $t = 0.2, 0.4, 0.6, 0.8$ and $1$ as well as the reference solution obtained by the $S_N$ method, and the numerical solution to the diffusion limit at $t = 1$. The numerical solution fits well with the reference solution by the $S_N$ method. Compared to the Marshak wave-2B problem, the radiation wave in the Marshak wave-2A problem propagates much more slowly. Moreover, the radiation wave by the diffusion limit moves much faster than that of the Marshak wave-2A problem.

\begin{figure}[!hptb]
  \centering
  \subfloat[Radiation temperature $T_r$]{
     \includegraphics[width=0.4\textwidth]{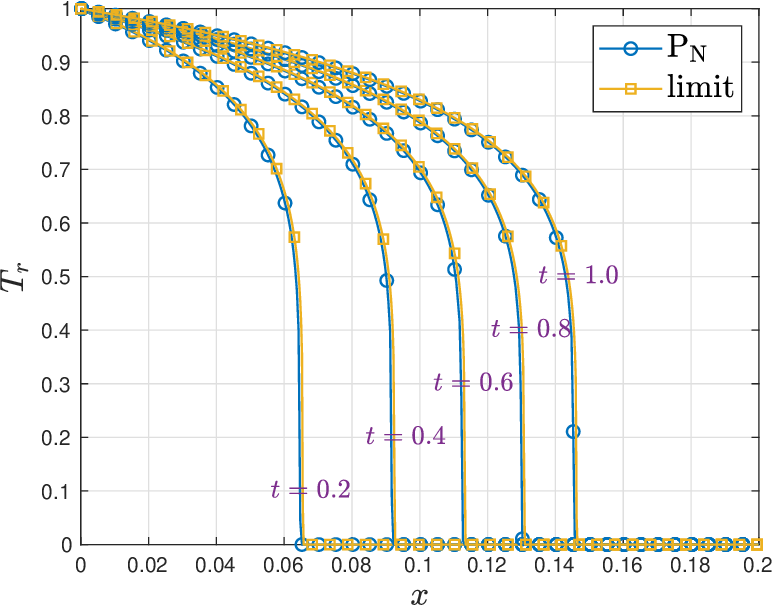}
     } \hfill
  \subfloat[Material temperature $T$]{    
    \includegraphics[width=0.4\textwidth]{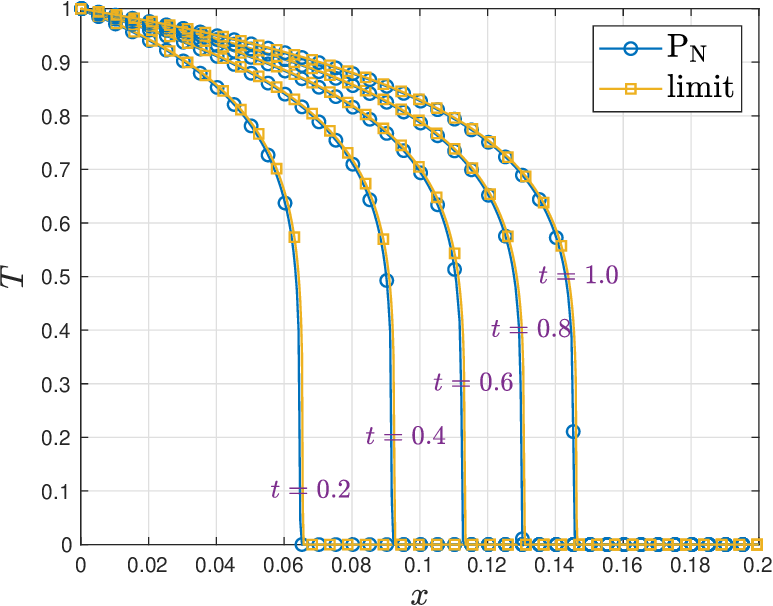}} 
    \caption{(Marshak wave-2A problem in Sec. \ref{sec:ex3}) The radiation temperature $T_r$ and material temperature $T$ of the Marshak wave-2A problem with $\epsilon = 10^{-6}$ at times $t = 0.2, 0.4, 0.6, 0.8$ and $1.0$. The blue circle lines are the numerical solution obtained by IMEX-IM, and the yellow square lines are the solution of the diffusion limit equation. (a) Radiation temperature $T_r$. (b) Material temperature $T$.}
\label{fig:2A_limit}
\end{figure}

To validate the AP property of IMEX-IM, we set $\epsilon = 10^{-6}$, and other parameters remain the same. The numerical solutions of $T_r$ and $T$ for RTE and the diffusion limit at time $t = 0.2, 0.4, 0.6, 0.8$ and $1$ are plotted in Fig. \ref{fig:2A_limit}. It shows that the numerical solution of RTE is almost the same as that of the diffusion limit, which also validates the AP property of the new numerical scheme  \eqref{eq:first_RTE_non} and \eqref{eq:convection}.

\subsection{2D line source problem}
\label{sec:ex4}
In this section, the 2D line source problem is studied, another benchmark problem first proposed by Ganapol \cite{Ganapol1999}, and is also studied in \cite{Hu2021lowrank, Comparison2013}. This line source problem is the 2D version of the 1D plane source problem in Sec. \ref{sec:ex2}, but it is more difficult due to the curse of dimensionality. Here, a Gaussian function is adopted to approximate the initial Delta function as 
\begin{equation}
    \label{ex4:ini}
 I =  \frac{1}{4\pi}\left(\frac{1}{4\pi\theta}\exp\left(-\frac{x^2 + y^2}{4 \theta}\right)\right), \qquad \theta = 1.6 \times 10^{-4}
\end{equation}
with the other parameters the same as in \eqref{eq:ex2_ini_1}.

\begin{figure}[!hptb]
  \centering
  \subfloat[$E_r$ by 1st-order scheme]{
    \label{fig:linesource_1st_ep1}
    \includegraphics[width=0.3\textwidth]{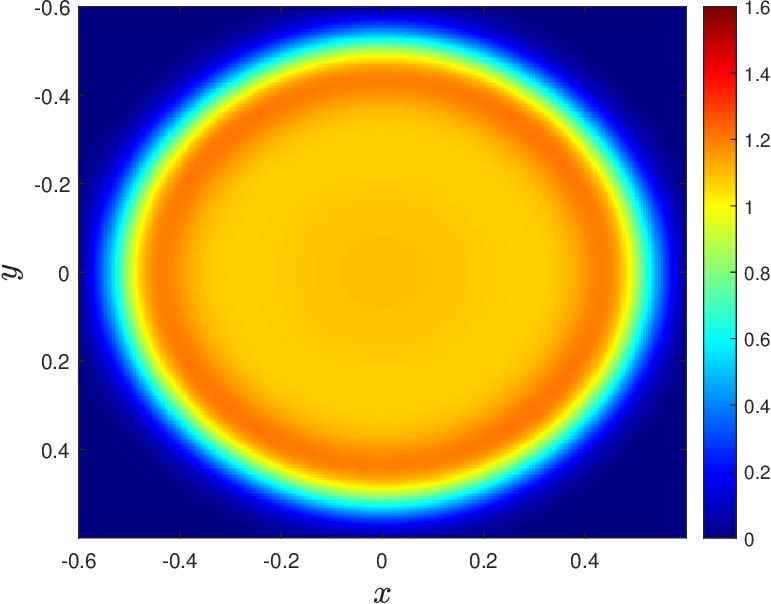}}
  \hfill 
  \subfloat[$E_r$ by 2nd-order scheme]{
    \label{fig:linesource_2nd_ep1}
    \includegraphics[width=0.3\textwidth]{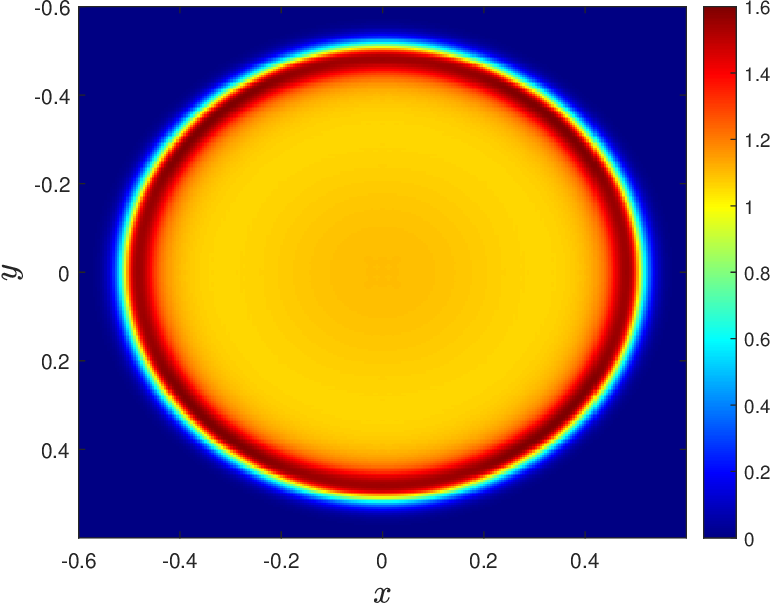}}
  \hfill 
  \subfloat[ $E_r$ of reference solution]{
    \label{fig:linesource_cut}
    \includegraphics[width=0.3\textwidth]{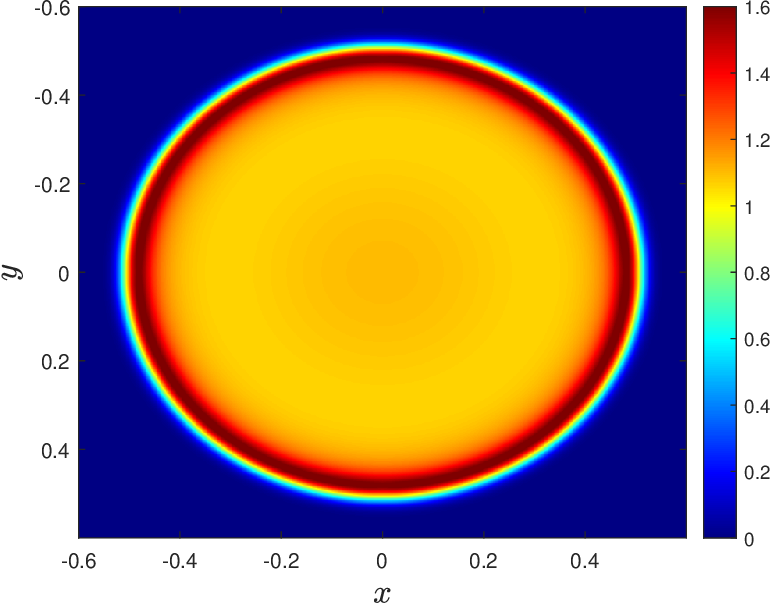}}
  \caption{(Line source problem in Sec. \ref{sec:ex4}) The contour plot of the numerical solution $E_r$ at $t = 0.5$ with $\epsilon = 1$. Here, the mesh size is $N_x = N_y = 256$. (a) Contour plot of $E_r$ by the 1st-order scheme. (b) Contour plot of $E_r$ by 2nd-order scheme. (c) Contour plot of $E_r$ of the reference solution obtained by the semi-analytical solution \cite{Ganapol2001, Comparison2013}.}
    \label{fig:linesource_ep_1}
\end{figure}

\begin{figure}[!hptb]
  \centering
    \subfloat[ $E_r$ when $\epsilon = 1$]{
    \label{fig:linesource_t_05_CFL_04_ep_1_x_0_line_2}
    \includegraphics[width=0.45\textwidth]{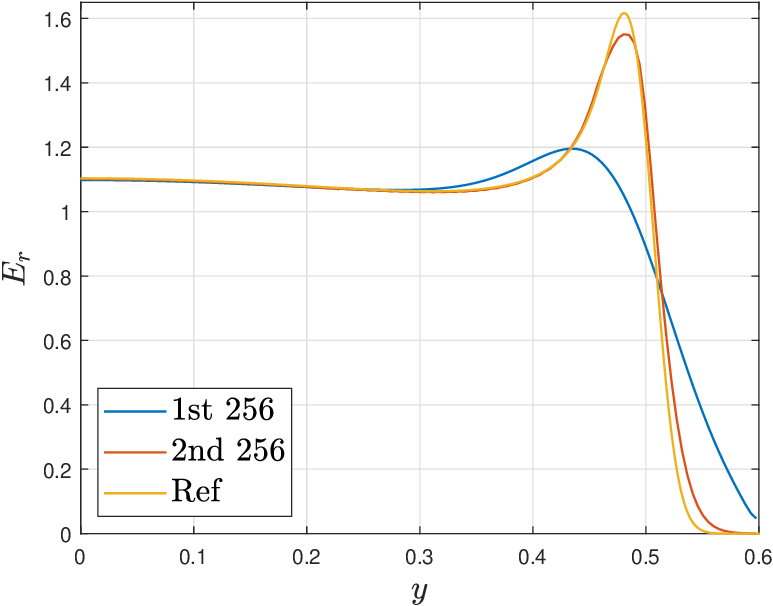}}
    \hfill
  \subfloat[Convergence behavior of $E_r$ when $\epsilon = 1$]{
    \label{fig:linesource_t_05_CFL_04_ep_1_x_star_6}
    \includegraphics[width=0.45\textwidth]{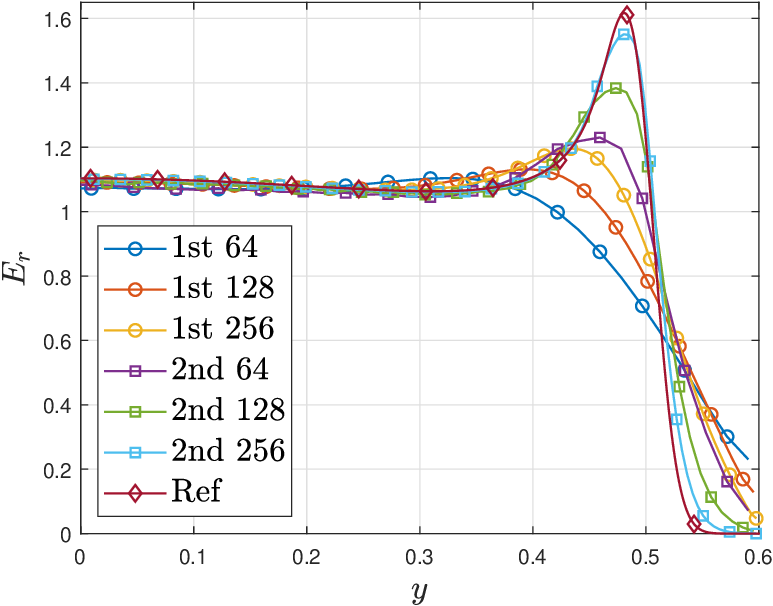}}
    \caption{(Line source problem in Sec. \ref{sec:ex4}) The numerical solution $E_r$ at $t = 0.5$ along $x = 0$ when $\epsilon = 1$.  (a) The comparison of the numerical solution by first and second-order schemes, and the reference solution obtained by the semi-analytical solution. Here, the mesh size is $N_x = N_y = 256$. (b) The convergence behavior of $E_r$. The mesh sizes are set as $N_x = N_y = 64, 128, 256$. }
    \label{fig:linesource_cut_ep_1}
\end{figure}

In the simulation, we first set $\epsilon = 1$, and the computational domain is chosen as $[-0.6, 0.6] \times [-0.6, 0.6]$. The expansion order of the $P_N$ method is set as $M = 59$. To validate the high efficiency of the second-order scheme \eqref{eq:first_order_scheme_linear_time_2_order}, the numerical results of the radiation density $E_r$ \eqref{eq:energy_Er} by the first-order scheme \eqref{eq:first_linear_I} and \eqref{eq:first_linear_II} and the second-order scheme \eqref{eq:first_order_scheme_linear_time_2_order} with mesh size $N_x = N_y = 256$ at $t = 0.5$ are plotted in Fig. \ref{fig:linesource_ep_1}, where the semi-analytical solution \cite{Ganapol2001, Comparison2013} is also plotted as the reference solution. It shows clearly the numerical solution matches well with the reference solution and the second-order numerical solution indicates high resolution compared to the first-order scheme. Fig. \ref{fig:linesource_ep_1} also indicates that IMEX-IM keeps the rotation invariant, and there is no non-physical periodic phenomenon such as the ray effect \cite{Frank2020Rayeffectartifical2020} appearing. To compare the numerical solution and the reference more clearly, the numerical solution $E_r$ along $x= 0$ is plotted in Fig. \ref{fig:linesource_cut_ep_1}, where Fig. \ref{fig:linesource_t_05_CFL_04_ep_1_x_0_line_2} shows that the second-order numerical solution almost matches the reference solution, while there is a large difference between the first-order numerical solution and the reference solution. The convergence of the numerical solution is illustrated in Fig. \ref{fig:linesource_t_05_CFL_04_ep_1_x_star_6}, where the second-order numerical solution also behaves much better than the first-order scheme.

% where different mesh number $N = N_x = N_y = 64, 128$ and $256$ are utilized. For the first-order scheme, there is no spatial reconstruction, while for the second-order scheme, linear reconstruction is applied for the expansion coefficients for the first two orders, where the linear equation system is solved to update the first two order expansion coefficients, and the fifth-order WENO reconstruction is applied to the other expansion coefficients. The semi-analytical solution \cite{Ganapol2001, Comparison2013} and the numerical solution by $S_N$ method are taken as the reference solution. 

% \wyl{add SN} It indicates that this numerical scheme keeps the rotation invariant, and there is no non-physical periodic phenomenon such as the ray effect \cite{Frank2020Rayeffectartifical2020} appearing. Moreover, the numerical solution together with both reference solutions along $x = 0$ is plotted in Fig. \ref{fig:linesource_cut_ep_1}, where the numerical solution of $E_r$ obtained by the first and second-order scheme with different mesh sizes are shown in Fig. \ref{fig:linesource_t_05_CFL_04_ep_1_x_star_6}. It illustrates that the numerical solution is converging to the reference solution with the increasing mesh size, and the high efficiency of the second-order scheme is clearly seen in Fig. \ref{fig:linesource_t_05_CFL_04_ep_1_x_0_line_2}. 

\begin{figure}[!hptb]
  \centering
  \subfloat[$E_r$ by 1st-order scheme]{
    \label{fig:linesource_1st_ep1e1}
    \includegraphics[width=0.3\textwidth]{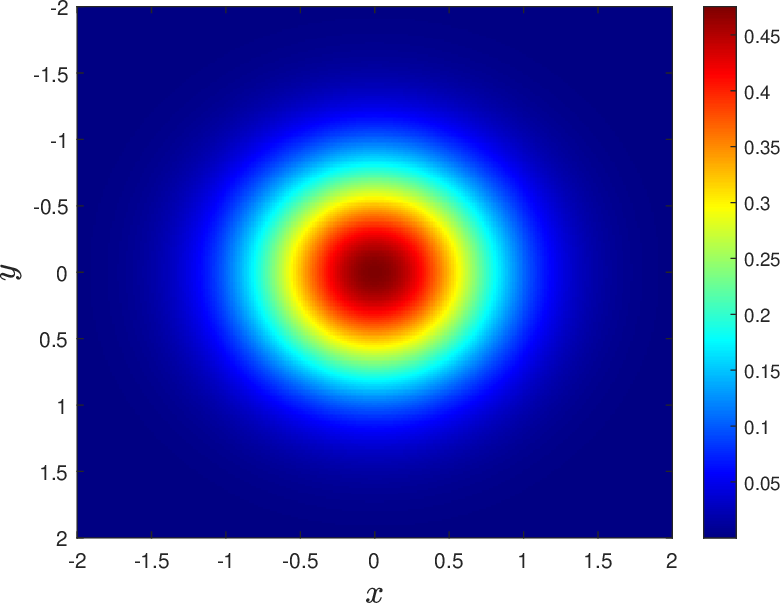}}
  \hfill 
  \subfloat[$E_r$ by 2nd-order scheme]{
    \label{fig:linesource_2nd_ep1e1}
    \includegraphics[width=0.3\textwidth]{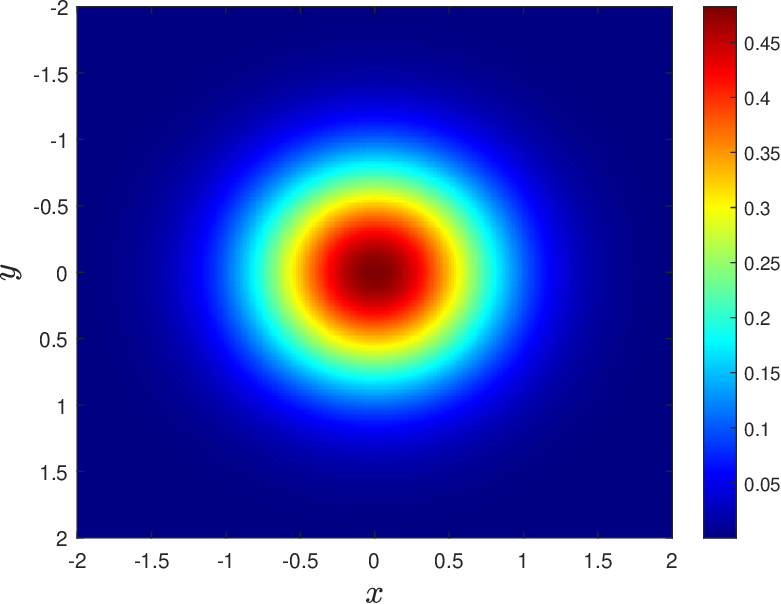}}
\hfill
    \subfloat[$E_r$ by $S_N$ method]{
    \label{fig:linesource_Sn_ep1e1}
    \includegraphics[width=0.3\textwidth]{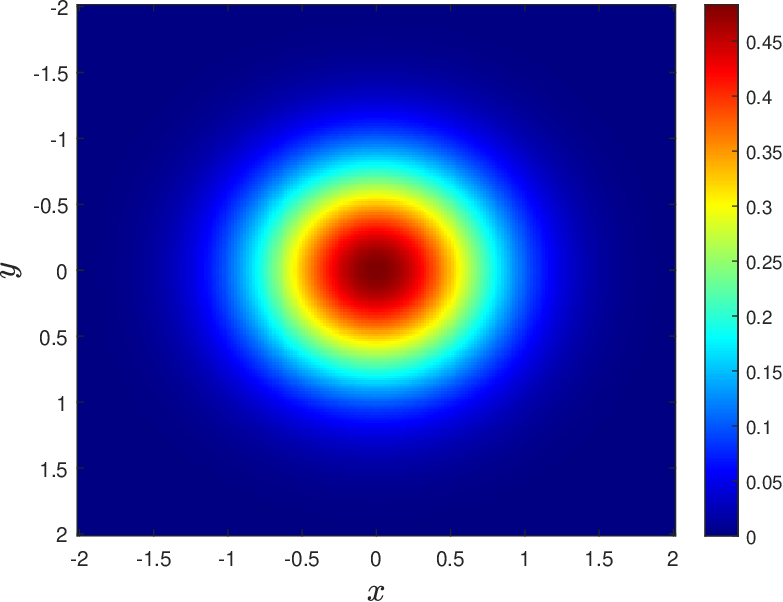}}
  \caption{(Line source problem in Sec. \ref{sec:ex4}) The contour plot of the numerical solution $E_r$ at $t = 0.5$ with $\epsilon = 0.1$. Here, the mesh size is $N_x = N_y = 256$. (a) Contour plot of $E_r$ by the 1st-order scheme. (b) Contour plot of $E_r$ by the 2nd-order scheme. (c) Contour plot of $E_r$ of the reference solution obtained by the $S_N$ method.}
  \label{fig:linesource_limit_ep_1e1}
\end{figure}

\begin{figure}[!hptb]
  \centering
    \subfloat[$E_r$ when $\epsilon = 0.1$]{
    \label{fig:linesource_t_05__CFL_04_ep_1e1_x_0_line_2}
    \includegraphics[width=0.45\textwidth]{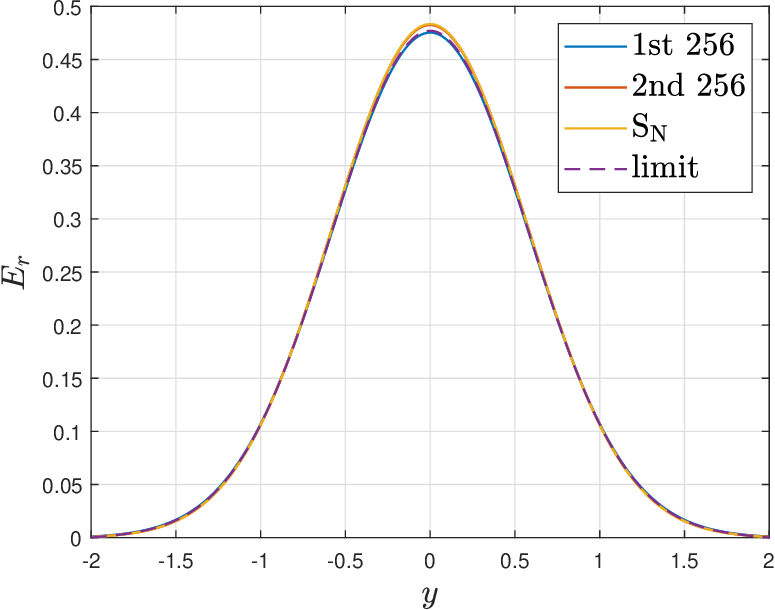}}
  \hfill  
   \subfloat[$E_r$  when $\epsilon = 0.1$]{
    \label{fig:linesource_t_05__CFL_04_ep_1e1_x_star_6}
    \includegraphics[width=0.505\textwidth]{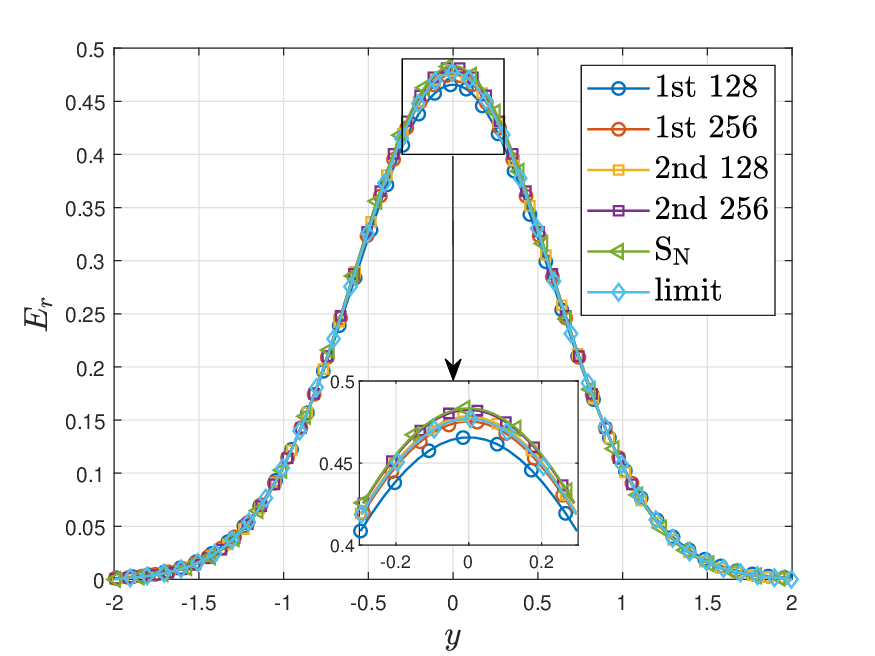}}
    \caption{(Line source problem in Sec. \ref{sec:ex4})
    The numerical solution $E_r$ at $t = 0.5$ along $x = 0$ when $\epsilon = 0.1$. (a) The comparison of the numerical solution by first and second-order schemes, and the reference solution obtained by the $S_N$ method. Here, the mesh size is $N_x = N_y = 256$. (b) The convergence behavior of $E_r$. Here, the mesh sizes are chosen as $N_x = N_y = 128, 256$.}
    \label{fig:linesource_cut_limit_ep_1e1}
\end{figure}

To validate the AP property of IMEX-IM, we first follow the tests in \cite{Hu2021lowrank} and set $\epsilon = 0.1$. In this case, the computational region is $[-2, 2]\times [-2, 2]$, and the expansion order of the $P_N$ method is set as $M = 7$, while the other parameters remain the same as $\epsilon = 1$. The numerical solution of $E_r$ by the first- and second-order schemes with mesh size $N_x = N_y = 256$ at $t = 0.5$ is shown in Fig. \ref{fig:linesource_limit_ep_1e1}, where the reference solution obtained by the $S_N$ method is also plotted. In the $S_N$ method, $30$ Gauss quadrature points in $[-1,1]$ and $60$ equally spaced points in $[0,2\pi]$ are utilized in velocity space with a mesh size of $200 \times 200$ used in the spatial space. With this small $\epsilon$, IMEX-IM still keeps the rotation invariance well. The numerical solution of $E_r$, together with the reference solution obtained by the $S_N$ method and the numerical solution of the diffusion limit equation \eqref{eq:linear_limit} is shown in Fig. \ref{fig:linesource_t_05__CFL_04_ep_1e1_x_0_line_2}, where the numerical solution obtained by the second-order scheme matches well with the reference solution, while there is a small error between the numerical solution when $\epsilon = 0.1$ and the numerical solution of the diffusion limit equation \eqref{eq:linear_limit}. Moreover, the convergence of the numerical solution by the first and second-order schemes is plotted in Fig. \ref{fig:linesource_t_05__CFL_04_ep_1e1_x_star_6}. When $\epsilon$ approaches zero, the first-order scheme will converge to the second-order five-point difference scheme of the diffusion limit equation, which is the reason why the behavior of the first and second-order schemes is almost the same when $\epsilon$ is small. Though these numerical solutions are almost on top of each other, we still find that the second-order scheme behaves better than the first-order scheme.

\begin{figure}[!hptb]
  \centering
  \subfloat[$E_r$ by 1st-order scheme]{
    \label{fig:linesource_t_05_CFL_04_ep_1e6_1order_imagesc}
    \includegraphics[width=0.3\textwidth]{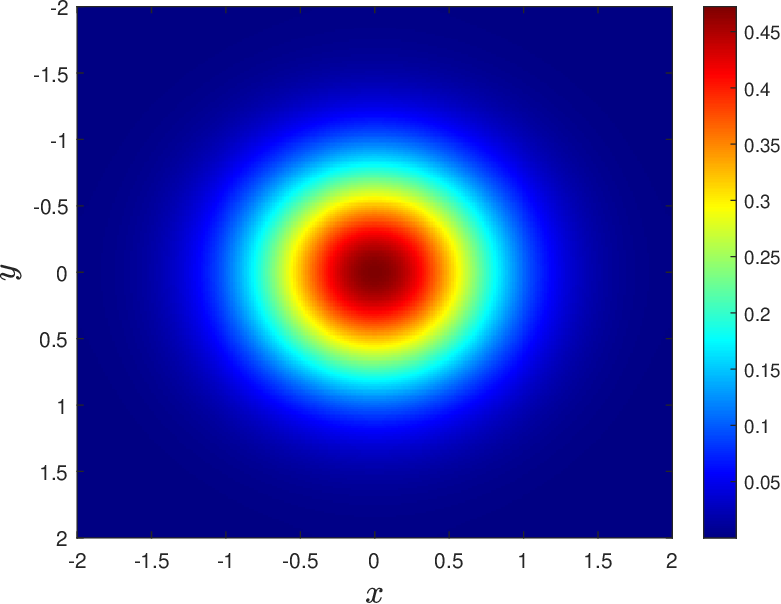}}
  \hfill 
  \subfloat[$E_r$ by 2nd-order scheme]{
    \label{fig:linesource_t_05__CFL_04_ep_1e6_2order_imagesc}
    \includegraphics[width=0.3\textwidth]{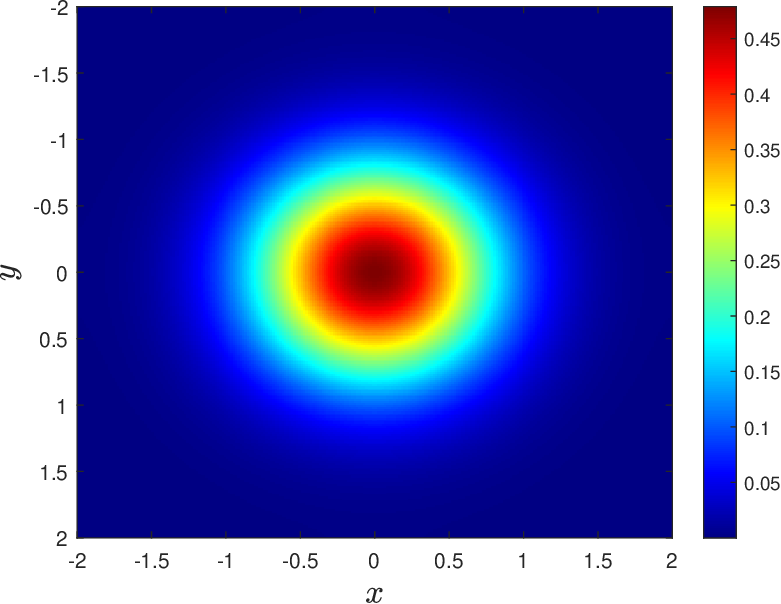}}
     \hfill 
  \subfloat[$E_r$ of the diffusion limit equation]{
    \label{fig:linesource_t_05_CFL_04_ep_limit_imagesc}
    \includegraphics[width=0.3\textwidth]{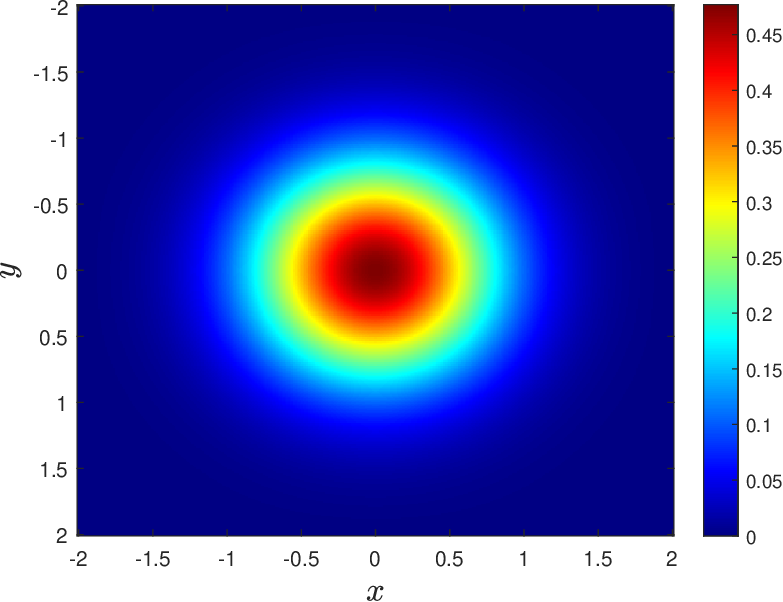}}
  \caption{(Line source problem in Sec. \ref{sec:ex4}) The contour plot of the numerical solution $E_r$ at $t = 0.5$ with $\epsilon = 10^{-6}$. Here, the mesh size is $N_x = N_y = 256$. (a) Contour plot of $E_r$ by the 1st-order scheme. (b) Contour plot of $E_r$ by the 2nd-order scheme. (c) Contour plot of $E_r$ of the diffusion limit equation. }
  \label{fig:linesource_limit_ep_1e6}
  \end{figure}

\begin{figure}[!hptb]
  \centering
    \subfloat[$E_r$ when $\epsilon = 10^{-6}$]{
    \label{fig:linesource_t_05__CFL_04_ep_1e6_x_0_line_2}
    \includegraphics[width=0.45\textwidth]{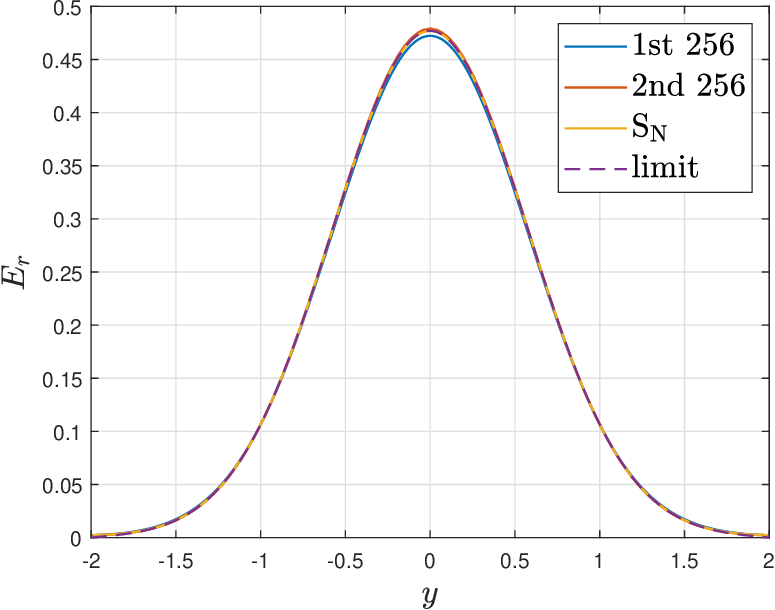}}
  \hfill 
   \subfloat[Convergence behavior of  $E_r$ when $\epsilon = 10^{-6}$]{
    \label{fig:linesource_t_05__CFL_04_ep_1e6_x_star_6}
    \includegraphics[width=0.505\textwidth]{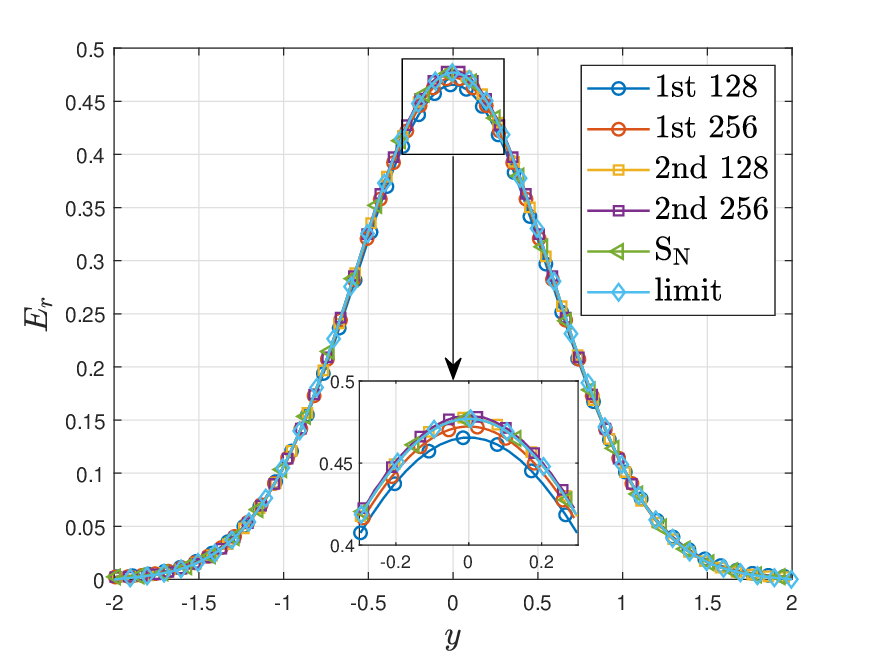}}
   \caption{(Line source problem in Sec. \ref{sec:ex4}) 
    The numerical solution $E_r$ at $t = 0.5$ along $x = 0$ when $\epsilon = 10^{-6}$. (a) The comparison of the numerical solution by first and second-order schemes, and the numerical solution of the diffusion limit equation. Here, the mesh size is $N_x = N_y = 256$. (b) The convergence behavior of $E_r$. Here, the mesh sizes are set as $N_x = N_y = 128, 256$.}
    \label{fig:linesource_cut_limit_ep_1e6}
\end{figure}

To further validate the AP property of IMEX-IM, we reduce $\epsilon$ to $\epsilon = 10^{-6}$ and keep all the other parameters the same as $\epsilon = 0.1$. The numerical solution $E_r$ obtained by the first and second-order schemes with mesh size $N_x = N_y = 256$ at $t = 0.5$ is plotted in Fig. \ref{fig:linesource_limit_ep_1e6}, where the solution of the diffusion limit equation \eqref{eq:linear_limit} is also shown. It indicates that for this even smaller $\epsilon$, IMEX-IM can keep the rotation invariance and the numerical solution behaves almost the same as the reference solution. The numerical solution along $x = 0$ is plotted in Fig. \ref{fig:linesource_t_05__CFL_04_ep_1e6_x_0_line_2}, where the reference solution obtained by the $S_N$ method and the numerical solution of the diffusion limit equation \eqref{eq:linear_limit} are also plotted. When $\epsilon = 10^{-6}$, the numerical solution matches well with the diffusion limit as well as the reference solution. The convergence behavior of IMEX-IM when $\epsilon = 10^{-6}$ is almost the same as $\epsilon = 0.1$, which is plotted in Fig. \ref{fig:linesource_t_05__CFL_04_ep_1e6_x_star_6}.

\subsection{2D lattice problem}
\label{sec:ex5}

\begin{figure}[!hptb]
  \centering
  \includegraphics[width=0.4\textwidth]{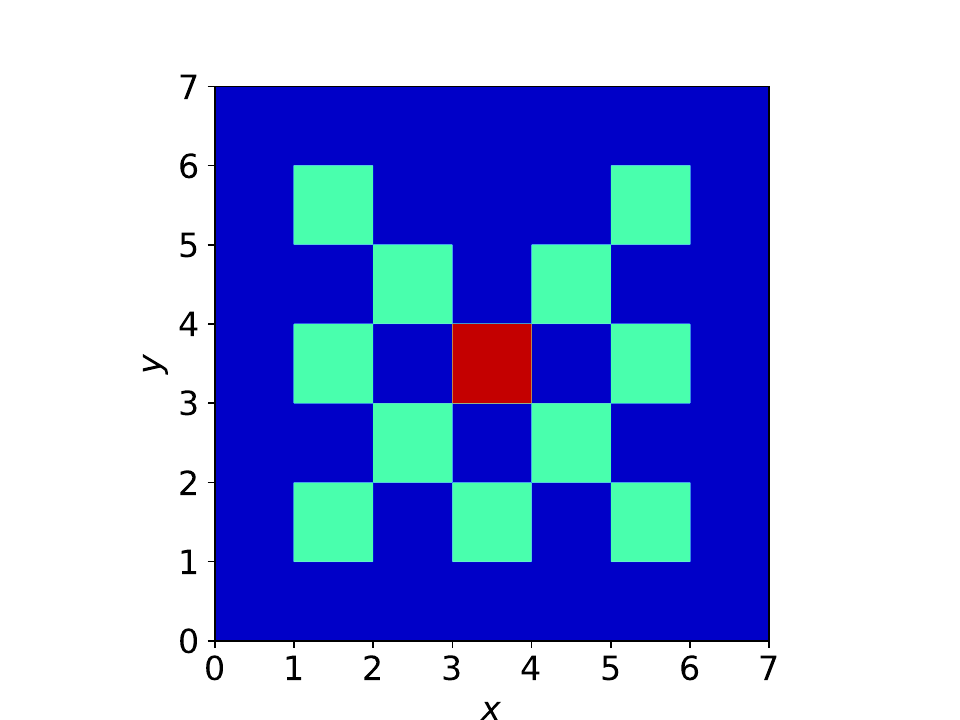}
  \caption{Layout of the 2D lattice problem. }
 \label{fig:lattice_mesh}
\end{figure}

In this section, the 2D lattice benchmark problem is studied. In this problem, the computation region is a checkerboard where $11$ absorption regions are embedded in the background of the scattering region, with a source at the center of the region. The computation region size is $[0,7]\times [0,7]$. Fig. \ref{fig:lattice_mesh} shows the specific layout of the problem. The purely scattering regions are colored blue and dark red, with $\sigma_s = 1$ and $\sigma_a = 0$, and the purely absorbing regions are colored in green with $\sigma_s = 0$ and $\sigma_a = 10$. An isotropic source $G = \dfrac{1}{4\pi}$ is turned on in the dark red region and $G$ is set to $0$ elsewhere. At the initial time $t = 0$, the specific intensity is $10^{-10} / 4\pi$. Extrapolation boundary conditions are imposed. The details of extrapolation boundary conditions are given in \cite{seibold2014starmap}.

\begin{figure}[!hptb]
  \centering
  \subfloat[$\log_{10}E_r (\epsilon = 1)$, 1st order]{
    \label{fig:lattice_t_32_CFL_04_ep_1_1order_imagesc}
    \includegraphics[width=0.30\textwidth]{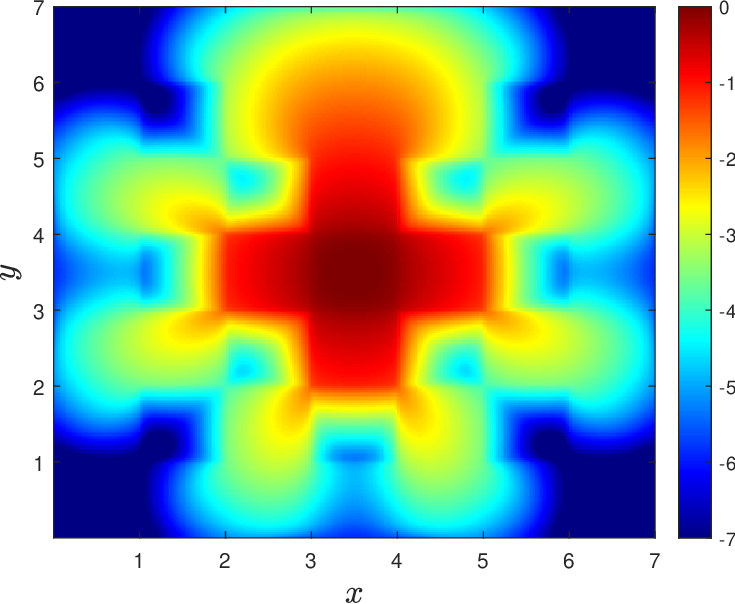}}
  \hfill
   \subfloat[$\log_{10}E_r (\epsilon = 1)$, 2nd order]{
    \label{fig:lattice_t_32__CFL_04_ep_1_2order_imagesc}
    \includegraphics[width=0.30\textwidth]{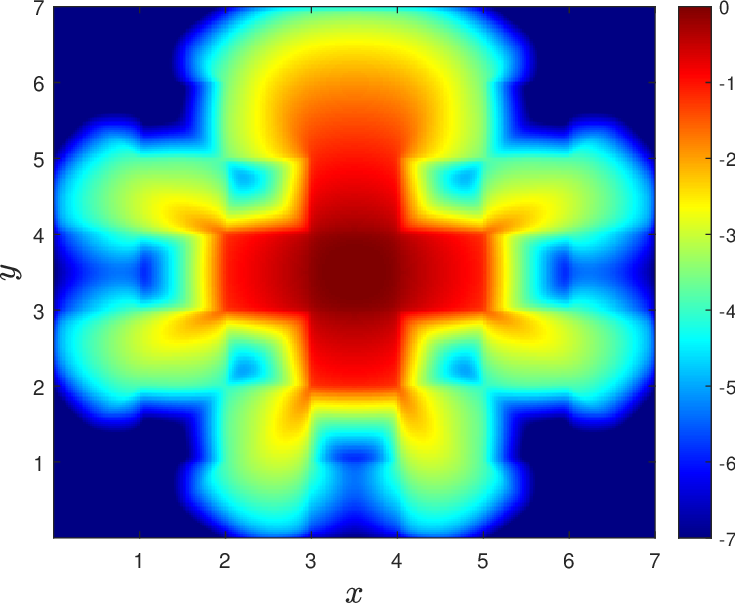}}
   \hfill
   \subfloat[$\log_{10}E_r$, Starmap]{
    \label{fig:lattice_t_32_CFL_04_ep_1_Starmap_imagesc}
    \includegraphics[width=0.30\textwidth]{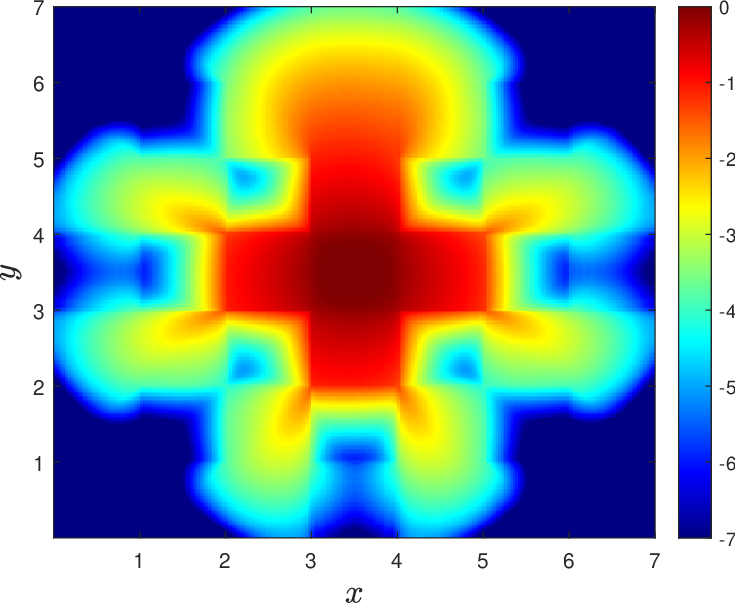}}
    \caption{(2D lattice problem in Sec. \ref{sec:ex5}) The contour plot of the log density for the
    lattice problem at $t = 3.2$ with $\epsilon = 1$. Here, the mesh size is $N_x = N_y = 280$. (a) Contour plot of $\log_{10} E_r$ by 1st-order scheme. (b) Contour plot of $\log_{10} E_r$ by 2nd-order scheme. (c) Contour plot of $\log_{10} E_r$ by Starmap.}
    \label{fig:lattice_ep_1_imagesc}
\end{figure}

\begin{figure}[!hptb]
  \centering
     \subfloat[ $\log_{10} E_r (\epsilon = 1), x = 3.5$]{
    \label{fig:lattice_t_32__CFL_04_ep_1_y_35_line_2}
    \includegraphics[width=0.45\textwidth]{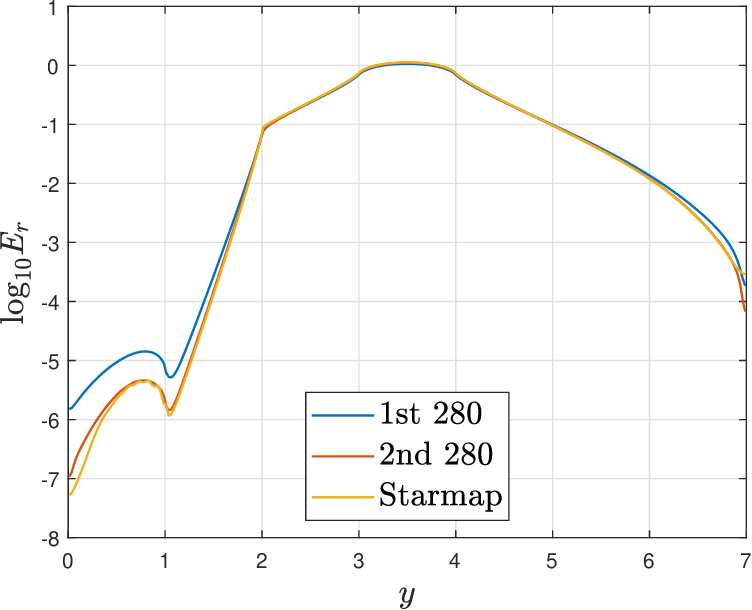}}
    \hfill
  \subfloat[ $\log_{10} E_r (\epsilon = 1), x = 3.5$]{
    \label{fig:lattice_t_32__CFL_04_ep_1_y_35_star_6}
    \includegraphics[width=0.45\textwidth]{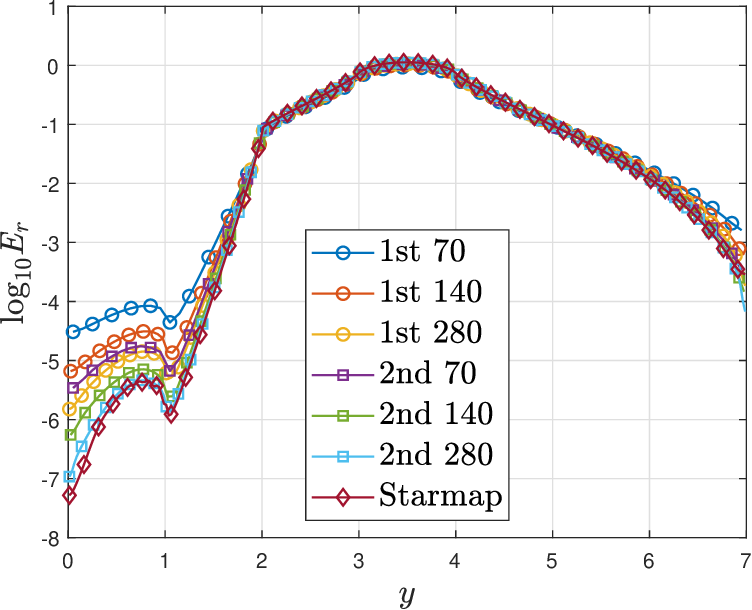}} \\
     \subfloat[ $\log_{10} E_r (\epsilon = 1), y = 3.5$]{
    \label{fig:lattice_t_32__CFL_04_ep_1_x_35_line_2}
    \includegraphics[width=0.45\textwidth]{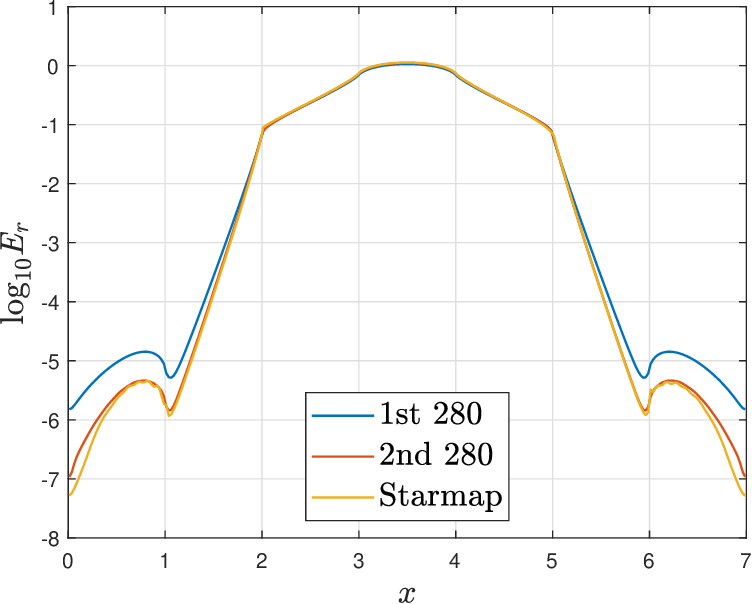}}
    \hfill
     \subfloat[ $\log_{10} E_r (\epsilon = 1), y = 3.5$]{
    \label{fig:lattice_t_32__CFL_04_ep_1_x_35_star_6}
    \includegraphics[width=0.45\textwidth]{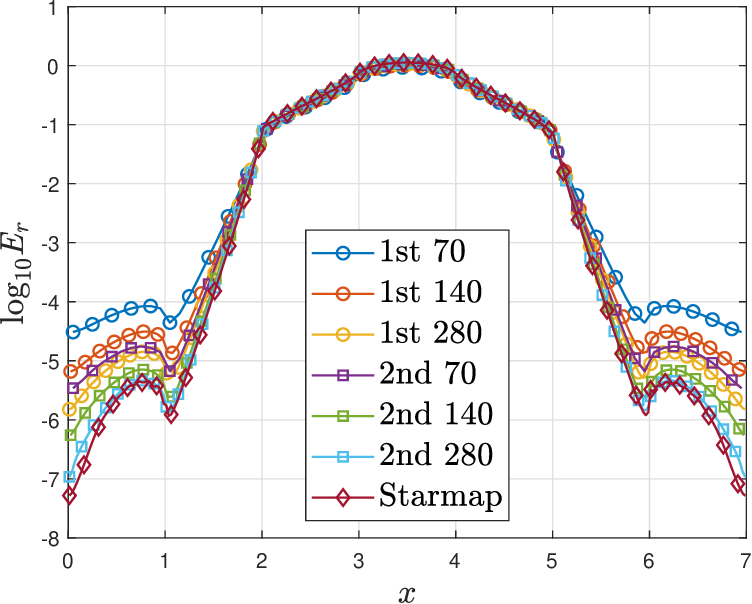}}
   \caption{(2D lattice problem in Sec. \ref{sec:ex5}) 
   The slice cut of the numerical solution $\log_{10} E_r$ at $t = 3.2$ when $\epsilon = 1$. (a) The slice cut along $x = 3.5$ of $\log_{10} E_r$. Here, the mesh size is $N_x = N_y = 280$. (b) The convergence behavior of $\log_{10} E_r$ along $x = 3.5$. Here, the mesh sizes are set as $N_x = N_y = 70, 140, 280$. (c) The slice cut along $y = 3.5$ of $\log_{10} E_r$. Here, the mesh size is $N_x = N_y = 280$. (d) The convergence behavior of $\log_{10} E_r$ along $y = 3.5$. Here, the mesh sizes are set as $N_x = N_y = 70, 140, 280$.}   
   \label{fig:lattice_cut_ep_1}
\end{figure}

We first set the parameter $\epsilon = 1$. The mesh size is chosen as $N_x = N_y = 280$ with the expansion number $M = 39$ adopted here. Fig. \ref{fig:lattice_ep_1_imagesc} shows the numerical results of the logarithm of radiation energy density $\log_{10}E_r$ \eqref{eq:energy_Er} obtained by the first and second-order schemes at $t = 3.2$ as well as the reference solution obtained by StarMAP \cite{seibold2012starmap, seibold2014starmap}. We can see that the particle beams leak between the corners of the absorption region, which is the same as that in \cite{brunner2002forms, seibold2014starmap}. The numerical solution matches well with the reference solution. To compare the numerical solution and the reference solution more clearly, the numerical solution along $x = 3.5$ and $y = 3.5$ is plotted in Fig. \ref{fig:lattice_cut_ep_1}. From it, we can see there are some oscillations in the reference solution near the absorption corner while the numerical solution holds smooth. The convergence of the numerical solution with $N_x = N_y = 70, 140$ and $280$ is illustrated in Fig. \ref{fig:lattice_t_32__CFL_04_ep_1_y_35_star_6} and \ref{fig:lattice_t_32__CFL_04_ep_1_x_35_star_6}, where the second-order scheme behaves much better than the first-order scheme and they are all converging to the reference solution. 

\begin{figure}[!hptb]
  \centering
  \subfloat[$\log_{10} E_r (\epsilon = 0.1)$, 1st order]{
    \label{fig:lattice_t_05_CFL_04_ep_1e1_1order_imagesc}
    \includegraphics[width=0.30\textwidth]{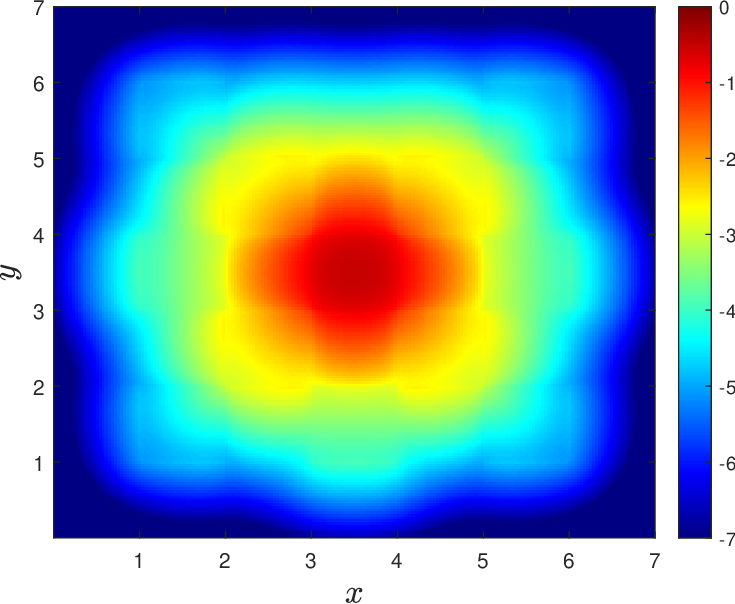}}
  \hfill
   \subfloat[$\log_{10} E_r (\epsilon = 0.1)$, 2nd order]{
    \label{fig:lattice_t_05__CFL_04_ep_1e1_2order_imagesc}
    \includegraphics[width=0.30\textwidth]{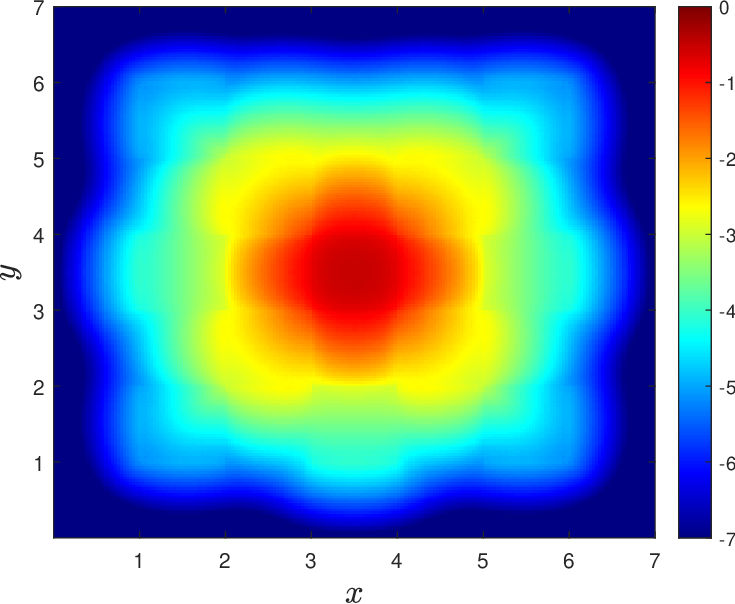}}
     \hfill
   \subfloat[$\log_{10} E_r (\epsilon = 0.1)$, $S_N$ method]{
    \label{fig:lattice_t_05__CFL_04_ep_1e1_Sn_imagesc}
    \includegraphics[width=0.30\textwidth]{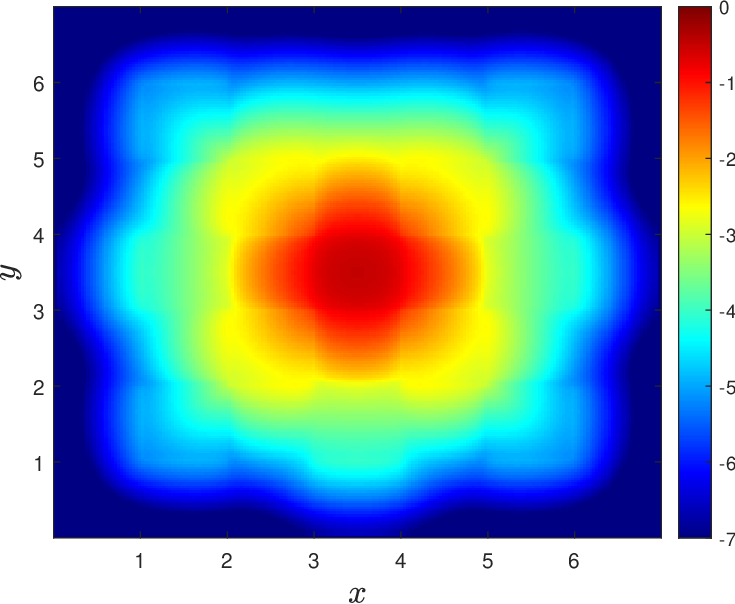}} 
  \caption{(2D lattice problem in Sec. \ref{sec:ex5}) The contour plot of the log density for the lattice problem at $t = 0.5$ with $\epsilon = 0.1$. Here, the mesh size is $N_x = N_y = 200$. (a) Contour plot of $\log_{10} E_r$ by 1st-order scheme. (b) Contour plot of $\log_{10} E_r$ by 2nd-order scheme. (c) Contour plot of $\log_{10} E_r$ by $S_N$ method.}
    \label{fig:lattice_ep_01_imagesc}
\end{figure}

\begin{figure}[!hptb]
  \centering
  \subfloat[ $\log_{10} E_r (\epsilon = 0.1), x= 3.5$]{
    \label{fig:lattice_t_05__CFL_04_ep_1e1_y_35_line_2}
    \includegraphics[width=0.4\textwidth]{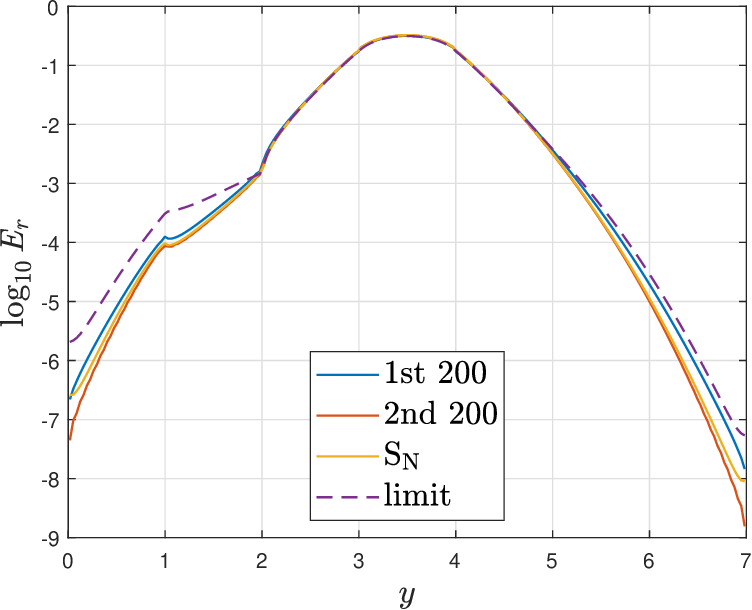}}  
    \hfill
    \subfloat[ $\log_{10} E_r (\epsilon = 0.1), y = 3.5$]{
    \label{fig:lattice_t_05__CFL_04_ep_1e1_x_35_line_2}
    \includegraphics[width=0.4\textwidth]{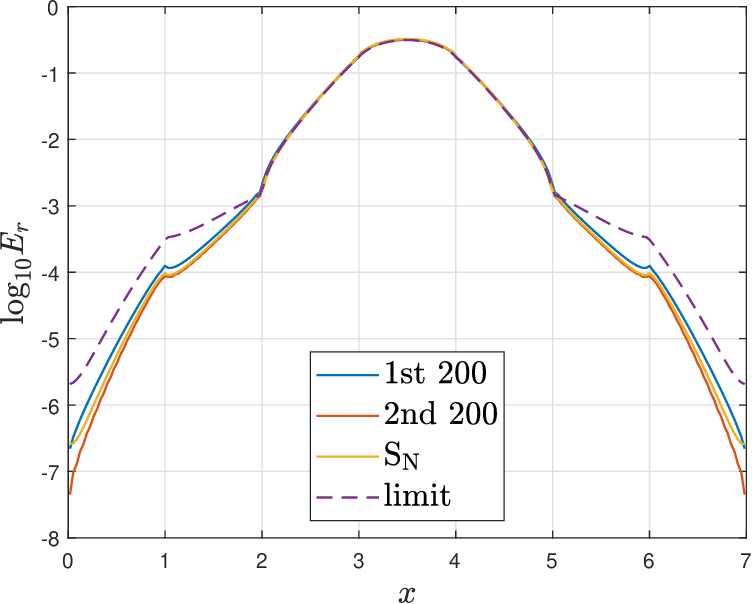}}
    \caption{(2D lattice problem in Sec. \ref{sec:ex5}) 
     The slice cut of the numerical solution $\log_{10} E_r$ at $t = 0.5$ when $\epsilon = 0.1$. Here, the mesh size is $N_x = N_y = 200$. (a) The slice cut along $x = 3.5$ of $\log_{10} E_r$. (b) The slice cut along $y = 3.5$ of $\log_{10} E_r$. }
    \label{fig:lattice_cut_ep_1e1}
\end{figure}

To validate the AP property of IMEX-IM, we decrease $\epsilon$ to $0.1$. In this case, the pure absorption regions become the strong absorption regions, with $\sigma_s = 0.1$ and $\sigma_a = 9.9$, while other parameters are the same as the case $\epsilon = 1$. We set the mesh size as $N_x = N_y = 200$, the expansion order as $M = 7$. Fig. \ref{fig:lattice_ep_01_imagesc} shows the numerical solution of $\log_{10} E_r$ obtained by the first- and second-order schemes at $t = 0.5$ as well as the reference solution obtained by the $S_N$ method. The slice cut along $x = 3.5$ and $y = 3.5$ of $\log_{10} E_r$ are displayed in Fig. \ref{fig:lattice_cut_ep_1e1} as well as the reference solution obtained by the $S_N$ method, and the solution of the diffusion limit \eqref{eq:linear_limit}. The numerical solution by the second-order scheme behaves much better than that by the first-order scheme and the difference with the reference solution by the $S_N$ method is quite small. It also clearly shows some distance between the numerical solution of RTE \eqref{eq:linear_RTE} and that of the diffusion limit when $\epsilon = 0.1$.

\begin{figure}[!hptb]
  \centering
       \subfloat[$\log_{10}E_r (\epsilon = 10^{-6})$, 1st order]{
    \label{fig:lattice_t_05_CFL_04_ep_1e6_1order_imagesc}
    \includegraphics[width=0.30\textwidth]{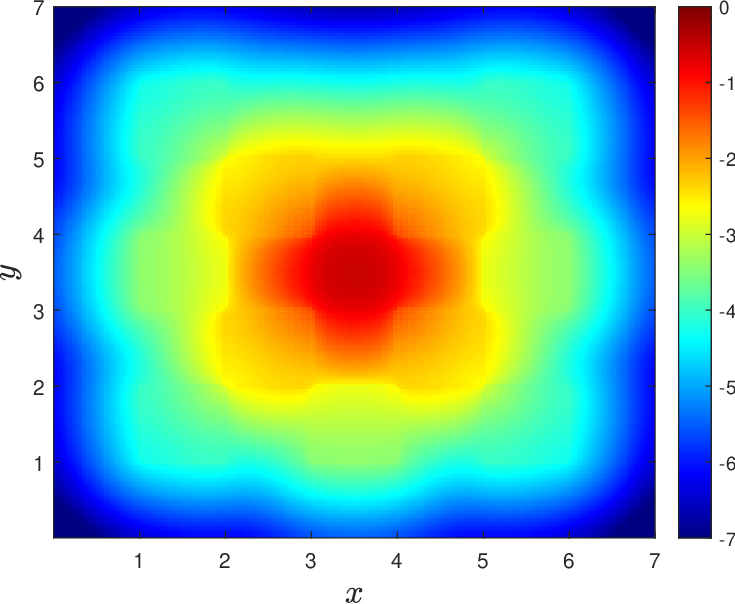}}
  \hfill
   \subfloat[$\log_{10}E_r (\epsilon = 10^{-6})$, 2nd order]{
    \label{fig:lattice_t_05__CFL_04_ep_1e6_2order_imagesc}
    \includegraphics[width=0.30\textwidth]{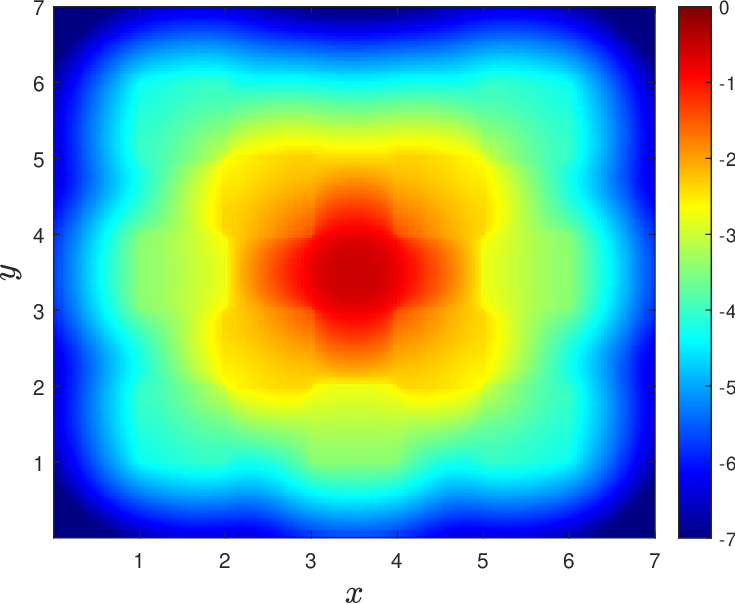}}
   \hfill
   \subfloat[$\log_{10}E_r$ of the diffusion limit]{
    \label{fig:lattice_t_05_CFL_04__limit_imagesc}
    \includegraphics[width=0.30\textwidth]{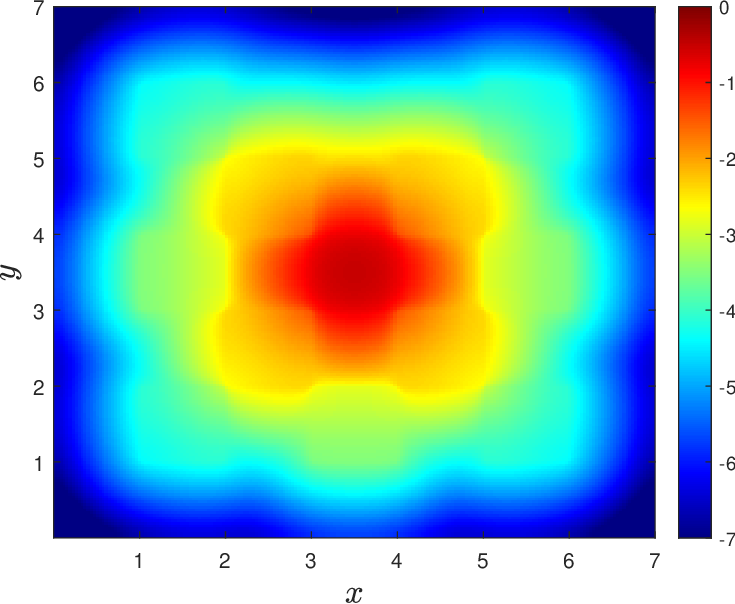}}
  \caption{(2D lattice problem in Sec. \ref{sec:ex5}) The contour plot of the log density for the lattice problem at $t = 0.5$ with $\epsilon = 10^{-6}$. Here, the mesh size is $N_x = N_y = 200$. (a) Contour plot of $\log_{10} E_r$ by 1st-order scheme. (b) Contour plot of $\log_{10} E_r$ by 2nd-order scheme. (c) Contour plot of $\log_{10} E_r$ of the diffusion limit.}
\label{fig:lattice_ep_06_imagesc}
\end{figure}

\begin{figure}[!hptb]
  \centering
    \subfloat[ $\log_{10} E_r (\epsilon = 10^{-6}), x = 3.5$]{
    \label{fig:lattice_t_05__CFL_04_ep_1e6_y_35_line_2}
    \includegraphics[width=0.4\textwidth]{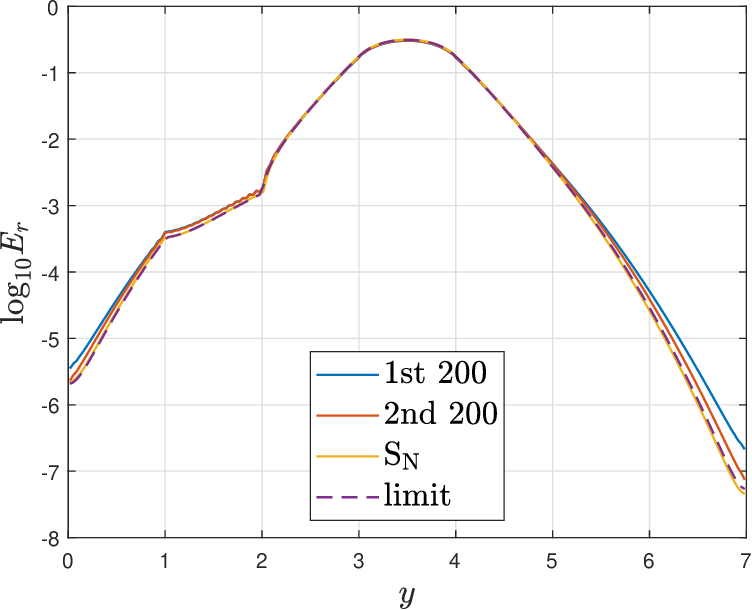}} 
    \hfill
    \subfloat[ $\log_{10} E_r (\epsilon = 10^{-6}), y = 3.5$]{
    \label{fig:lattice_t_05__CFL_04_ep_1e6_x_35_line_2}
    \includegraphics[width=0.4\textwidth]{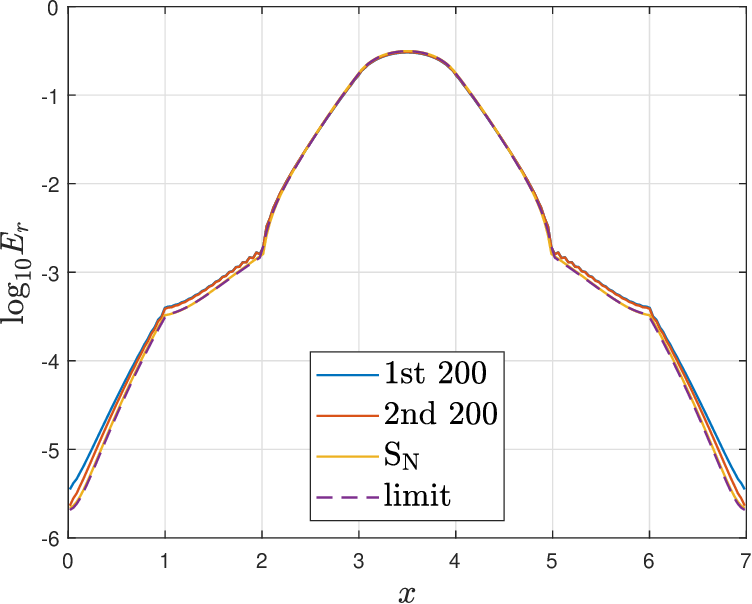}}
    \caption{(2D lattice problem in Sec. \ref{sec:ex5}) 
    The slice cut of the numerical solution $E_r$ at $t = 0.5$ when $\epsilon = 10^{-6}$. Here, the mesh size is $N_x = N_y = 200$. (a) The slice cut along $x = 3.5$ of $\log_{10} E_r$. (b) The slice cut along $y = 3.5$ of $\log_{10} E_r$.}
    \label{fig:lattice_cut_ep_1e6}
\end{figure}

To further validate the AP property of IMEX-IM, $\epsilon$ is reduced to $\epsilon = 10^{-6}$ with other parameters the same as that of $\epsilon = 0.1$. Similarly, the numerical solution by the first- and second-order schemes at $t = 0.5$ as well as the numerical solution of the diffusion limit \eqref{eq:linear_limit} is plotted in Fig. \ref{fig:lattice_ep_06_imagesc} with the slice cut along $x = 3.5$ and $y = 3.5$ shown in Fig. \ref{fig:lattice_cut_ep_1e6}. We can find that the numerical solution of RTE \eqref{eq:linear_RTE} is almost the same as that of the diffusion limit \eqref{eq:linear_limit}, and the performance of the second-order scheme is better than the first-order scheme.

%%%%%%%%%%%%%%%%%%%%%%%%%%%%%%%%%%%%%%%%%%%%%%%%%%%%%%%%%%%%%

\subsection{2D Riemann problem}
\label{sec:ex6}
\begin{figure}[!hptb]
  \centering
   \subfloat[Layout of the 2D Riemann problem]{
    \label{fig:Rie_2D_set}
    \includegraphics[width=0.4\textwidth]{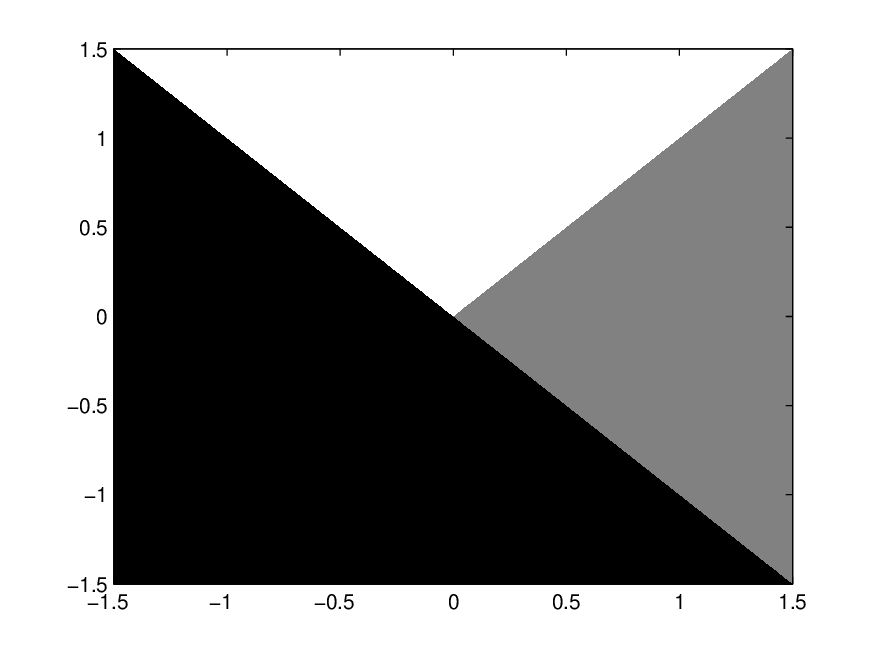}}
    \hfill
     \subfloat[ Layout of the 1D Riemann problem]{
    \label{fig:Rie_1D_set}
    \includegraphics[width=0.4\textwidth]{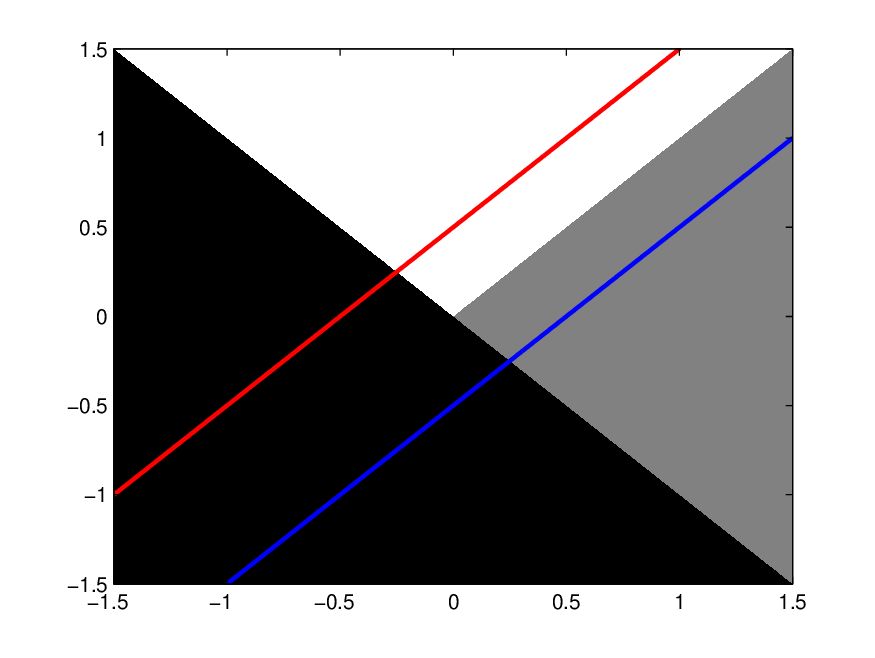}}
  \caption{(2D Riemann problem in Sec. \ref{sec:ex6}) The problem setting and computational regime for the 2D Riemann problem. (a) The layout of the 2D Riemann problem. (b) The layout of the 1D Riemann problem along the red line $y = x + 0.5$ and the blue line $y = x - 0.5$.}
 \label{fig:Riemann_set}
\end{figure}

In this section, a Riemann problem is simulated with a nonuniform absorption cross-section. The computational domain is specified as $[-1.5,1.5] \times [-1.5,1.5]$. Fig. \ref{fig:Rie_2D_set}  illustrates the layout of the problem, where black and gray regions denote the medium with $\sigma_a = 1$ and the white region indicates the medium with $\sigma_a = 10$. The initial temperature is set to $T = 1$ in the black region and $T = 0.1$ in gray and white regions, with the radiation density in equilibrium 
\begin{equation}
I = \frac{1}{4\pi}acT^4.
\end{equation}
Similar to 1D Marshak wave problems, a constant isotropic incident radiation intensity with a Planckian distribution at $T = 1$ is set to the left and bottom boundaries, while $T = 0.1$ is set to the right and top boundaries, with a Marshak-type boundary condition \cite{semi2008Ryan} applied. The other parameters are set as $C_v = 1, a = 1$ and $c = 1$.

\begin{figure}[!hptb]
  \centering
  \subfloat[$T_r (\epsilon = 1)$, 1st order]{
    \label{fig:Riemann_Tr_Pn_ep_1_1order}
    \includegraphics[width=0.3\textwidth]{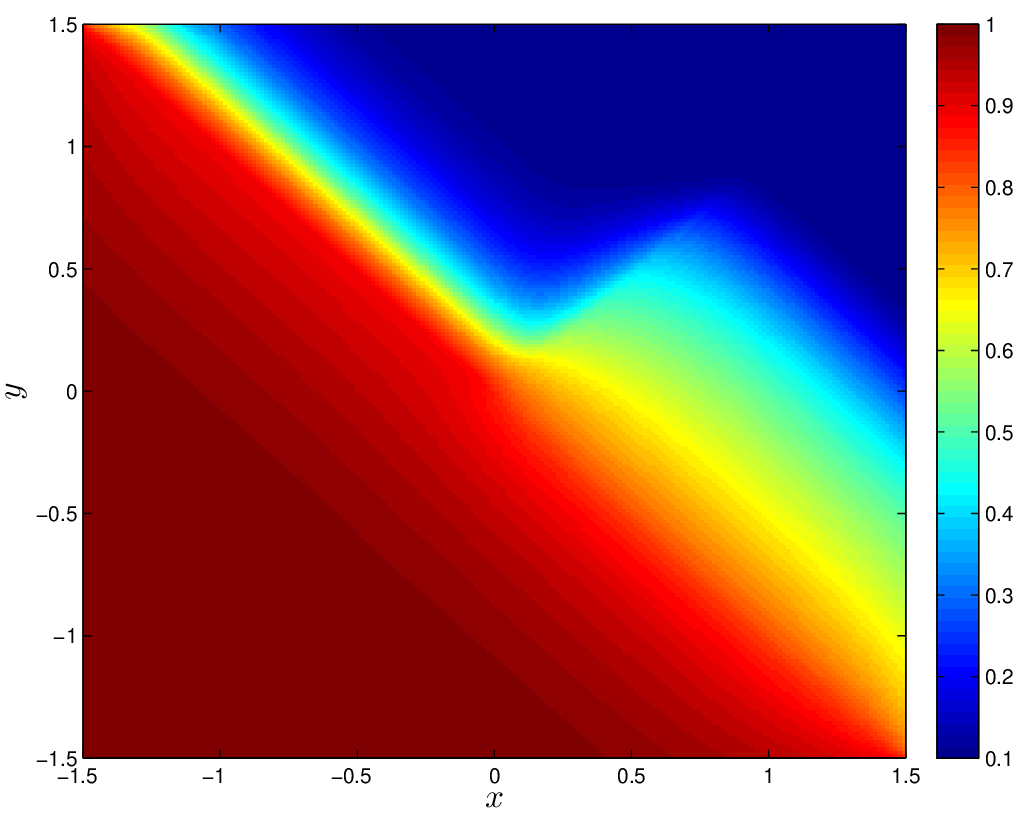}}
     \hfill 
    \subfloat[$T_r (\epsilon = 1)$, 2nd order]{
    \label{fig:Riemann_Tr_Pn_ep_1_2order}
    \includegraphics[width=0.3\textwidth]{./images/Rie_t_1__CFL_04_ep_1_x_Tr_imagesc_1order.eps}}
     \hfill 
  \subfloat[$T_r (\epsilon = 1)$, $S_N$ method]{
    \label{fig:Riemann_Tr_Sn_ep_1}
    \includegraphics[width=0.3\textwidth]{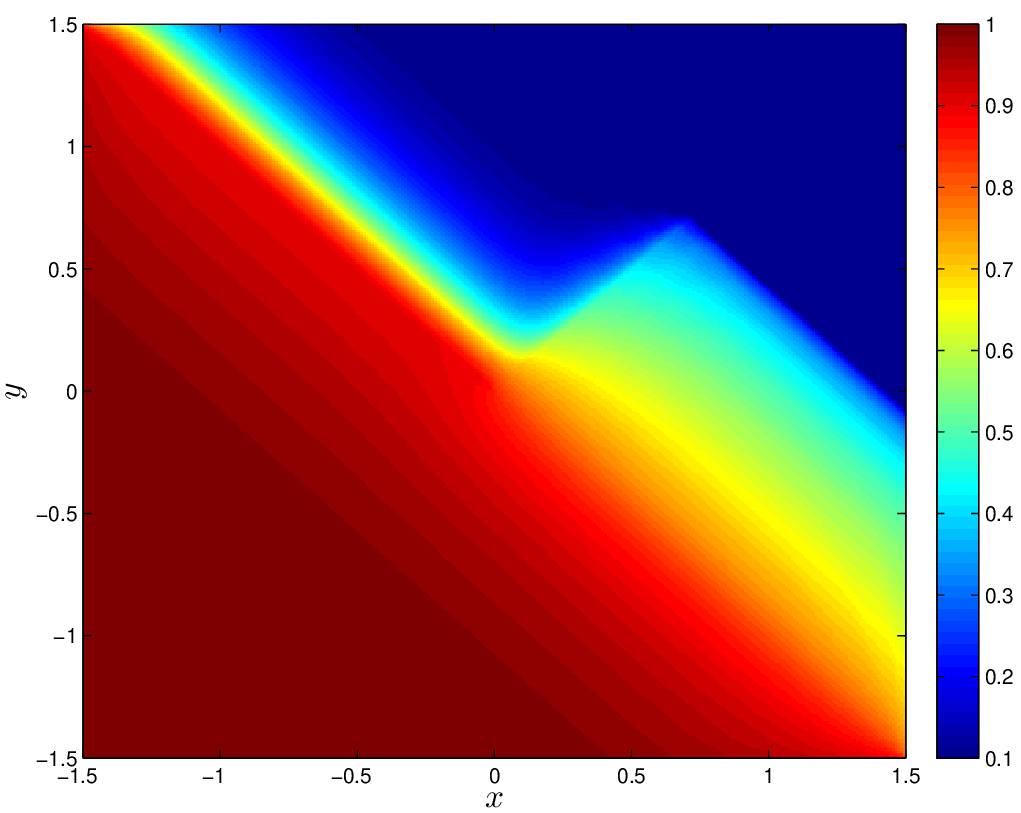}}
    \hfill 
  \subfloat[$T (\epsilon = 1)$, 1st order]{
    \label{fig:Riemann_T_Pn_ep_1_1order}
    \includegraphics[width=0.3\textwidth]{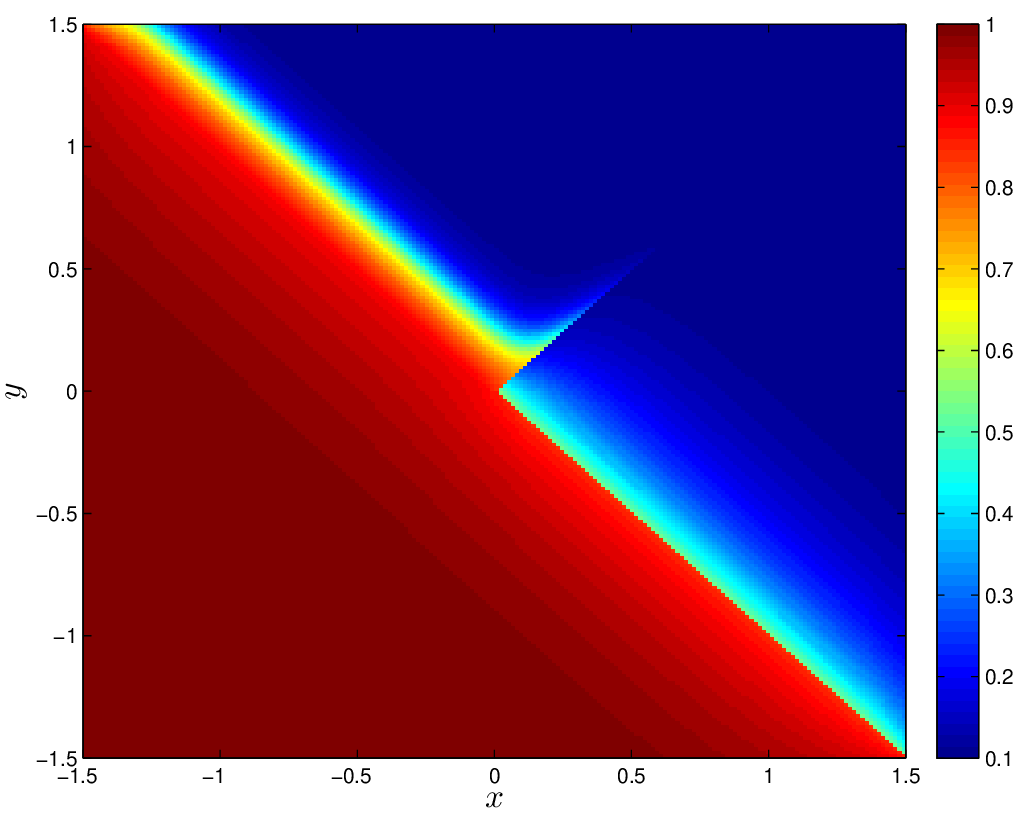}}
    \hfill 
  \subfloat[$T (\epsilon = 1)$, 2nd order]{
    \label{fig:Riemann_T_Pn_ep_1_2order}
    \includegraphics[width=0.3\textwidth]{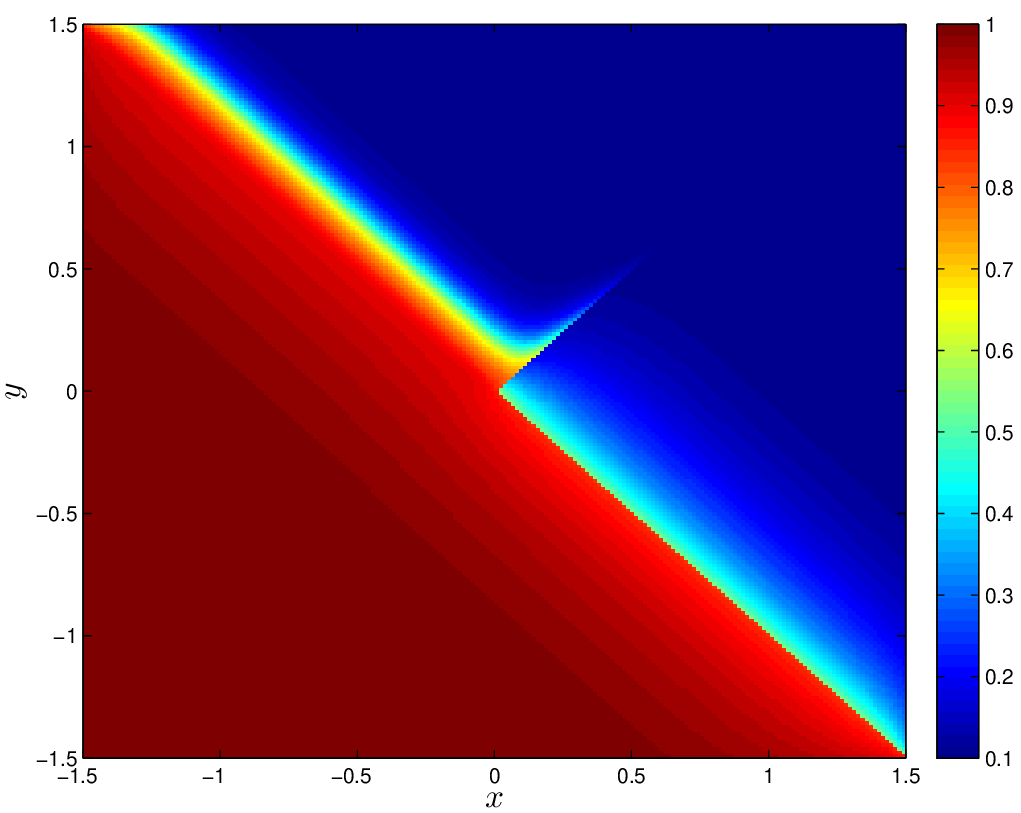}}
  \hfill 
   \subfloat[$T (\epsilon = 1)$, $S_N$ method]{
    \label{fig:Riemann_T_Sn_ep_1}
    \includegraphics[width=0.3\textwidth]{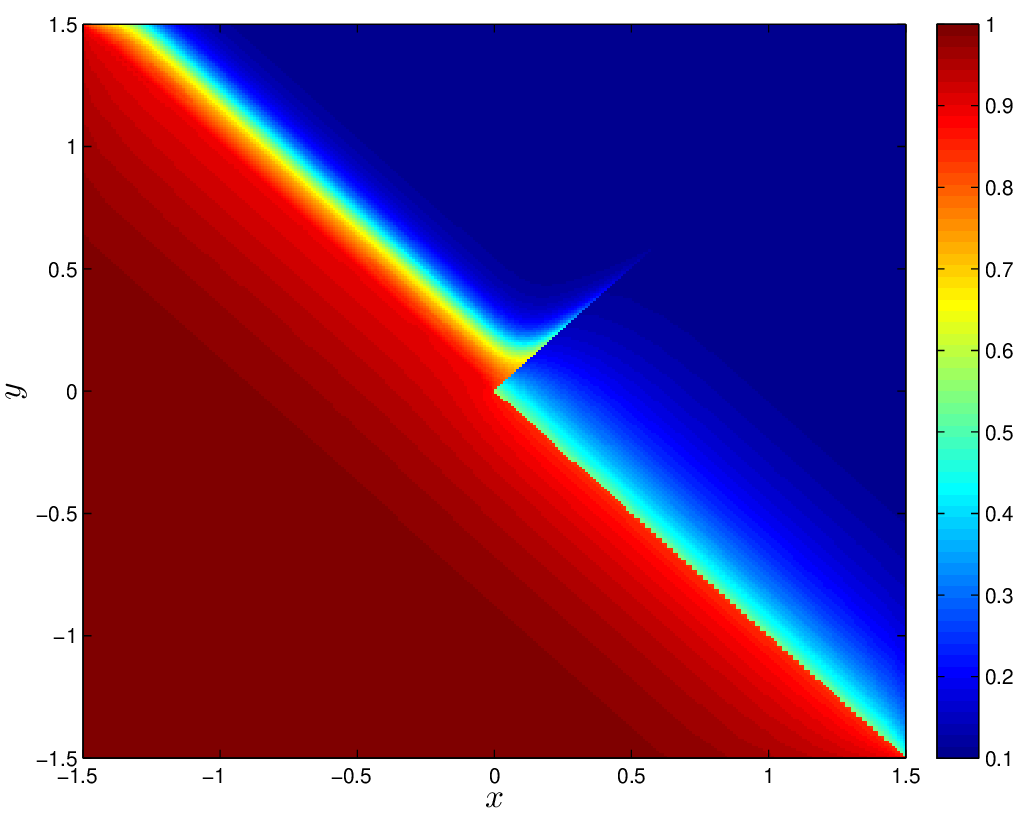}}
  \hfill 
  \caption{(2D Riemann problem in Sec. \ref{sec:ex6}) The numerical solution of $T_r$ and $T$ for the 2D Riemann problem with $\epsilon = 1$ at $t = 1$. The first row is that for $T_r$ and the second row is that for $T$. }
    \label{fig:Riemann_1}
\end{figure}

\begin{figure}[!htbp]
 \flushleft
  \subfloat[$T_r (\epsilon = 1, y = x + 0.5)$]{
    \label{fig:Riemann_Tr_Pn_plus_ep_1_line3}
    \includegraphics[width=0.4\textwidth]{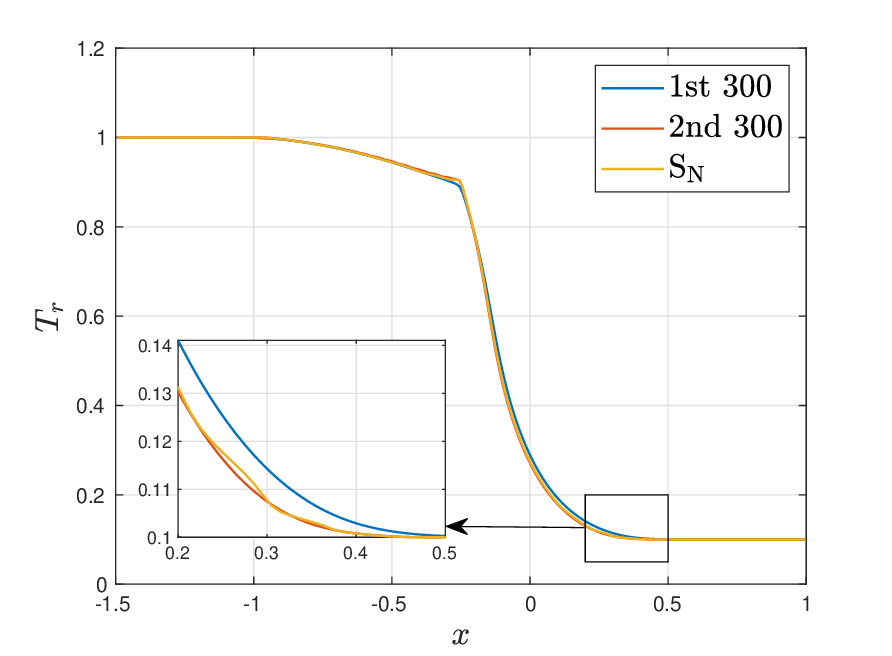}}
  \hfill 
    \subfloat[$T_r (\epsilon = 1, y = x - 0.5)$.]{
    \label{fig:Riemann_Tr_Pn_minus_ep_1_line3}
    \includegraphics[width=0.4\textwidth]{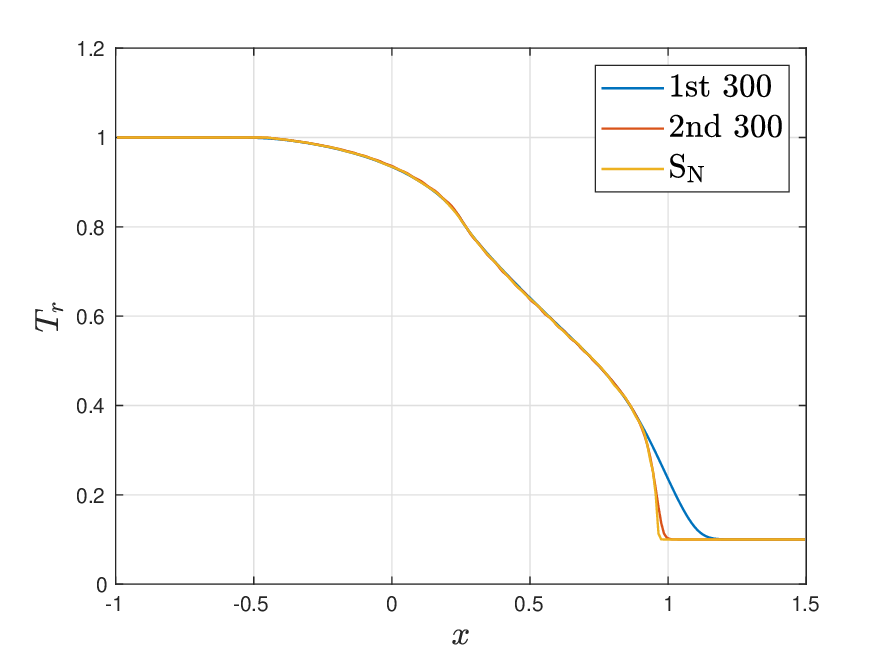}}
  \hfill 
  \subfloat[$T (\epsilon = 1, y = x + 0.5)$.]{
    \label{fig:Riemann_T_Pn_plus_ep_1_line3}
    \includegraphics[width=0.4\textwidth]{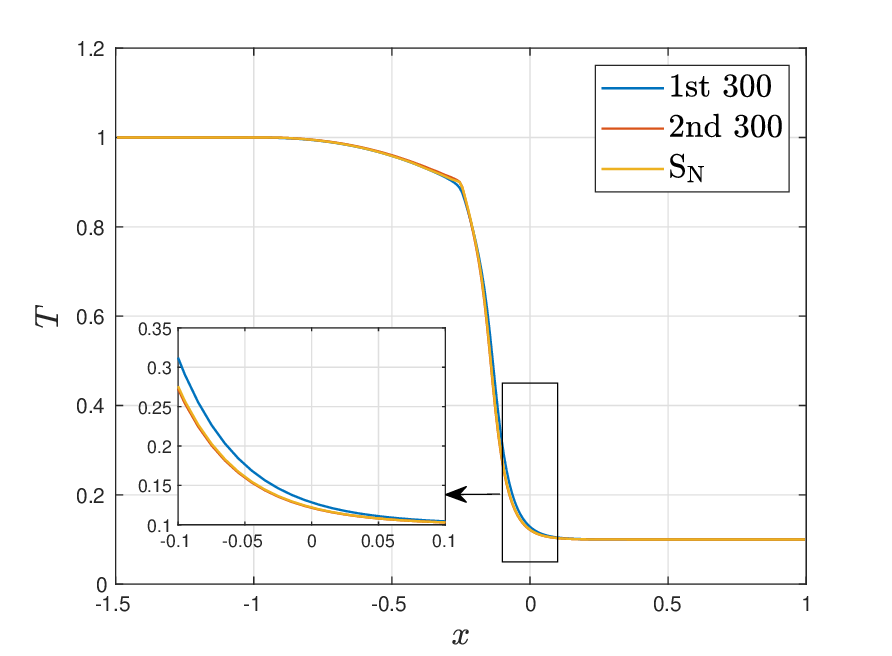}}
    \hfill
    \subfloat[$T (\epsilon = 1, y = x - 0.5)$.]{
    \label{fig:Riemann_T_Pn_minus_ep_1_line3}
    \includegraphics[width=0.4\textwidth]{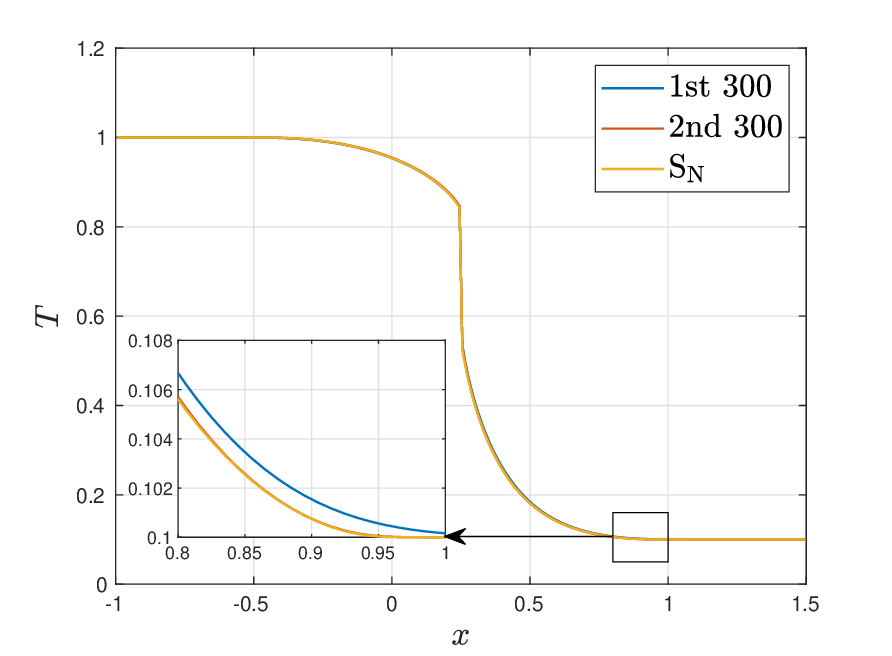}}
  \caption{(2D Riemann problem in Sec. \ref{sec:ex6}) The slice plot of $T_r$ and $T$ for the 2D Riemann problem with $\epsilon = 1$ at $t = 1$. (a) $T_r$ along $y = x + 0.5$. (b)  $T_r$ along $y = x - 0.5$. (c) $T$ along $y = x + 0.5$. (d) $T$ along $y = x - 0.5$. }
    \label{fig:Riemann_ep_1_cut}
\end{figure}

We first set $\epsilon = 1$. Then, the expansion order is chosen as $M = 39$ for the $P_N$ method, and a mesh of $300 \times 300$ is utilized in spatial space. The final computation time is $t = 1$. The reference solution is provided by the $S_N$ method. Fig. \ref{fig:Riemann_1} shows the contour plots of the radiation temperature $T_r$ and material temperature $T$, where the radiation temperature is low and the material temperature is high in the white area, while the radiation temperature is high and the material temperature is low in the gray area. This result is consistent with the behavior of the opacity $\sigma$ and validates the accuracy of the numerical solution.
Fig. \ref{fig:Riemann_ep_1_cut} displays $T_r$ and $T$ along $y = x + 0.5$ and $y = x - 0.5$, perpendicular to the black area's interface region. The numerical results from the second-order scheme exhibit higher resolution and align more closely with the reference solution than those from the first-order scheme on the same $300 \times 300$ mesh, demonstrating a high level of consistency between the numerical and reference solutions.
\begin{figure}[!hptb]
  \centering
  \subfloat[$T_r (\epsilon = 0.1)$, 1st order]{
    \label{fig:Riemann_Tr_Pn_ep_1e1_1order}
    \includegraphics[width=0.23\textwidth]{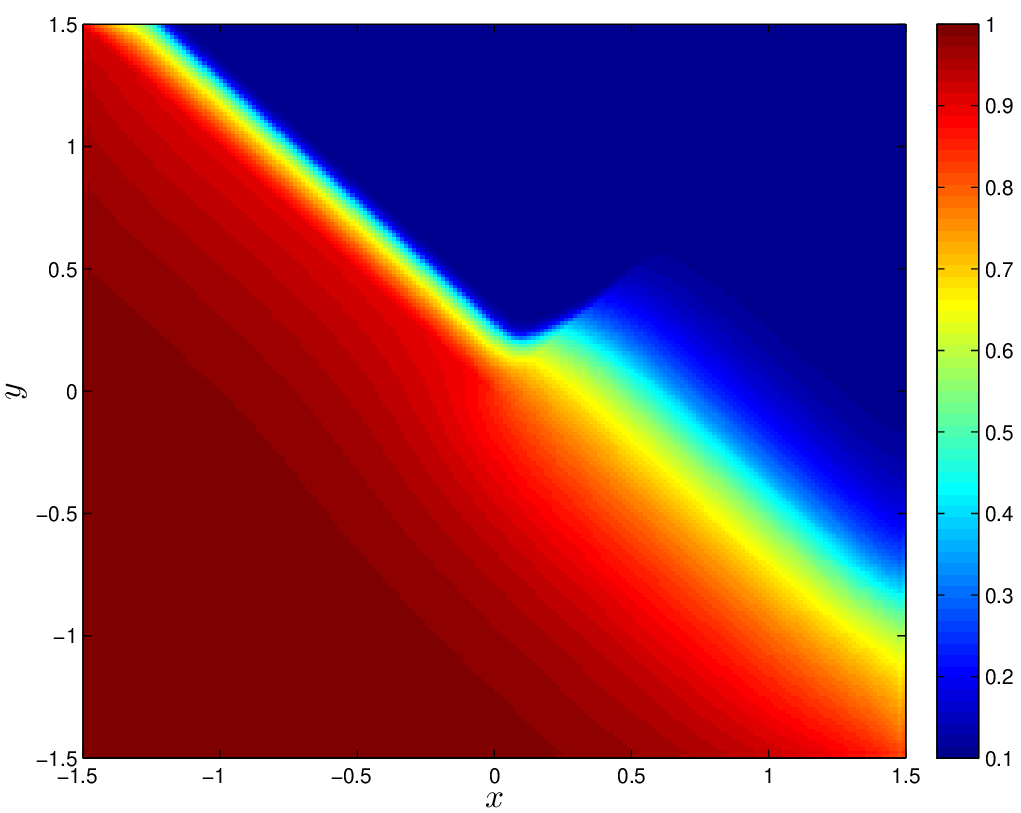}}
    \hfill 
  \subfloat[$T_r (\epsilon = 0.1)$, 2nd order]{
    \label{fig:Riemann_Tr_Pn_ep_1e1_2order}
    \includegraphics[width=0.23\textwidth]{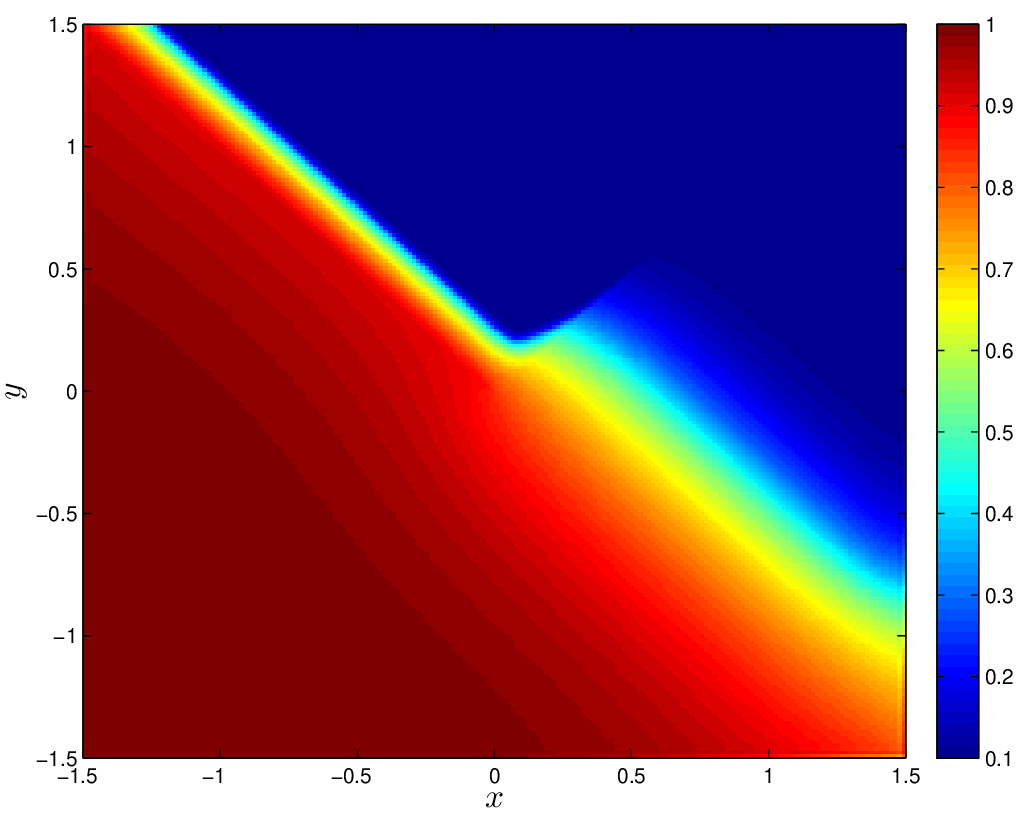}}
  \hfill 
   \subfloat[$T (\epsilon = 0.1)$, 1st order]{
    \label{fig:Riemann_T_Pn_ep_1e1_1order}
    \includegraphics[width=0.23\textwidth]{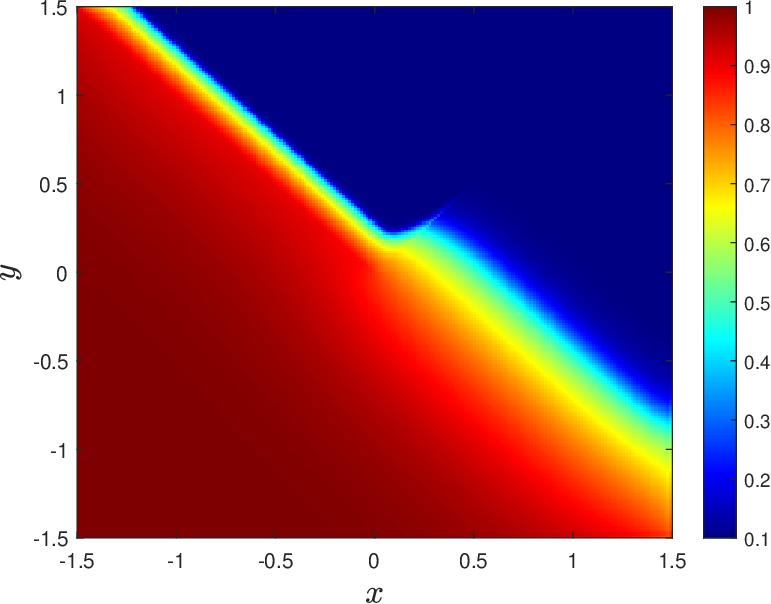}}
    \hfill 
  \subfloat[$T (\epsilon = 0.1)$, 2nd order]{
    \label{fig:Riemann_T_Pn_ep_1e1_2order}
    \includegraphics[width=0.23\textwidth]{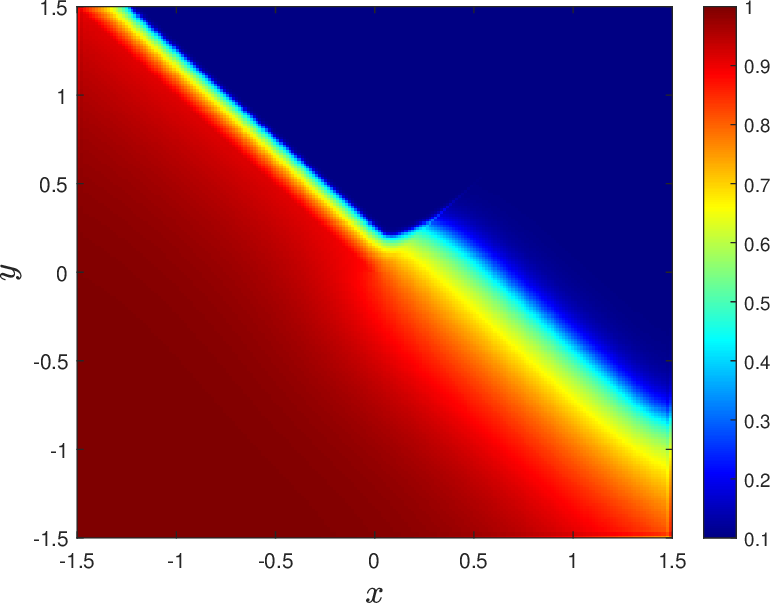}}
  \caption{(2D Riemann problem in Sec. \ref{sec:ex6}) 
  The numerical solution of $T_r$ and $T$ for the 2D Riemann problem with $\epsilon = 0.1$ at $t = 0.5$. (a) Contour plot of $T_r$ by 1st-order scheme. (b) Contour plot of $T_r$ by 2nd-order scheme. (c) Contour plot of $T$ by 1st-order scheme. (d) Contour plot of $T$ by 2nd-order scheme.}
    \label{fig:Riemann_ep_1e1}
\end{figure}

\begin{figure}[!htbp]
  \centering
  \subfloat[$T_r (\epsilon = 0.1, y = x + 0.5)$.]{
    \label{fig:Riemann_Tr_Pn_plus_ep_1e1_line3}
    \includegraphics[width=0.23\textwidth]{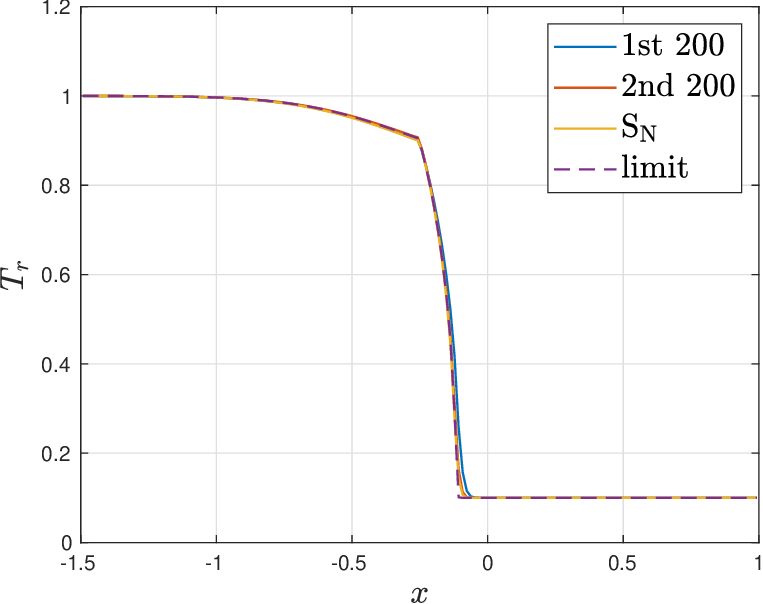}}
  \hfill     
    \subfloat[$T_r (\epsilon = 0.1, y = x - 0.5)$.]{
    \label{fig:Riemann_Tr_Pn_minus_ep_1e1_line3}
    \includegraphics[width=0.23\textwidth]{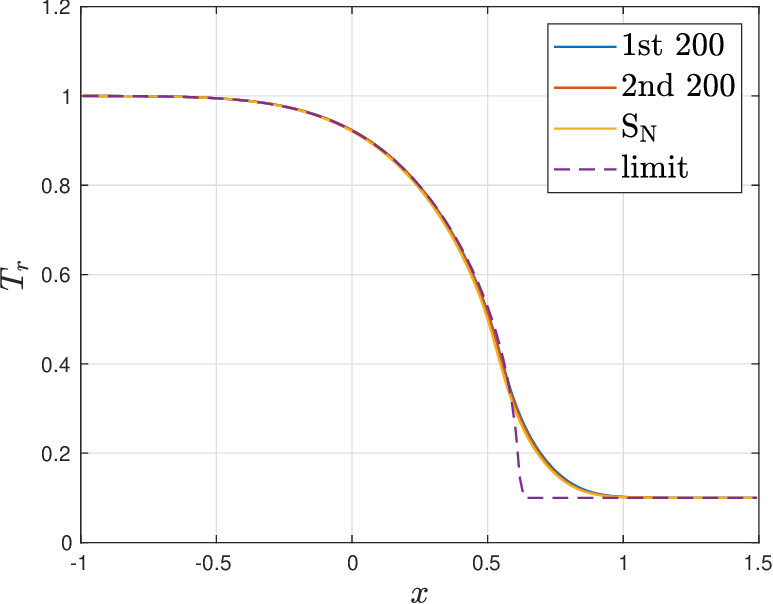}}
    \hfill
    \subfloat[$T (\epsilon = 0.1, y = x + 0.5)$.]{
    \label{fig:Riemann_T_Pn_plus_ep_1e1_line3}
    \includegraphics[width=0.23\textwidth]{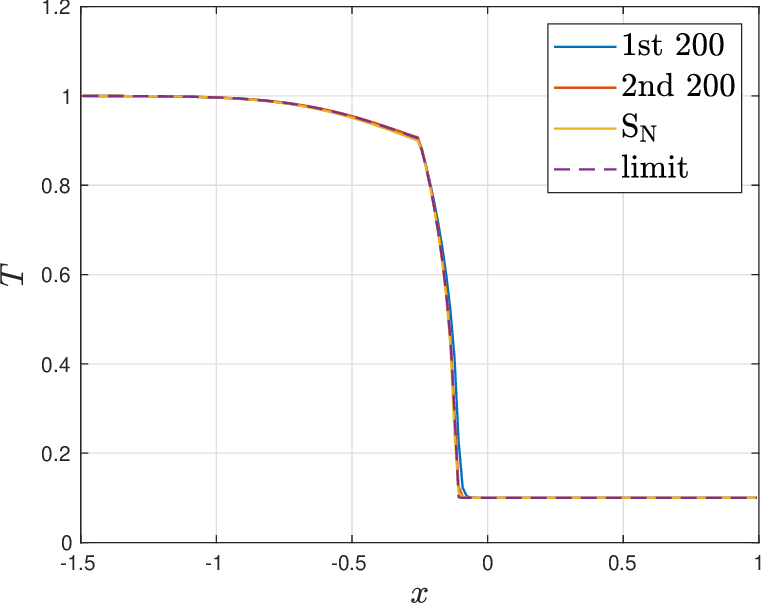}}
   \hfill 
    \subfloat[$T (\epsilon = 0.1, y = x - 0.5)$.]{
    \label{fig:Riemann_T_Pn_minus_ep_1e1_line3}
    \includegraphics[width=0.23\textwidth]{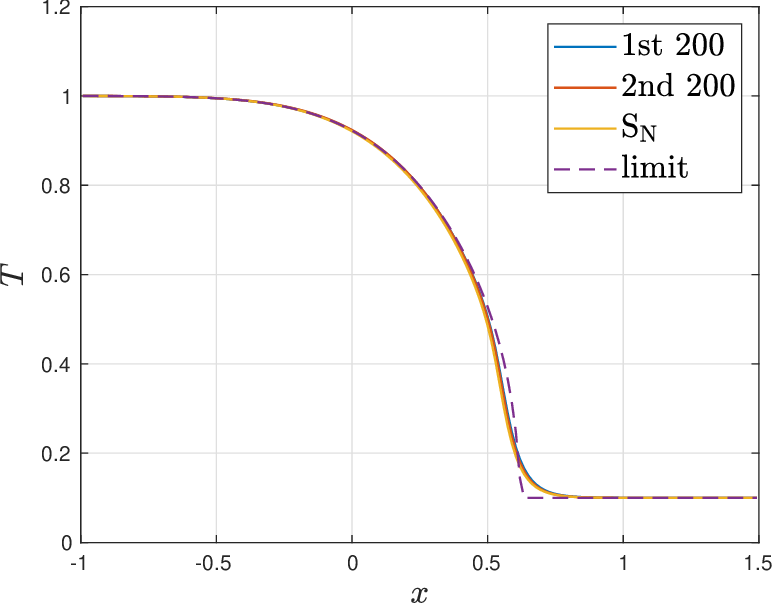}}
  \caption{(2D Riemann problem in Sec. \ref{sec:ex6}) The slice plot of $T_r$ and $T$ for the 2D Riemann problem with $\epsilon = 0.1$ at $t = 0.5$. (a) $T_r$ along $y = x + 0.5$. (b)  $T_r$ along $y = x - 0.5$. (c) $T$ along $y = x + 0.5$. (d) $T$ along $y = x - 0.5$.}
    \label{fig:Riemann_cut_1e1}
\end{figure}

\begin{figure}[!hptb]
  \centering
  \subfloat[$T_r (\epsilon = 10^{-6})$, 1st order]{
    \label{fig:Riemann_Tr_Pn_ep_1e6_1order}
    \includegraphics[width=0.23\textwidth]{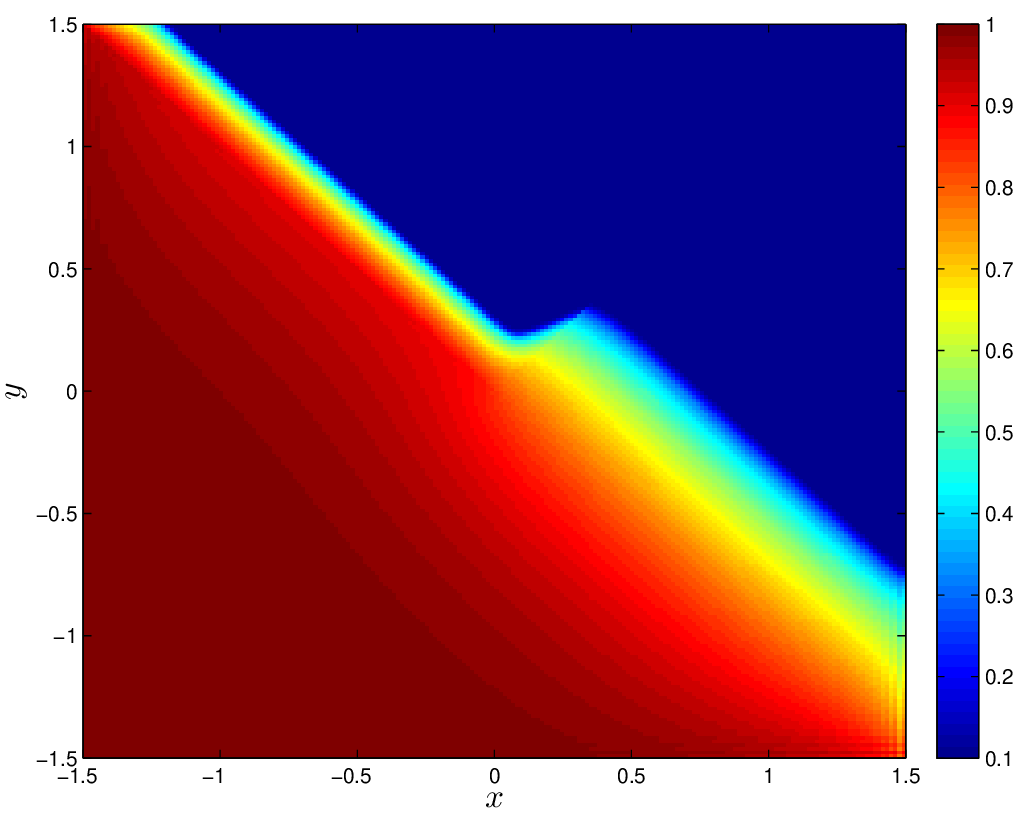}}
    \hfill 
  \subfloat[$T_r (\epsilon = 10^{-6})$, 2nd order]{
    \label{fig:Riemann_Tr_Pn_ep_1e6_2order}
    \includegraphics[width=0.23\textwidth]{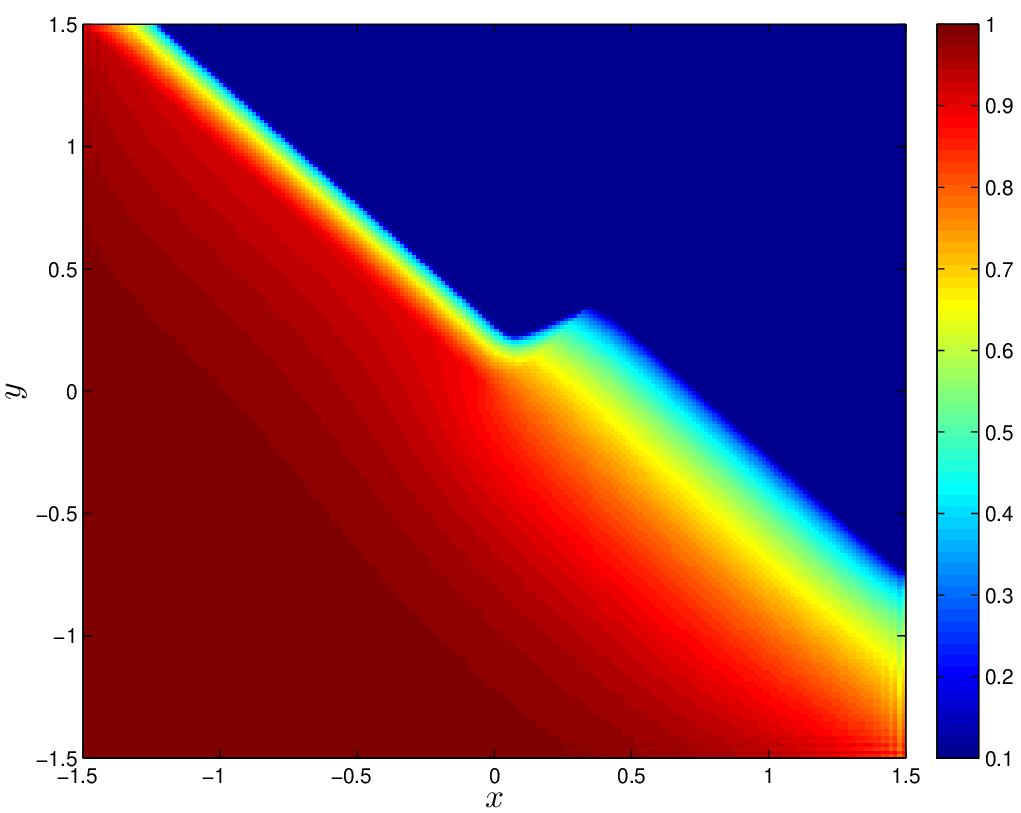}}
  \hfill 
   \subfloat[$T (\epsilon = 10^{-6})$, 1st order]{
    \label{fig:Riemann_T_Pn_ep_1e6_1order}
    \includegraphics[width=0.23\textwidth]{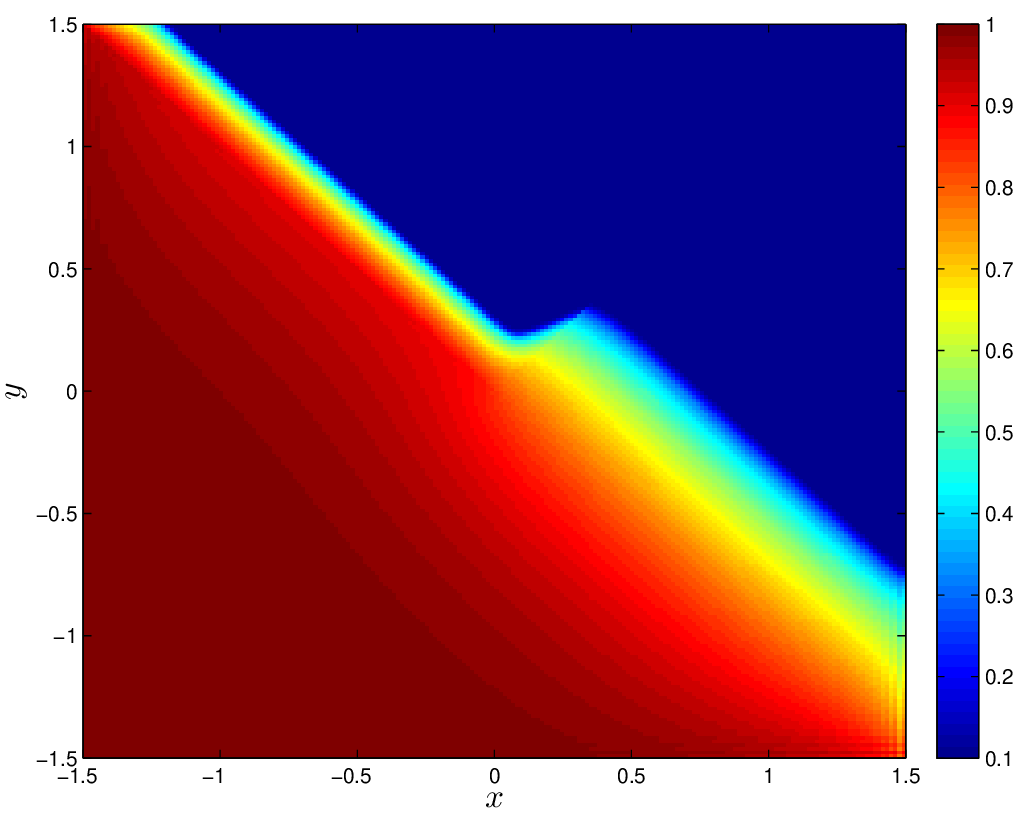}}
    \hfill 
  \subfloat[$T (\epsilon = 10^{-6})$, 2nd order]{
    \label{fig:Riemann_T_Pn_ep_1e6_2order}
    \includegraphics[width=0.23\textwidth]{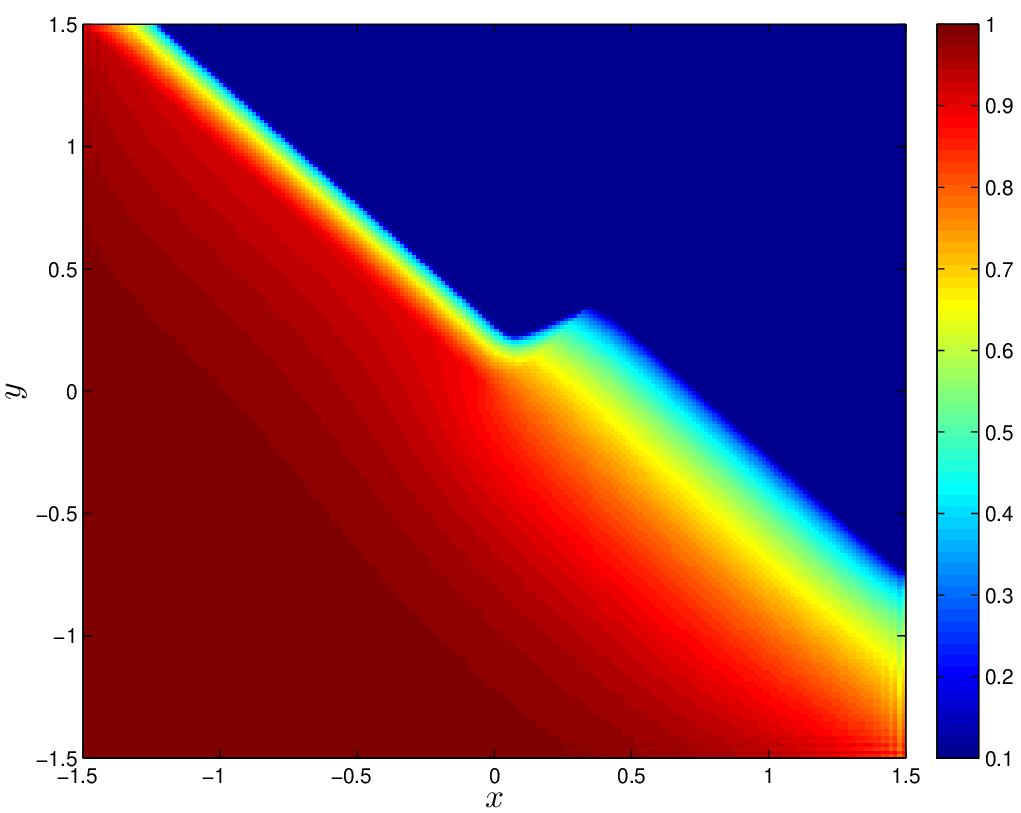}}
  \caption{(2D Riemann problem in Sec. \ref{sec:ex6}) The numerical solution of $T_r$ and $T$ for the 2D Riemann problem with $\epsilon = 10^{-6}$ at $t = 0.5$. (a) Contour plot of $T_r$ by 1st-order scheme. (b) Contour plot of $T_r$ by 2nd-order scheme. (c) Contour plot of $T$ by 1st-order scheme. (d) Contour plot of $T$ by 2nd-order scheme.}
    \label{fig:Riemann_ep_1e6}
\end{figure}

\begin{figure}[!htbp]
  \centering
  \subfloat[$T_r (\epsilon = 10^{-6}, y = x + 0.5)$]{
    \label{fig:Riemann_Tr_Pn_plus_ep_1e6_line3}
    \includegraphics[width=0.23\textwidth]{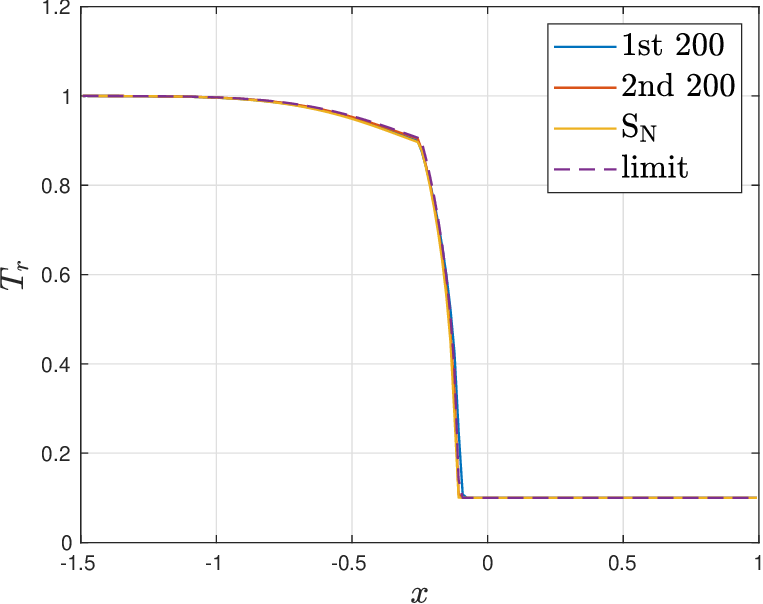}}
   \hfill    
    \subfloat[$T_r (\epsilon = 10^{-6}, y = x - 0.5)$]{
    \label{fig:Riemann_Tr_Pn_minus_ep_1e6_line3}
    \includegraphics[width=0.23\textwidth]{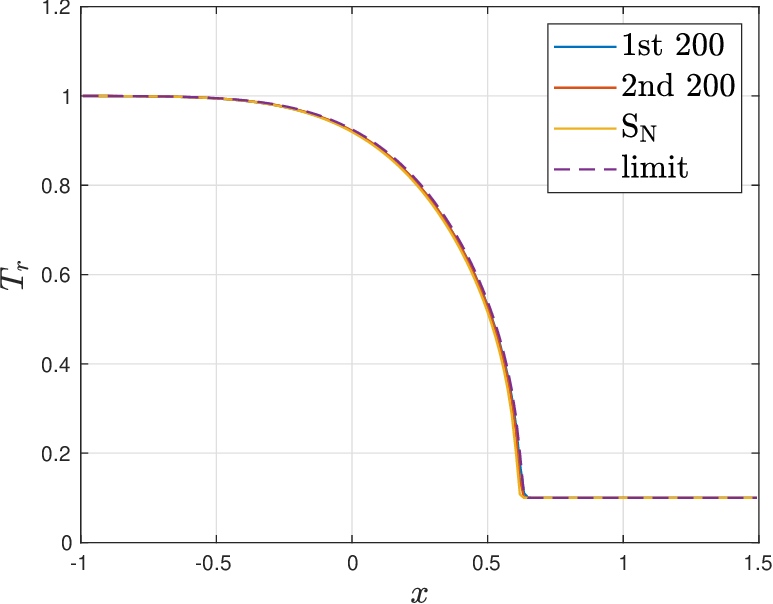}}
    \hfill
     \subfloat[$T (\epsilon = 10^{-6}, y = x + 0.5)$]{
    \label{fig:Riemann_T_Pn_plus_ep_1e6_line3}
    \includegraphics[width=0.23\textwidth]{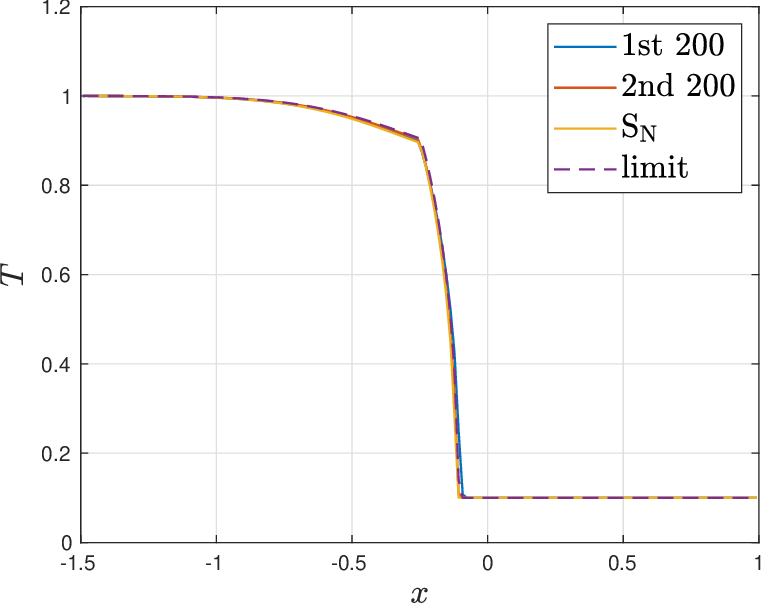}}
   \hfill 
    \subfloat[$T (\epsilon = 10^{-6}, y = x - 0.5)$]{
    \label{fig:Riemann_T_Pn_minus_ep_1e6_line3}
    \includegraphics[width=0.23\textwidth]{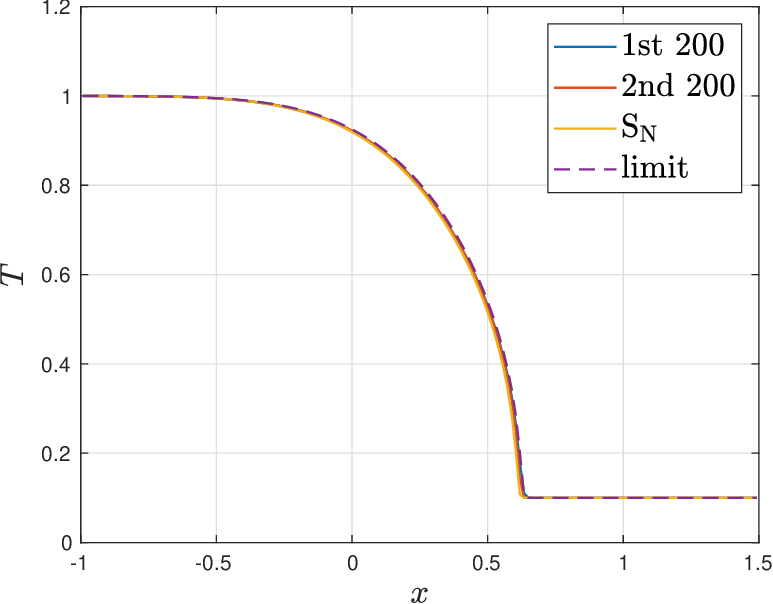}}
   \caption{(2D Riemann problem in Sec. \ref{sec:ex6}) The slice plot of $T_r$ and $T$ for the 2D Riemann problem with $\epsilon = 10^{-6}$ at $t = 0.5$. (a) $T_r$ along $y = x + 0.5$. (b)  $T_r$ along $y = x - 0.5$. (c) $T$ along $y = x + 0.5$. (d) $T$ along $y = x - 0.5$.}
    \label{fig:Riemann_cut_1e6}
\end{figure}

To verify the AP properties of IMEX-IM, we reduce $\epsilon$ to $0.1$, and the expansion order of the $P_N$ method is set to $M = 7$ and the mesh size is $200 \times 200$, with all other parameters the same as $\epsilon = 1$. Fig. \ref{fig:Riemann_ep_1e1} presents the numerical solution of $T_r$ and $T$ at $t = 0.5$, where the behavior of $T_r$ is similar to that of $\epsilon = 1$, while the behavior of $T$ is quite different, as $T$ is high in the gray area, and low in the white area. The numerical solution of $T_r$ and $T$ along $y = x + 0.5$ and $y = x - 0.5$ is provided in Fig. \ref{fig:Riemann_cut_1e1} with the reference solution obtained by the $S_N$ method and the numerical solution of the diffusion limit equation \eqref{eq:limit}, which shows that the resolution of the numerical solution obtained by the second-order scheme is much better than that of the first-order scheme, and matches well with the reference solution. Moreover, there still exists some distance between the numerical solution of RTE and that of the diffusion limit equation when $\epsilon = 0.1$. Then we reduce $\epsilon$ to $10^{-6}$. With the same settings as $\epsilon = 0.1$, the numerical solution of $T_r$ and $T$ at $t = 0.5$ is plotted in Fig. \ref{fig:Riemann_ep_1e6} with the numerical solution along $y = x + 0.5$ and $y = x - 0.5$ presented in Fig. \ref{fig:Riemann_cut_1e6} as well as the reference solution obtained by the $S_N$ method and the numerical solution of the diffusion limit equation \eqref{eq:limit}. In this case, the resolution of the numerical solution obtained by the first- and the second-order schemes is almost the same, and they are well correlated with the reference solution by the $S_N$ method, and that of the diffusion limit equation. This also indicates the high efficiency of IMEX-IM.

\begin{figure}[!htbp]
  \centering
  \subfloat[$T_r (\epsilon = 1, y = x + 0.5)$]{
    \label{fig:Riemann_Tr_Pn_plus_ep_1_1D}
    \includegraphics[width=0.23\textwidth]{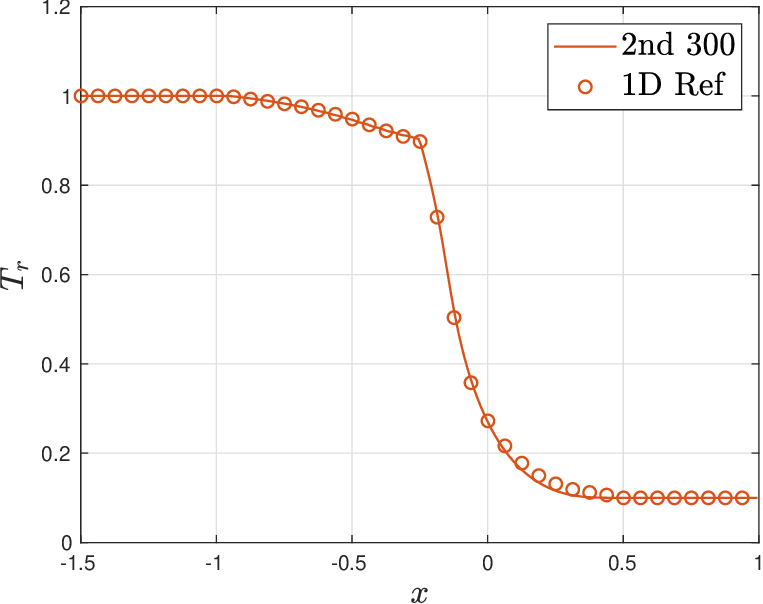}}
    \hfill
    \subfloat[$T_r (\epsilon = 1, y = x - 0.5)$]{
    \label{fig:Riemann_Tr_Pn_minus_ep_1_1D}
    \includegraphics[width=0.23\textwidth]{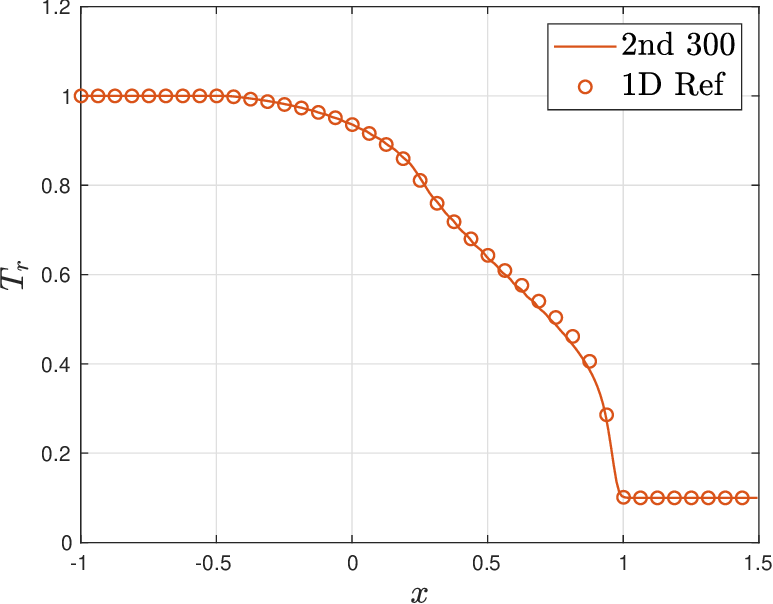}}
    \hfill 
    \subfloat[$T (\epsilon = 1, y = x + 0.5)$]{
    \label{fig:Riemann_T_Pn_plus_ep_1_1D}
    \includegraphics[width=0.23\textwidth]{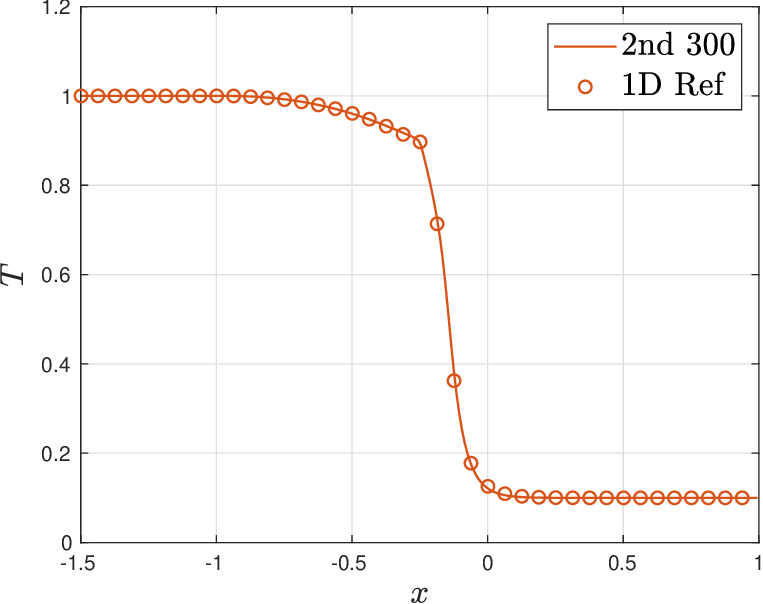}}
    \hfill 
   \subfloat[$T (\epsilon = 1, y = x - 0.5)$]{
    \label{fig:Riemann_T_Pn_minus_ep_1_1D}
    \includegraphics[width=0.23\textwidth]{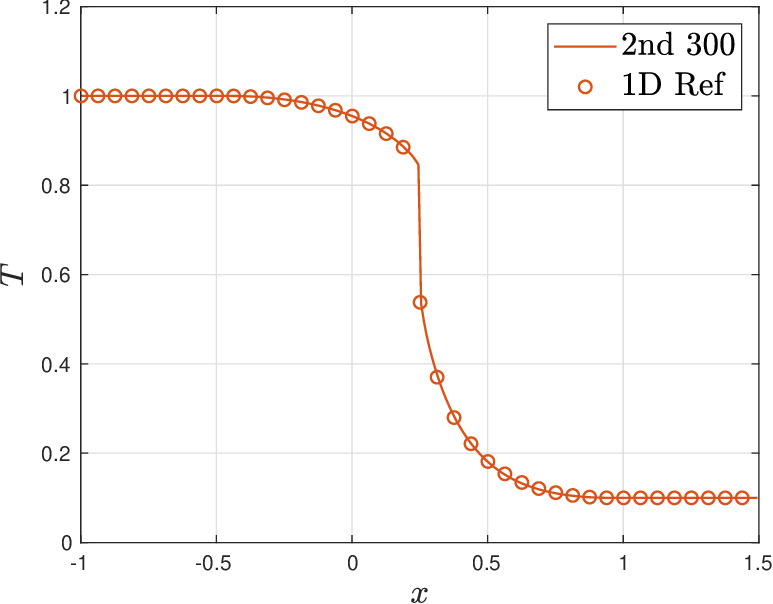}} \\
    \subfloat[$T_r (\epsilon = 0.1, y = x + 0.5)$]{
    \label{fig:Riemann_Tr_Pn_plus_ep_1e1_1D}
    \includegraphics[width=0.23\textwidth]{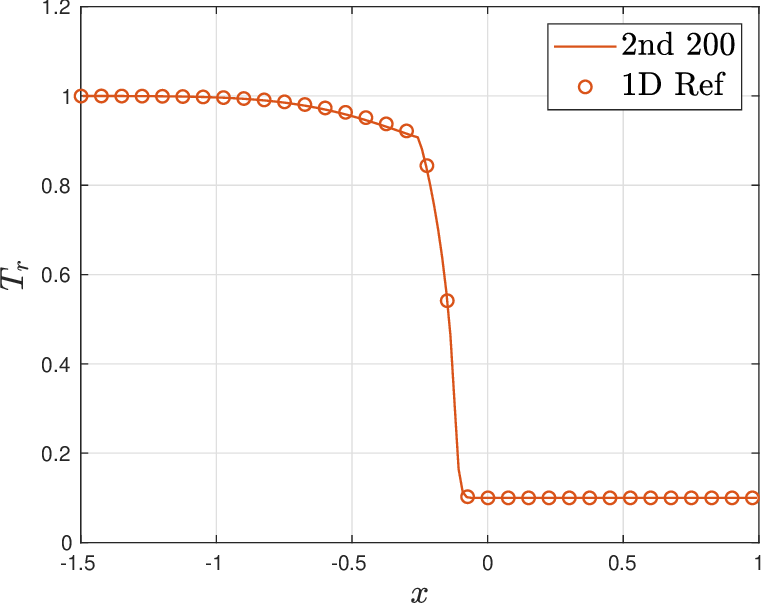}}
    \hfill    
    \subfloat[$T_r (\epsilon = 0.1, y = x - 0.5)$]{
    \label{fig:Riemann_Tr_Pn_minus_ep_1e1_1D}
    \includegraphics[width=0.23\textwidth]{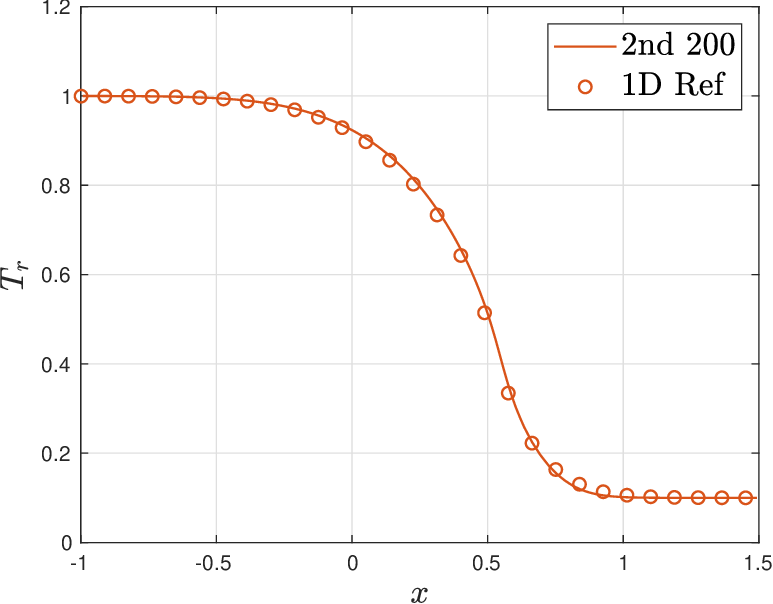}}
    \hfill 
     \subfloat[$T (\epsilon = 0.1, y = x + 0.5)$]{
    \label{fig:Riemann_T_Pn_plus_ep_1e1_1D}
    \includegraphics[width=0.23\textwidth]{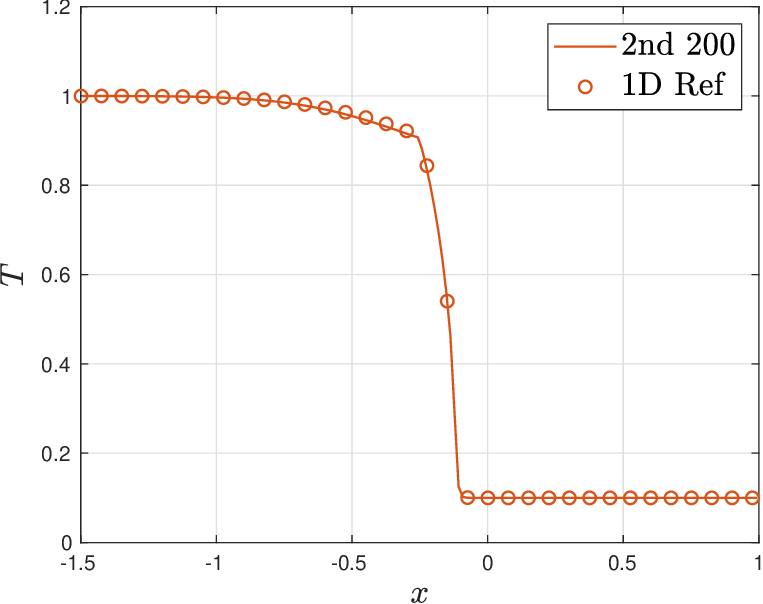}}
    \hfill
  \subfloat[$T (\epsilon = 0.1, y = x - 0.5)$]{
    \label{fig:Riemann_T_Pn_minus_ep_1e1_1D}
    \includegraphics[width=0.23\textwidth]{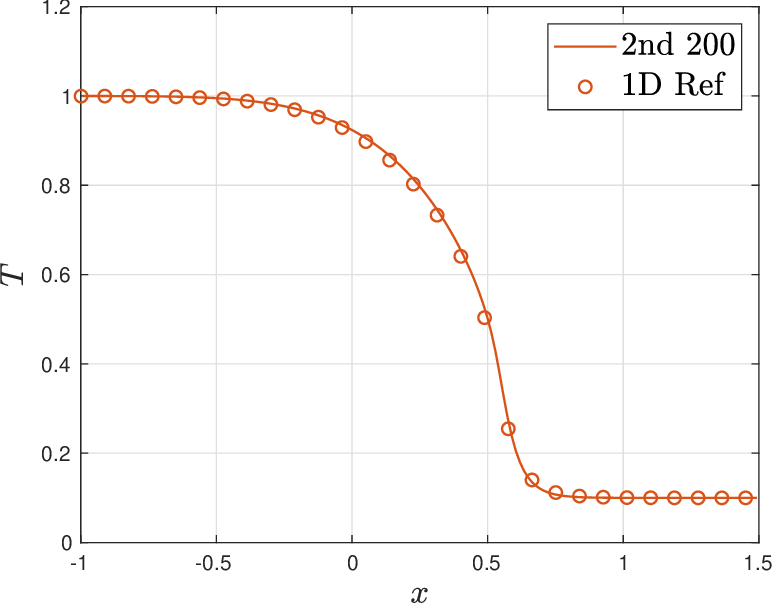}} \\
    \subfloat[$T_r (\epsilon = 10^{-6}, y = x + 0.5)$]{
    \label{fig:Riemann_Tr_Pn_plus_ep_1e6_1D}
    \includegraphics[width=0.23\textwidth]{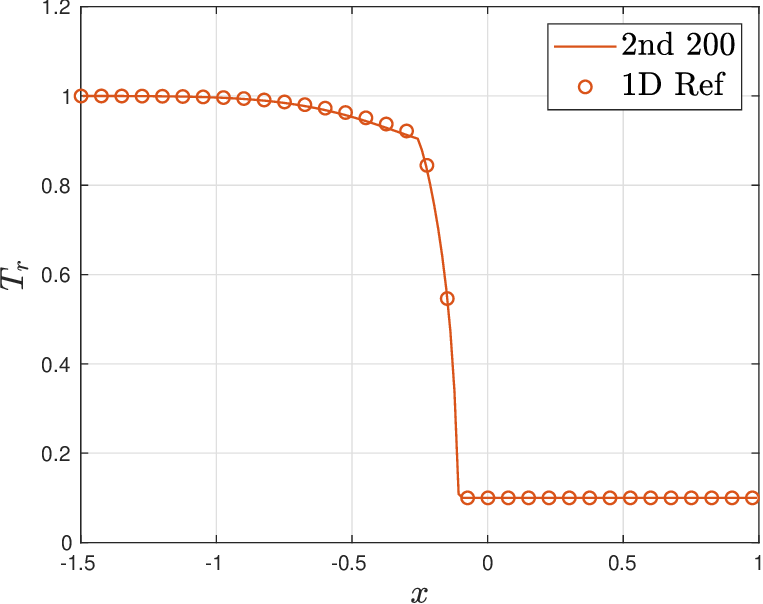}}
    \hfill
    \subfloat[$T_r (\epsilon = 10^{-6}, y = x - 0.5)$]{
    \label{fig:Riemann_Tr_Pn_minus_ep_1e6_1D}
    \includegraphics[width=0.23\textwidth]{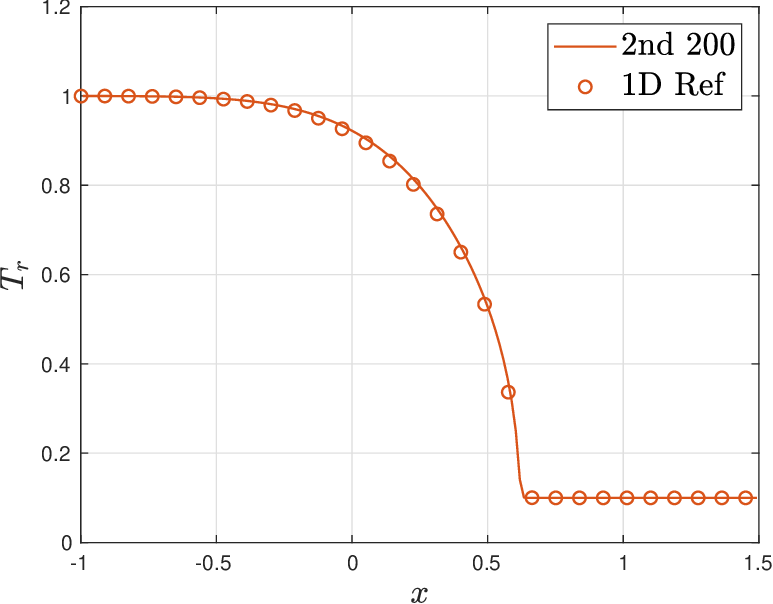}}
    \hfill 
    \subfloat[$T (\epsilon = 10^{-6}, y = x + 0.5)$]{
    \label{fig:Riemann_T_Pn_plus_ep_1e6_1D}
    \includegraphics[width=0.23\textwidth]{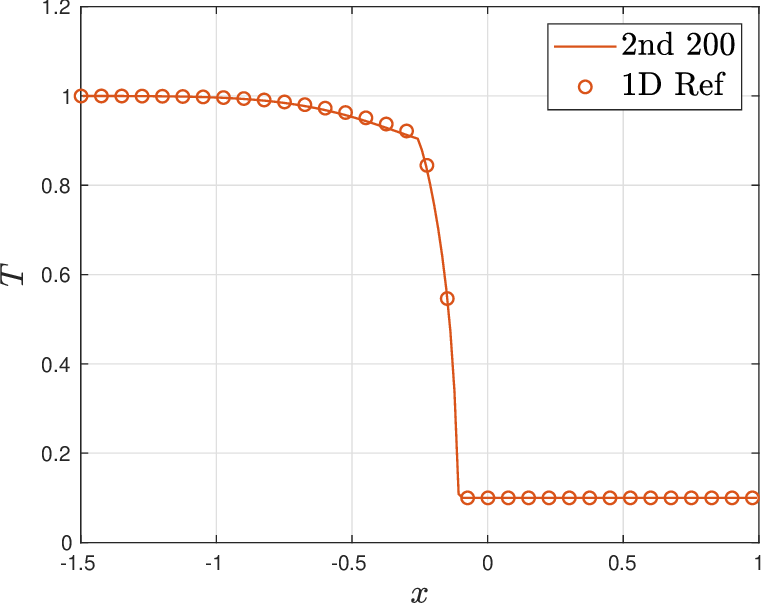}}
    \hfill 
  \subfloat[$T (\epsilon = 10^{-6}, y = x - 0.5)$]{
    \label{fig:Riemann_T_Pn_minus_ep_1e6_1D}
    \includegraphics[width=0.23\textwidth]{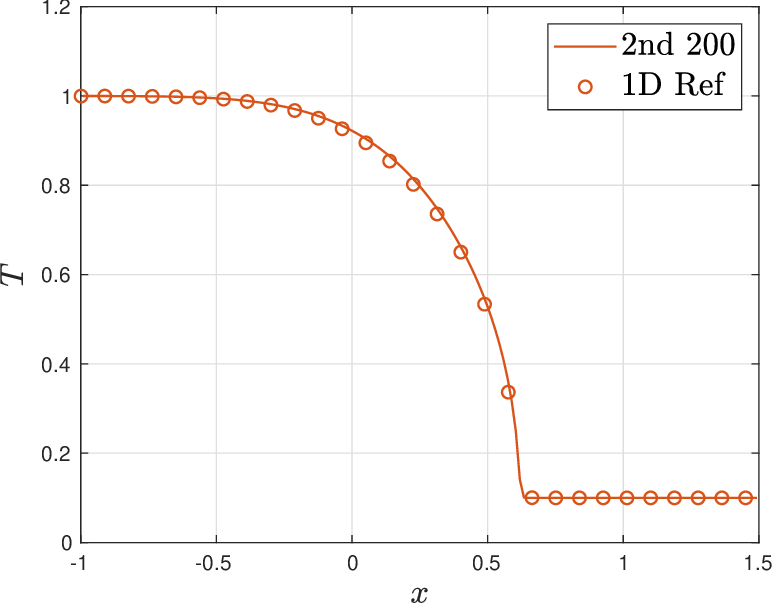}}
  \caption{(2D Riemann problem in Sec. \ref{sec:ex6}) Comparison between the numerical solution to the 2D Riemann problem and the 1D Riemann problem. The first row is that for $\epsilon = 1$, the second is that for $\epsilon = 0.1$, and the third is that for $\epsilon = 10^{-6}$. }
    \label{fig:Riemann_cut_1D}
\end{figure}

 Finally, to further verify the correctness of the numerical solution to this 2D Riemann problem, two 1D problems along the red line $y = x + 0.5$ and blue line $y = x - 0.5$ shown in Fig. \ref{fig:Rie_1D_set} are studied. It is expected that the numerical solution to the two 1D problems is similar to that of the 2D Riemann problem along $y = x + 0.5$ and $y = x - 0.5$. The initial and boundary conditions for the 1D problems are consistent with the 2D Riemann problem. For the 1D problem, the second-order numerical scheme with the grid size $N = 1000$ is utilized. The expansion number for $\epsilon = 1$ is set as $M = 39$ with that for $\epsilon = 0.1$ and $10^{-6}$ set as $M = 7$. Fig. \ref{fig:Riemann_cut_1D} presents the numerical results of $T_r$ and $T$  for the 2D Riemann problem and the 1D Riemann problem with different $\epsilon$. It shows that the two numerical solutions are almost identical for different $\epsilon$, which also validates the efficiency of IMEX-IM.

\subsection{Study of efficiency}
\label{sec:eff}
In this section, the efficiency of this IMEX-IM method is validated by comparing the computational time for some classical benchmark problems with the IMEX numerical method in \cite{IMEX2022} (named IMEX-EX here). For all the numerical tests, the first-order scheme for both methods is utilized, and they are conducted on the CPU model Intel(R) Xeon(R) Gold 5218 CPU@2.30 GHz. 

\begin{figure}[!hptb]
  \centering
   %  \subfloat[IMEX-EX $(\epsilon = 1)$ ]{
   %  \label{fig:linesource_t_05_xian_ep1}
   %  \includegraphics[width=0.30\textwidth]{./IM-EX-compare/linesource_t_05_xian_ep1.eps}}
   %  \hfill
   % \subfloat[IMEX-IX $(\epsilon = 1)$]{
   %  \label{fig:linesource_t_05_yin_ep1}
   %  \includegraphics[width=0.30\textwidth]{./Linesource_fig/linesource_t_05_CFL_04_ep_1_1order_imagesc.eps}}
   %  \hfill
   % \subfloat[$x = 0  (\epsilon = 1)$]{
   %  \label{fig:linesource_t_05_xian_yin_ep1_x}
   %  \includegraphics[width=0.30\textwidth]{./IM-EX-compare/linesource_t_05_xian_yin_ep1_x.eps}}
   %   \hfill
   %    \subfloat[IMEX-EX, $(\epsilon = 10^{-1})$]{
   %  \label{fig:linesource_t_05_xian_ep01}
   %  \includegraphics[width=0.30\textwidth]{./IM-EX-compare/linesource_t_05_xian_ep01.eps}}
   %   \hfill
   % \subfloat[IMEX-IM $(\epsilon = 10^{-1})$]{
   %  \label{fig:linesource_t_05_yin_ep01}
   %  \includegraphics[width=0.30\textwidth]{./Linesource_fig/linesource_t_05_CFL_04_ep_1e1_1order_imagesc.eps}}
   %   \hfill
   % \subfloat[$x = 0 (\epsilon = 10^{-1})$]{
   %  \label{fig:linesource_t_05_xian_yin_ep01}
   %  \includegraphics[width=0.30\textwidth]{./IM-EX-compare/linesource_t_05_xian_yin_ep01_x.eps}}
   %   \hfill
    \subfloat[IMEX-EX]{
    \label{fig:linesource_t_05_xian_ep1e6}
    \includegraphics[width=0.30\textwidth]{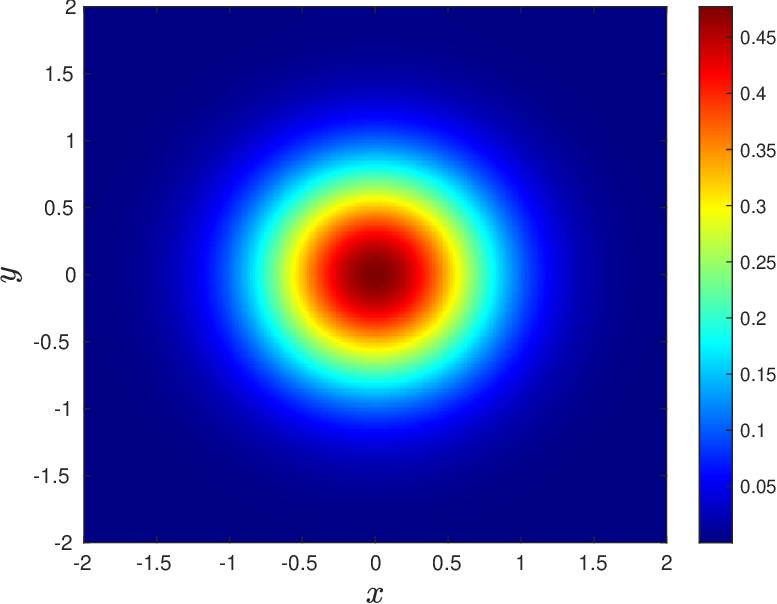}}
    \hfill
      \subfloat[IMEX-IM]{
    \label{fig:linesource_t_05_yin_ep1e6}
    \includegraphics[width=0.30\textwidth]{./images/linesource_t_05_CFL_04_ep_1e6_1order_imagesc.eps}}
  \hfill
   \subfloat[$x = 0$]{
    \label{fig:linesource_t_05_xian_yin_ep1e6_x}
    \includegraphics[width=0.30\textwidth]{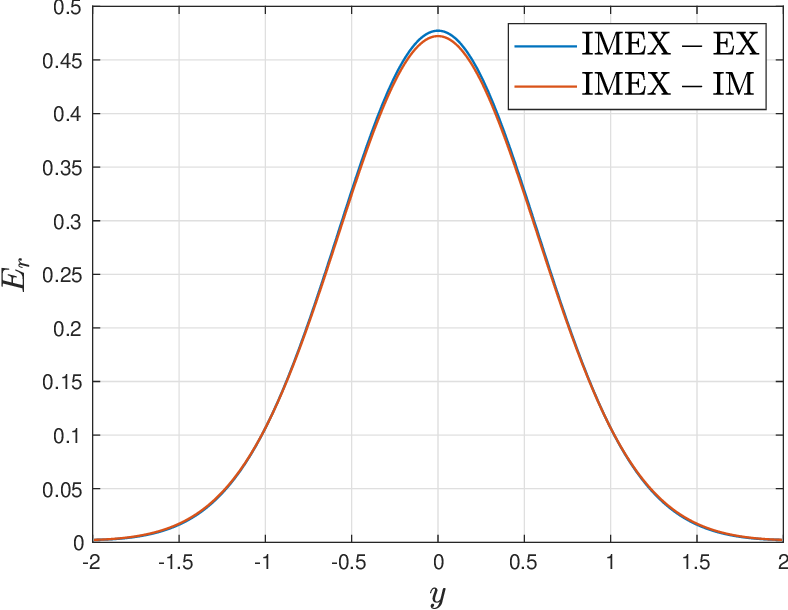}}
  \caption{(Study of efficiency Sec. \ref{sec:eff}) Comparison of numerical solution $E_r$ for the line source problem between IMEX-EX and IMEX-IM for $\epsilon = 10^{-6}$ with $N_x = N_y = 256$. }
\label{fig:linesource_ep_com}
\end{figure}

First, the line source problem in Sec. \ref{sec:ex4} with $\epsilon = 10^{-6}$ is studied with other computational parameters the same as those in Sec. \ref{sec:ex4}. The time step length for IMEX-EX is $\Delta t = 
C\sigma_{\min}\Delta x^2 / c$ with $\sigma_{\min} = 1$ here, while it is set as $\Delta t = C \Delta x / c$ 
for IMEX-IM. The numerical solution $E_r$ for both methods at $t = 0.5$ with mesh size $N_x = N_y = 256$ is plotted in Fig. \ref{fig:linesource_ep_com}, where the two solutions are almost on top of each other, which means that the two methods are of the same resolution. To show the improvement of the efficiency clearly, the study of different mesh sizes is conducted. The total computational time for both methods with mesh size $N_x = N_y = 64, 128$ and $256$ is listed in Tab. \ref{table:compute_time_linesource_ms} as well as the length of the time steps, and the total time steps. For all the mesh sizes, the time step length for IMEX-IM is much larger than that of IMEX-EM. Especially when $N_x = N_y = 256$, the time step is enlarged for more than $60$ times, and the larger the mesh size is, the more increment the efficiency is. 
% \begin{table}[htbp]
% \centering
% \def\arraystretch{1.5}
% \scalebox{0.8}{
% \begin{tabular}{c|c|cccc}
% \multicolumn{6}{c}{{\bf line source problem}}\\
% \hline
% mesh size & $\epsilon$ & method & expansion order & total time  steps & total time (s)\\
% \hline
%      \multirow{6}{*}{$256 \times 256$}&
%      \multirow{2}{*}{$1$} & IMEX-EX & $M = 59$ & 267&  65507\\
%                        &   & IMEX-IM & $M = 59$ & 267& 73384 \\ 
%                        \cline{2-6}
%       &\multirow{2}{*}{$0.1$} & IMEX-EX & $M = 7$ & 800& 1874 \\
%                      &   & IMEX-IM & $M = 7$ & 800& 3313 \\ 
%                        \cline{2-6}
%      &\multirow{2}{*}{$10^{-6}$} & IMEX-EX & $M = 7$ & 5120& 11901  \\
%              &   & IMEX-IM & $M = 7$ & ${\bm{80}}$&  ${\bm{452}}$           
%     \end{tabular}
%     }
%     \caption{(Study of efficiency in Sec. \ref{sec:eff}) Comparison of computation costs between IMEX-EX and IMEX-IM for the line source problem with $\epsilon = 1, 0.1$ and $10^{-6}$ when mesh size $N_x = N_y = 256$.}
% \label{table:compute_time_linesource_ep}
% \end{table}

\begin{table}[htbp]
\centering
\def\arraystretch{1.5}
\scalebox{0.8}{
\begin{tabular}{c|cccc}
\multicolumn{5}{c}{{\bf line source problem}}\\
\hline
mesh size & method & $\Delta t$ & total time steps & total time (s)\\
\hline
    \multirow{2}{*}{$64\times 64$}&
    IMEX-EX & $1.56 \times 10^{-3}$ & 320& 45 \\
    & IMEX-IM & $2.50 \times 10^{-2}$ & 20& 5\\
    \hline
    \multirow{2}{*}{$128\times 128$}&
    IMEX-EX & $3.91 \times 10^{-4}$ & 1280& 692 \\
    & IMEX-IM & $1.25 \times 10^{-2}$ & 40 & 44 \\
    \hline
    \multirow{2}{*}{$256\times 256$}&
    IMEX-EX & $9.77 \times 10^{-5}$ & 5120& 11901\\
    & IMEX-IM & $6.25 \times 10^{-3}$ & 80&  452\\  
\end{tabular}
}
\caption{(Study of efficiency in Sec. \ref{sec:eff}) Comparison of computation costs between IMEX-EX and IMEX-IM for the line source problem with mesh size $N_x = N_y = 64, 128$ and $256$ when $\epsilon = 10^{-6}$.
}
\label{table:compute_time_linesource_ms}
\end{table}

\begin{figure}[!hptb]
  \centering
  %   \subfloat[IMEX-EX$(\epsilon = 1)$ ]{
  %   \label{fig:lattice_t_32_xian_ep1}
  %   \includegraphics[width=0.23\textwidth]{./IM-EX-compare/lattice_t_32_xian_ep1.eps}}
  % \hfill
  %  \subfloat[IMEX-IM$(\epsilon = 1)$]{
  %   \label{fig:lattice_t_32_yin_ep1}
  %   \includegraphics[width=0.23\textwidth]{./IM-EX-compare/lattice_t_32_yin_ep1.eps}}
  %  \hfill
  %  \subfloat[$x = 3.5  (\epsilon = 1)$]{
  %   \label{fig:lattice_t_32_xian_yin_ep1_x_35}
  %   \includegraphics[width=0.23\textwidth]{./IM-EX-compare/lattice_t_32_xian_yin_ep1_x_35.eps}}
  %   \hfill
  %   \subfloat[$y = 3.5  (\epsilon = 1)$]{
  %   \label{figlattice_t_32_xian_yin_ep1_y_35}
  %   \includegraphics[width=0.23\textwidth]{./IM-EX-compare/lattice_t_32_xian_yin_ep1_y_35.eps}}
  %   \hfill
  %     \subfloat[IMEX-EX$(\epsilon = 10^{-1})$]{
  %   \label{fig:lattice_t_05_xian_ep01}
  %   \includegraphics[width=0.23\textwidth]{./IM-EX-compare/lattice_t_05_xian_ep01.eps}}
  % \hfill
  %  \subfloat[IMEX-IM$(\epsilon = 10^{-1})$]{
  %   \label{fig:lattice_t_32_yin_ep01}
  %   \includegraphics[width=0.23\textwidth]{./Lattice_fig/lattice_t_05_CFL_04_ep_1e1_1order_imagesc_N_200.eps}}
  %  \hfill
  %  \subfloat[$x = 3.5 (\epsilon = 10^{-1})$]{
  %   \label{fig:lattice_t_05_xian_yin_ep01_x_35}
  %   \includegraphics[width=0.23\textwidth]{./IM-EX-compare/lattice_t_05_xian_yin_ep01_x_35.eps}}
  %   \hfill
  %   \subfloat[$y = 3.5 (\epsilon = 10^{-1})$]{
  %   \label{fig:lattice_t_05_xian_yin_ep01_y_35}
  %   \includegraphics[width=0.23\textwidth]{./IM-EX-compare/lattice_t_05_xian_yin_ep01_y_35.eps}}
  %    \hfill
      \subfloat[IMEX-EX]{
    \label{fig:lattice_t_05_ex_ep1e6}
    \includegraphics[width=0.23\textwidth]{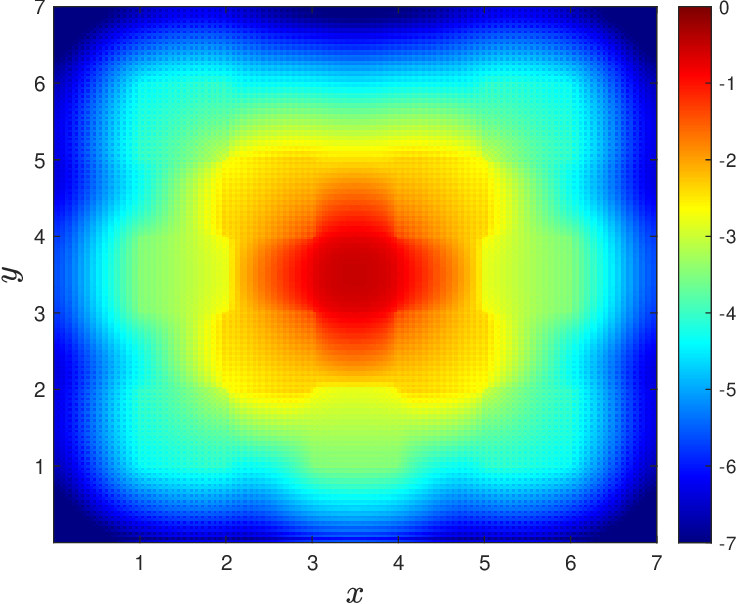}}
  \hfill
   \subfloat[IMEX-IM]{
    \label{fig:lattice_t_32_im_ep1e6}
    \includegraphics[width=0.23\textwidth]{./images/lattice_t_05_CFL_04_ep_1e6_1order_imagesc_N_200.eps}}
   \hfill
   \subfloat[$x = 3.5$]{
    \label{fig:lattice_t_05_ex_im_ep1e6_x_35}
    \includegraphics[width=0.23\textwidth]{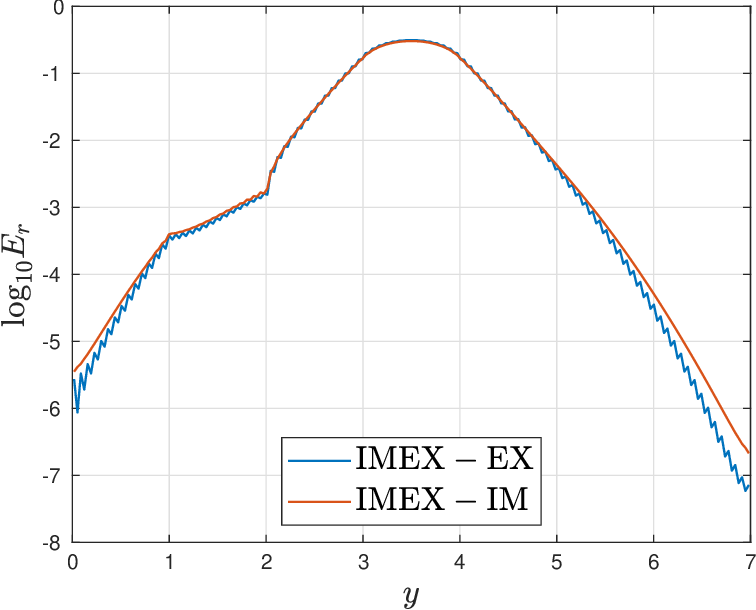}}
    \hfill
    \subfloat[$y = 3.5$]{
    \label{fig:lattice_t_05_ex_im_ep1e6_y_35}
    \includegraphics[width=0.23\textwidth]{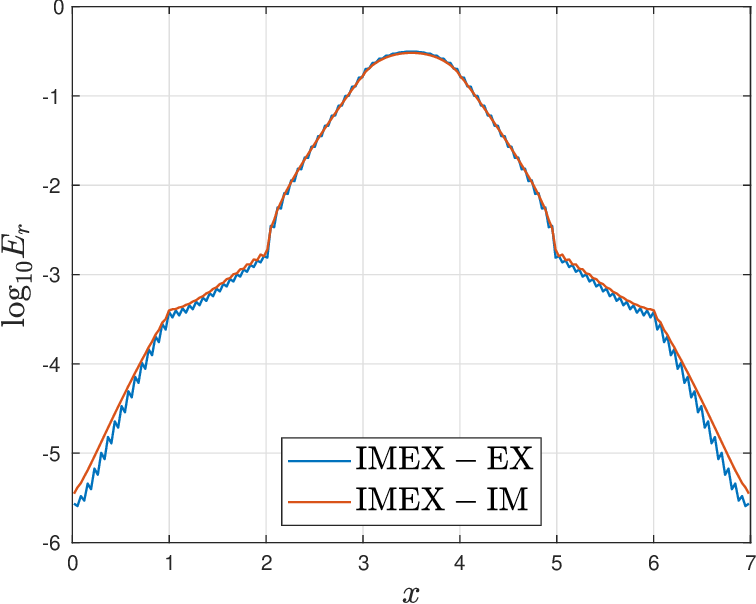}}
  \caption{(Study of efficiency in Sec. \ref{sec:eff}) Comparison of numerical solution $\log_{10} E_r$ for the lattice problem between IMEX-EX and IMEX-IM for  $\epsilon = 10^{-6}$ with $N_x = N_y = 200$.}
\label{fig:lattice_ep_com}
\end{figure}

Next, the lattice problem in Sec. \ref{sec:ex5} is studied in the diffusion regime with $\epsilon = 10^{-6}$. The same time step length as that for the line source problem is chosen with $\sigma_{\min} = 0.1$ in this problem. The numerical solution of $\log_{10} E_r$ obtained by these two methods with $N_x = N_y = 200$  is illustrated in Fig. \ref{fig:lattice_ep_com}. Fig. \ref{fig:lattice_t_05_ex_im_ep1e6_x_35} and \ref{fig:lattice_t_05_ex_im_ep1e6_y_35} show that there exist small oscillations for IMEX-EX near the boundary while that of IMEX-IM keeps smooth. Tab. \ref{table:compute_time_lattice_ms} presents the total computational time as well as the time step length and the total time steps for both methods with mesh size $N_x = N_y = 100$ and $200$. It indicates that for both mesh sizes, the time step length for IMEX-IM is much larger than IMEX-EX. Especially, when $N_x = N_y = 200$, it is enlarged approximately $300$ times, and the total computational time is shortened almost $100$ times accordingly, which all reveals the high efficiency of IMEX-IM compared to IMEX-EX.

\begin{table}[htbp]
\centering
\def\arraystretch{1.5}
\scalebox{0.8}{
\begin{tabular}{c|cccc}
\multicolumn{5}{c}{{\bf lattice problem}}\\
\hline
mesh size & method & $\Delta t$ & total time steps & total time (s)\\
\hline
    \multirow{2}{*}{$100\times 100$}&
    IMEX-EX & $1.96 \times 10^{-4}$ & 2552& 718 \\
    & IMEX-IM & $2.80 \times 10^{-2}$ & 18& 12\\
    \hline
    \multirow{2}{*}{$200\times 200$}&
    IMEX-EX & $4.90 \times 10^{-5}$ & 10204& 12037 \\
    & IMEX-IM & $1.40 \times 10^{-2}$ & 36 & 122 \\
\end{tabular}
}
\caption{(Study of efficiency in Sec. \ref{sec:eff}) Comparison of computation costs between IMEX-EX and IMEX-IM for the lattice problem with mesh size $N_x = N_y = 100$ and $200$ when $\epsilon = 10^{-6}$.
}
\label{table:compute_time_lattice_ms}
\end{table}

%% file: article_conclusion.tex
\section{Conclusion}
\label{sec:conclusion}
We proposed an IMEX numerical scheme for the radiative transfer equations in the framework of the $P_N$ method, where the first- and second-order schemes in the temporal discretization for the linear RTE model and the gray model of RTE are discussed. Its AP property and the numerical stability analyzed by the Fourier analysis are also presented. Several numerical examples are studied to verify the AP property numerically and validate the efficiency of this numerical scheme. 

\section*{Acknowledgments}
This work of Yanli Wang is partially supported by the National Natural Science Foundation of China (Grant No. 12171026, U2230402, and 12031013), and the Foundation of the President of China Academy of Engineering Physics (YZJJZQ2022017). Tao Xiong is partially supported by NSFC grant No. 92270112, and NSF of Fujian Province grant No. 2023J02003.
Weiming Li is partially supported by NSFC grant No. 12471365. Juan Cheng is partially supported by NSFC grant No. 12031001.